\documentclass{amsart}

\usepackage[all,knot]{xy}
\xyoption{arc}
\usepackage{latexsym}
\usepackage{amssymb}
\usepackage{amsmath}
\usepackage{amsthm}
\usepackage{amscd}
\usepackage{fancyhdr}
\usepackage{fullpage}
\usepackage{mathrsfs}
\input cyracc.def

\newfam\cyifam
  \font\tencyi=wncyi10
  \font\sevencyi=wncyi7
  \font\fivecyi=wncyi5
  
  \textfont\cyifam=\tencyi \scriptfont\cyifam=\sevencyi
    \scriptscriptfont\cyifam=\fivecyi

\def\id{{\mbox{1 \hskip -7pt 1}}}
\def\ids{{1\hspace{-0.97mm} 1}}
\newcommand{\sgn}{{\mathit s  \mathit g\mathit  n}}
 \newcommand{\lon}{\longrightarrow}
 \newcommand{\bu}{\bullet}
 
 \newcommand{\rar}{\rightarrow}

\newcommand{\p}{{\partial}}
\newcommand{\Id}{{\mathrm I\mathrm d}}

\newcommand{\Mor}{{\mathrm{Mor}}}

\newcommand{\GCor}{\mbox{\sf GC}^{\text{or}}}
\newcommand{\fGC}{\mbox{\sf fGC}}
\newcommand{\fcGC}{\mbox{\sf fcGC}}
\newcommand{\GC}{\mbox{\sf GC}}

\newcommand{\dfGC}{\mbox{\sf dfGC}}
\newcommand{\grt}{\fg\fr\ft}

\newcommand{\B}{{\mathbb B}}
\newcommand{\Q}{{\mathbb Q}}
 \newcommand{\Z}{{\mathbb Z}}
 \newcommand{\bS}{{\mathbb S}}
 \renewcommand{\P}{{\mathbb P}}
 \newcommand{\C}{{\mathbb C}}
 \newcommand{\R}{{\mathbb R}}
 \newcommand{\N}{{\mathbb N}}
 \newcommand{\K}{{\mathbb K}}

 \newcommand{\bbH}{{\mathbb H}}

 \newcommand{\Ker}{{\mathsf{Ker}}}
 \newcommand{\Img}{{\mathsf I\mathsf m}\, }
 \newcommand{\Hom}{{\mathrm H\mathrm o\mathrm m}}

\newcommand{\ccdot}{\hspace{-1mm}\cdot\hspace{-1mm}}

\newcommand{\Def}{{\mathsf{Def}}}
 \newcommand{\LB}{\mathcal{L}\mathit{ieb}}

\newcommand{\wHoLBcd}{\widehat{\mathcal{H}\mathit{olieb}}_{c,d}}

\newcommand{\HoLB}{\mathcal{H}\mathit{olieb}}

 \newcommand{\ot}{\otimes}
 
 \newcommand{\wh}{\widehat}

\newcommand{\sC}{{\mathsf C}}
\newcommand{\sD}{{\mathsf D}}
\newcommand{\sE}{{\mathsf E}}

\newcommand{\sG}{{\mathsf G}}

\newcommand{\sS}{{\mathsf S}}
\newcommand{\sfd}{{\mathsf d}}
\newcommand{\sfb}{{\mathsf b}}

\newcommand{\sfo}{{\mathsf o}}

 \newcommand{\Beq}{\begin{equation}}
 \newcommand{\Eeq}{\end{equation}}
 \newcommand{\Beqr}{\begin{eqnarray}}
 \newcommand{\Eeqr}{\end{eqnarray}}
 \newcommand{\Beqrn}{\begin{eqnarray*}}
 \newcommand{\Eeqrn}{\end{eqnarray*}}
 \newcommand{\Ba}{\begin{array}}
 \newcommand{\Ea}{\end{array}}
 \newcommand{\Bi}{\begin{itemize}}
 \newcommand{\Ei}{\end{itemize}}
 \newcommand{\Bc}{\begin{center}}
 \newcommand{\Ec}{\end{center}}
  
  \newcommand{\fb}{{\mathfrak b}}
 \newcommand{\fg}{{\mathfrak g}}
 \newcommand{\fh}{{\mathfrak h}}
 \newcommand{\fm}{{\mathfrak m}}

\newcommand{\fr}{{\mathfrak r}}

\newcommand{\ft}{{\mathfrak t}}
\newcommand{\fp}{{\mathfrak p}}

\newcommand{\ffree}{{\mathfrak{lie}}}
\newcommand{\fass}{{\mathfrak{ass}}}

 \newcommand{\f}{{\mathcal O}}
 \newcommand{\cA}{{\mathcal A}}
 \newcommand{\cB}{{\mathcal B}}
 \newcommand{\cC}{{\mathcal C}}
 \newcommand{\caD}{{\mathcal D}}
 \newcommand{\cE}{{\mathcal E}}
 \newcommand{\cF}{{\mathcal F}}
 \newcommand{\cG}{{\mathcal G}}
 
 \newcommand{\cI}{{\mathcal I}}

 \newcommand{\caL}{{\mathcal L}}
 \newcommand{\cM}{{\mathcal M}}
 
 \newcommand{\cP}{{\mathcal P}}
 
 \newcommand{\cS}{{\mathcal S}}
 \newcommand{\cT}{{\mathcal T}}

 \newcommand{\al}{\alpha}
 \newcommand{\be}{\beta}
 \newcommand{\ga}{\gamma}
 
 \newcommand{\Ga}{\Gamma}

 \newcommand{\var}{\varepsilon}
 \newcommand{\la}{\lambda}

 \newcommand{\sip}{\smallskip}
 \newcommand{\bip}{\bigskip}
 \newcommand{\mip}{\vspace{2.5mm}}

\newtheorem{theorem}{Theorem}[section]

\theoremstyle{definition}
\newtheorem{definition}[theorem]{Definition}

\theoremstyle{remark}

\numberwithin{equation}{section}

\usepackage{tikz}
\usetikzlibrary{matrix}
\usetikzlibrary{shapes}
\usetikzlibrary{arrows}
\usetikzlibrary{calc,3d}
\usetikzlibrary{decorations,decorations.pathmorphing,decorations.pathreplacing,decorations.markings}
\usetikzlibrary{through}

\tikzset{ext/.style={circle, draw,inner sep=1pt},int/.style={circle,draw,fill,inner sep=1.4pt},nil/.style={inner sep=1pt}}
\tikzset{cy/.style={circle,draw,fill,inner sep=2pt},scy/.style={circle,draw,inner sep=2pt},scyx/.style={draw,cross out,inner sep=2pt},scyt/.style={draw,regular polygon,regular polygon sides=3,inner sep=0.95pt}}
\tikzset{exte/.style={circle, draw,inner sep=3pt},inte/.style={circle,draw,fill,inner sep=3pt}}
\tikzset{diagram/.style={matrix of math nodes, row sep=3em, column sep=2.5em, text height=1.5ex, text depth=0.25ex}}
\tikzset{diagram2/.style={matrix of math nodes, row sep=0.5em, column sep=0.5em, text height=1.5ex, text depth=0.25ex}}

 \tikzset{
  rightblack/.style={
    decoration={markings,mark=at position .8 with {\arrow[scale=1.2,black]{latex}}},
    postaction={decorate},
    shorten >=0.4pt}}
\tikzset{
  leftblack/.style={
    decoration={markings,mark=at position .55 with {\arrowreversed[scale=1.2,black]{latex}}},
    postaction={decorate},
    shorten >=0.4pt}}

\begin{document}

\title{Grothendieck-Teichm\"uller group,\\ operads and graph complexes: a survey}

\author{Sergei Merkulov}
 \address{Sergei~Merkulov:  Mathematics Research Unit, University of Luxembourg,  Grand Duchy of Luxembourg}
\email{sergei.merkulov@uni.lu}

\subjclass[2000]{Primary }
\date{}

\begin{abstract}
The paper offers a more or less self-contained introduction into the theory of the Grothendieck-Teichm\"uller group and  Drinfeld associators using the theory of operads and graph complexes.
\end{abstract}

\maketitle


 \sloppy

 \newenvironment{proo}{\begin{trivlist} \item{\sc {Proof.}}}
  {\hfill $\square$ \end{trivlist}}

\long\def\symbolfootnote[#1]#2{\begingroup%
\def\thefootnote{\fnsymbol{footnote}}\footnote[#1]{#2}\endgroup}

{\large
\section{\bf Introduction}
}
The Grothendieck-Teichm\"uller group $GT$ (and its graded version $GRT$) is one of the
most interesting and mysterious objects in modern mathematics. This group together with its principal homogeneous space --- the set of V.\ Drinfeld's associators ---
plays a central
role
in many seemingly unrelated areas of mathematics: V.\ Drinfeld pioneered its applications of in the number theory and the theory of quasi-Hopf algebras \cite{Dr}; P.\ Etingof and D.\ Kazhdan \cite{EK,ES} used it to solve the Drinfeld quantization
conjecture for Lie bialgebras; the formality theories of M.\ Kontsevich \cite{Ko1,Ko2,Ko3} and D.\ Tamarkin \cite{Ta1,Ta2} unravels
the role of the  group $GRT$ in deformation quantizations of Poisson structures; A.\ Alekseev and C.\ Torossian applied it to the solution of the
Kashiwara-Vergne
problem in the Lie theory \cite{AT}; the authors of \cite{MW3} found its interpretation as a symmetry group of the involutive Lie bialgebra properad playing a key role in the string topology;  A.\ Alekseev, N.\ Kawazumi, Y.\ Kuno and F.\ Naef \cite{AKKN1,AKKN2} and G.\ Massuyeau \cite{Mas} proved its importance in the Goldman-Turaev theory of spaces
of free homotopy loops in Rieman surfaces of genus $g$ with $n$ punctures; this group --- more precisely its Lie algebra $\grt$ --- plays a central role in the recent spectacular advances of M.\ Chan, S.\ Galatius and S.\ Payne \cite{CGP} in the theory of the cohomology groups of moduli spaces $\cM_g$ of genus $g$ algebraic curves. T.\ Willwacher \cite{Wi1,Wi2} established a very important link between $\grt$ and
cohomologies of some graph complexes which found many applications.

\sip

 Our main purpose is to give a more or less short and self-contained introduction into the theory of the Grothendieck-Teichm\"uller group and its applications through its operadic, properadic and graph complexes incarnations, and explain without proofs some of the results mentioned in the previous paragraph.

\sip

The Grothendieck-Teichm\"uller theory first appeared in A.\ Grothendieck's  famous {\em Esquisse d'un programme} (1981) with the purpose to give a new geometric description of the universal
Galois group $Gal(\Q)$ of the field of  rational numbers; A.\ Grothendieck's main tool to approach this problem was the Deligne-Mumford moduli stack $\cM$ of algebraic curves of arbitrary genus with marked points and a very insightful observation that the (outer) automorphism
group $Out(\pi_1^{geom}(\cM))$ of the geometric fundamental group of this stack can be described more or less explicitly as the set of
elements of a free profinite group on two generators satisfying a small number of equations (see \cite{L} and \cite{LS} for more details and references).
A few years later  (1989) V.\ Drinfeld introduced \cite{Dr} $\K$-prounipotent Grothendieck-Teichm\"uller groups
${GT}$ and $GRT$  for an arbitrary field $\K$ of characteristic zero while studying braided quasi-Hopf algebras and their universal deformations. The profinite version $\widehat{GT}$ of the first group contains  $Gal(\Q)$ and is essentially $Out(\pi_1^{geom}(\cM_0))$, where $\cM_0$ is the of moduli spaces of genus
zero curves with marked points.

\sip

V.\ Drinfeld has written down  explicitly \cite{Dr} systems of algebraic equations
defining elements of both prounipotent groups ${GT}$ and ${GRT}$. Thanks to the works of D.\ Bar-Natan
 \cite{BN} the algebro-geometric meanings of these two systems of equations are well-understood by now with the help of the theory of operads ---
${GT}$ is essentially the automorphism group of the topological operad  $\widehat{\cP a\cB}$ of parenthesized braids
 while ${GRT}$ is  the automorphism group of a much simpler graded analogue $\widehat{\cP a\cC \caD}$ of the operad $\widehat{\cP a\cB}$. The set of Drinfeld associators can be identified with the set of isomorphisms $\widehat{\cP a\cB}\rar
  \widehat{\cP a\cC \caD}$ of operads.
This part of the story  explaining $GT$ and $GRT$ as automorphism groups of operads is covered in \S\S 1-6 of this survey; it is based on the notes of a lecture course read by the author at the University of Stockholm in the fall of 2011. B.\ Fresse has written a book \cite{Fr} explaining this approach to
$GT$ and $GRT$ in much more detail.
 In the second part,
\S\S 7-8, we try to explain a very different  T.\ Willwacher's approach to the group $GRT$ via the graph complexes (and compactified configuration spaces), and its applications to the deformation quantization theory of Poisson structures and Lie bialgebras, Lie theory, Goldman-Turaev theory, and the theory of moduli spaces $\cM_{g}$. We explain various graph complexes
incarnations of the Lie algebra $\grt$ of $GRT$; some of them appear quite mysterious at present.

\bip

{\large
\section{\bf Completions of groups, Lie algebras and algebras}
}
\subsection{Completed filtered vector spaces, Lie algebras and algebras}\label{2: subsec on completed foltered}
 A topological ring (in particular, a field) is a ring $R$ which is also a topological space such that both the addition and the multiplication operations define {\em continuous}\, maps
$R \times R \rar R$.

\sip

A topological vector space $V$ is a vector space over a topological field $\K$ (say, $\R$ or $\C$ equipped with their standard topologies or with the discrete topology) which is endowed with a topology making the vector addition
$V \times V \rar V$ and scalar multiplication $\K \times V \rar V$ operations into continuous maps.
One can define similarly  {\em topological $\K$-(co)algebras},
{\em topological Lie algebras over $\K$}, {\em topological Hopf algebras over $\K$}, etc.
\sip

Topological vector spaces, rings and $\K$-algebras are  particular examples of uniform spaces in which it make sense to talk about Cauchy sequences and hence about their completeness.  Any such a space $V$  can be completed in an essentially unique way with the help  of equivalence classes of Cauchy sequences; $V$ sits in its completion $\widehat{V}$  as a dense subset and the induced topology coincides with the original one. The completion of the tensor product $\wh{V}\ot \wh{W}$ is denoted by $\wh{V}\wh{\ot} \wh{W}$ (sometimes we abbreviate the latter notation to  $\wh{V} {\ot} \wh{W}$
when no confusion may arise).

\sip

We shall be interested in this survey in a class of topological spaces, rings and (co)algebras whose topologies are determined by filtrations. Their completions --- often called {\em complete filtered} spaces, rings and (co)algebras --- can be constructed in an explicit way as inverse limits. We shall explain the phenomenon for algebras; the adoption of all the constructions below  to Lie algebras, coalgebras, etc.\ is straightforward (see, e.g., \cite{Ei}).

\sip

Let $A$ be a $\K$-algebra equipped with a descending filtration by ideals,
\Beq\label{2: descending filtration}
A=\fm_0 \supset \fm_1 \supset \fm_2 \supset \ldots
\Eeq
Hence we get a directed system of quotient rings and their morphisms,
$$
\ldots   \lon   A/\fm_{i+1}    \lon A/\fm_i   \lon   \ldots \lon   A/\fm_2 \lon   A/\fm_1 \lon   0.
$$
The associated completed filtered algebra is defined as the  inverse limit,
$$
\wh{A}:=\lim_{\longleftarrow} A/\fm_i:=
\lim_{\longleftarrow} A/\fm_i:=\left\{ a=(a_1,a_2,\ldots)\in \prod_{i=1}^\infty A/\fm_i\ |\ a_j= a_i \bmod \fm_i\ \ \mbox{for all}\ j> i \right\}
$$
The algebra $\wh{A}$  has an induced filtration by ideals,
$$
\wh{\fm}_i:=\left\{  a=(a_1,a_2,\ldots)\in \wh{A}\ |\  a_j= 0\ \ \mbox{for all}\ j\leq  i  \right\}.
$$
It is clear that the quotient algebras $A/\fm_i$ and $\wh{A}/\wh{\fm}_i$ can be identified. The completed tensor product of such completed algebras is defined by
$$
\wh{A}\wh{\ot}\wh{A'}:= \lim_{\longleftarrow\atop i,j} A/\fm_i \ot A'/\fm_i'.
$$
\sip

 The filtered algebra $A$ and and its completion $\wh{A}$ can be made into topological spaces by defining a basis of open neighborhoods of a point $a$ in $A$ or, respectively, in $\wh{A}$ to be
 $$
 \left\{a+ \fm^i\right\}_{i\in \N} \ \ \ \mbox{or, respectively},\ \ \ \
 \left\{a+ \wh{\fm}^i\right\}_{i\in \N}.
 $$
Such a topology on a filtered algebra is called the {\em Krull}\, topology. It is not hard to see that
 $\wh{A}$ is the completion of  $A$ as a topological algebra. Indeed, let $\{a_i\}_{i\geq 1}$ be a Cauchy sequence in $A$ equipped with the Krull topology, that is, a sequence which satisfies the condition:
  for any open neighborhood $U$ of zero in $A$ there is a number $N_U$ such that for any $i,j> N_U$ one has
  $a_i-a_j\in  U$. Equivalently, $\{a_i\}_{i\geq 1}$ is a Cauchy sequence in $A$ if and only if
  for any integer $n$ there exists
 an integer $N_n$ such that
 $$
 a_i-a_j\in \fm_i \ \ \  \mbox{for all}\ \ \ i,j>N_n.
$$
Such a sequence always converges in $\wh{A}$ to the point $a$ whose $n$-th coordinate in $\prod_{n}A/\fm_n$
is equal, by definition, to $a_{N_n}\mod \fm_n$. Reversely, any point in $\widehat{A}$ gives rise to a
Cauchy sequence in $A$.

\mip

We shall work below with filtrations generated by powers $\fm_i:=I^i$ of a fixed ideal $I$ in $A$.
The associated completion $\wh{A}$ is often called the $I$-{\em adic completion}\, of $A$, and the associated topology on $A$ is called the  $I$-{\em adic topology}.

\subsubsection{\bf Examples}\label{2: examples of completed filtered}

\Bi
\item[(i)]
Let $A=\K[x_1,\ldots,x_n]$ be a polynomial algebra over $\K$, and let $I$ be its maximal ideal.
The $I$-adic completion of $A$,
$$
\wh{A}=\K[[x_1,\ldots, x_n]],
$$
is the algebra of formal power series over $\K$.

\mip

\item[(ii)] Let $\fg$ be a positively graded Lie algebra, $\fg=\bigoplus_{i=1}^\infty \fg_i$,
with $\fg_i$ being finite dimensional. It can be equipped with a descending filtration as in (\ref{2: descending filtration}) given by the Lie ideals
$\fm_i:=\bigoplus_{j\geq i}\fg_j$. The completion of $\fg$ with respect to this filtration is then equal to
$$
\wh{\fg}= \prod_{i\geq 1} \fg_i
$$
and is called {\em the degree completion}\, of $\fg$ (similarly one can define the degree completion of positively graded vector spaces, rings, algebras, etc).
\sip

Consider next the universal enveloping algebra $U(\wh{\fg})$ of the completed graded Lie algebra $\wh{\fg}$. It is, by definition, the quotient of the tensor algebra $\ot^\bu \wh{\fg}$ by the ideal $J$
generated by all elements of the form $x\ot y - y\ot x - [x,y]$,  $x,y\in \wh{\fg}$. The tensor algebra inherits a positive gradation from $\wh{\fg}$ in the standard way,
$$
\ot^\bu \wh{\fg}=\bigoplus_{i\in 1}^\infty (\ot^\bu \wh{\fg})^i\ \ \ \ \mbox{with}\ \ \
(\ot^\bu \wh{\fg})^i:= \bigoplus_{i=i_1+\ldots+i_n\atop n\geq 0, i_1,\ldots, i_n\geq 1} \wh{\fg}_{i_1}\ot\ldots \ot \wh{\fg}_{i_n}.
$$
As $[\wh{\fg}_i,\wh{\fg}_j]\subset \wh{\fg}_{i+j}$, the ideal $J$ is homogeneous with respect to this
gradation, i.e.\ it is generated by the homogeneous elements. Therefore  $U(\wh{\fg})$
comes equipped with an induced gradation and we can define its graded completion $\wh{U}(\wh{\fg})$.
\Ei

\subsubsection{\bf Remark} To define the completed universal enveloping algebra $\wh{U}(\wh{\fg})$
we need only a descending filtration of $\fg$ as in (\ref{2: descending filtration}), a weaker structure than the positive gradation. However, in applications below such a filtration always comes from a positive gradation on $\fg$  as in the above Example {\ref{2: examples of completed filtered}}.

\subsubsection{\bf Exercise} Use the Krull topology on  $\wh{U}(\wh{\fg})$ to show that $\wh{U}(\wh{\fg})$ is a completed filtered Hopf algebra (with the coproduct, $\Delta: \wh{U}(\wh{\fg}) \rar \wh{U}(\wh{\fg})\wh{\ot}
    \wh{U}(\wh{\fg})$, taking values in the completed tensor product).

\subsection{Baker-Campbell-Hausdorff formula}
Recall that one can associate a
group $G:=\exp(\fg)$ to a completed positively graded Lie algebra over a field $\K$ of characteristic zero,
$$
\fg=\prod_{k=1}^\infty \fg_k,
$$
 with all graded components
being of finite dimension. The group $G$ coincides with $\fg$ as a set (and the identification map
 $\fg \rar G$ is denoted by $\exp$ while its inverse $G\rar \fg$ by $\log$), and the group multiplication  is
defined by the Campbell-Hausdorff formula:
$$
\exp(x)\cdot \exp(y):= \exp(\mathfrak{bch}(x,y)),
$$
where
\Beq\label{2: ch(c,y) series}
\mathfrak{bch}(x, y) =\log(e^x e^y)= x + y + \frac{1}{2}[x, y] +  \frac{1}{12}[x,[x, y]] -  \frac{1}{12}[y,[y, x]]\ldots,
\Eeq
and $e^x:=\sum_{k=0}^\infty \frac{x^k}{k!}$ and $\log(1+x):=\sum_{k=1}^\infty\frac{(-1)^k x^k}{k}$.
This formal power series satisfies obviously the associativity relation,
$$
\mathfrak{bch}(x, \mathfrak{bch}(y, z)) = \mathfrak{bch}(\mathfrak{bch}(x, y), z)=\log(e^xe^ye^z),
$$
in the ring $\K\langle\langle x,y,z\rangle\rangle$.

\subsubsection{\bf Example}\label{2: ass_n and Lie_n} Let $\ffree_n$ stand for the free Lie algebra over $\K$ generated by $n$ letters $x_1,\ldots, x_n$ and $\widehat{\ffree}_n$ for its degree completion (with degrees of the generators  $x_1,\ldots, x_n$ set to be $1$).
Let $\fass_n:=\K\langle x_1,\ldots,x_n \rangle$ be the free associative algebra
generated by $x_1,\ldots, x_n$ and  $\widehat{\fass}_n=\K\langle\langle x_1,\ldots,x_n \rangle\rangle$ its degree completion. One can make $\fass_n$ and $\widehat{\fass}_n$ into Lie algebras by setting
$[X,Y]=XY-YX$. Then   $\ffree_n\subset \fass_n$ and $\widehat{\ffree}_n\subset \widehat{\fass}_n$ are Lie subalgebras (in fact they are the smallest Lie subalgebras containing all the generators).

\sip

 One can make  $\widehat{\fass}_n$ into a bialgebra (in fact, a Hopf algebra) by defining a coproduct,
$$
\Delta: \widehat{\fass}_n \lon \widehat{\fass}_n \wh{\ot} \widehat{\fass}_n,
$$
by first setting
$$
\Delta(x_i):= 1\ot x_i + x_i\ot 1
$$
and then extending to arbitrary monomials in $\widehat{\fass}_n$ by the rule
$$
\Delta(x_{i_1}x_{i_2}\cdots x_{i_k})=\Delta(x_{i_1})\Delta(x_{i_2})\cdots \Delta(x_{i_k}).
$$
An element $g\in \widehat{\fass}_n$ is called {\em grouplike} if $\Delta(g)=g\ot g$.
An element $p\in \widehat{\fass}_n$ is called {\em primitive} if $\Delta(p)=1\ot p + p\ot 1$.
The following is true:
\begin{itemize}
\item Let $f$ be a formal power series from  $ \widehat{\fass}_n$ whose constant term is zero. Then  $e^f\in \widehat{\fass_n}$ is grouplike if and only if $f$ is primitive.
\item The set of primitive elements in  $\widehat{\fass}_n$ is a Lie subalgebra of $( \widehat{\fass}_n, [\ ,\ ])$ which can be identified with $\widehat{\ffree}_n$.
\item The set of grouplike elements  in  $\widehat{\fass}_n$ forms a group with respect to the multiplication in $\widehat{\fass}_n$ which can be identified with $\exp(\widehat{\ffree}_n)$.
\end{itemize}
As a Hopf algebra, $\widehat{\fass}_n$ is precisely the completed universal enveloping algebra of $\widehat{\ffree}_n$.

\bip

\subsection{Prounipotent completions} A {\em linear algebraic group}\, over a field $\K$ is, by definition, an algebraic group that is isomorphic to an algebraic subgroup of the group of invertible $n\times n$ matrices $GL(n,\K)$,  that is, to a subgroup defined by polynomial equations. Such a group $G$ is called {\em unipotent}\, if its image under the embedding $G\rar GL(n,\K)$ lies in the set $\{g\in GL(n,\K)| (g-1)^n=0\}$.  An example of a unipotent group is the group  of all upper-triangular matrices in $GL(n,\K)$  with $1$'s on the main diagonal.

\sip

A group $\widehat{G}$ is called {\em pro-unipotent}\, if it is equal to the inverse limit,
$$
\widehat{G}={\lim_{\longleftarrow}} G_i:=\left\{\vec{g}\in\prod_{i\in I}G_i \ | \ f_{ij}(g_j)=g_i\right\},
$$
of some series of unipotent groups $\{G_i\}_{i\in I}$ parameterized by a directed set $(I,\leq)$ and equipped with a system of homomorphisms
$$
\left\{f_{ij}: G_j\rar G_i\right\}_{ i,j \in I\ \mathrm{such\ that}\ i\leq j}
$$
satisfying the conditions\footnote{Such a system $\{G_i, f_{ij}\}$ is called an {\em inverse system}\, of groups. Analogously one defines an inverse system of algebras, modules etc.}: $f_{ii}=\Id$ and  $f_{ik}=f_{ij}\circ f_{jk}$ for any $i\leq j\leq k$. The Lie algebra $\fg$ of the prounipotent group $G$ is defined as the inverse limit,
$$
\widehat{\fg}=\lim_{\longleftarrow} \fg_i
$$
of the Lie algebras $\fg_i$ of the groups $G_i$. There exist mutually inverse maps $\exp: \widehat{\fg} \rar
 \widehat{G}$ and
$\log: \widehat{G} \rar \widehat{\fg}$.

\begin{definition} The {\em prounipotent completion  over a field $\K$}\, of an abstract group
$G$ is a prounipotent group $\widehat{G}(\K)$  together with a homomorphism $i: G\rar \widehat{G}(\K)$ satisfying the following universality property: if $h: G \rar \widehat{H}$ is a homomorphism of $G$ into a prounipotent group $\widehat{H}$ over $\K$, then there is a unique homomorphism $f:  \widehat{G}(\K)\rar \widehat{H}$
such that $h$ factors as the composition
$
h: G \stackrel{i}{\lon}  \widehat{G}(\K) \stackrel{f}{\lon}  \widehat{H}
$.
\end{definition}
Prounipotent completions are often called {\em Malcev completions}\, in the literature. One denotes sometimes the prounipotent completion over $\K$ of a group $G$ simply by $\wh{G}$ omitting thereby the reference to $\K$.

\subsection{Quillen's construction of the prounipotent completion }\label{2: Quillen's construction} Let $G$ be an abstract group, and
$$
\K[G]:=\left\{\sum_{g\in G} \la_g g \ |\ \la_g\in \K \ \text{such\ that\ only\ finitely many}\ \la_g\ \text{can be non-zero}\right\}
$$
 its group algebra over a field $\K$. There is an algebra homomorphism,
$$
\Ba{rccc}
\var: &\K[G] &\lon &\K\\
&:   \sum_{g\in G}\la_g  g &\lon & \sum_{g\in G} \la_g,
\Ea
$$
called the {\em augmentation}. Let $I_\K:=\Ker\, \var$ be the augmentation ideal, and let
 $\wh{\K[G]}$ be the  $I_\K$-adic completion\footnote{ The $I_\K$-adic topology on $\K[G]$ is Hausdorff as the intersection of all powers of $I$ is the zero ideal. In this case one can define the $I$-adic topology as the unique  {\em metric topology} associated with the following metric on $\K[G]$:
$$
d(x,x') := ||x-x'||,
$$
where, for any $x\in \K[G]$, its norm $||x||$ is defined to be $2^{-n}$ with $n$ being the largest
natural number such that $x\in I^n$. The algebra  $\wh{\K[G]}$ can be identified with the metric completion of  $\K[G]$.} of $\K[G]$ (see \S {\ref{2: subsec on completed foltered}}).
 The ideal $I_\K$ is generated by elements
 of the form $g-\id$, where $\id$ is the unit in $G$ (often denoted by $1$ by abuse of notation).

\mip

The algebra $\K[G]$ can be maid into a Hopf algebra with the coproduct defined as follows
$$
\Ba{rccc}
\Delta: &\K[G] &\lon &\K[G]\ot \K[G]\\
&  \sum_{g\in G}\la_g  g  &\lon & \sum_{g\in G} \la_g g\ot g,
\Ea
$$
This coproduct is continuous in the $I_\K$-adic topology on $\K[G]$ and hence induces by continuity a coproduct
$$
\Delta : \wh{\K[G]} \lon   \wh{\K[G]}\, \wh{\ot}\,  \wh{\K[G]}
$$
where $\wh{\ot}$ is the completed tensor product. This gives  $\wh{\K[G]}$ the structure
of a complete filtered Hopf algebra.

\sip

Consider the set of group-like elements,
\Beq\label{2: grouplike}
\wh{G}:= \left\{x \in  \wh{\K[G]}: \Delta(x) = x\ot x \ \mbox{and}\ \var(x) = 1\right\}.
\Eeq
This is a subgroup of the group of units in $\wh{\K[G]}$.

\subsubsection{\bf Remark}
 If $char(\K)=0$, then the $\K$-powers of elements $g\in G$ are well-defined in the ring $\wh{\K[G]}$,
$$
g^\la = (\id - (\id - g))^\la :=\sum_{i=0}^\infty\frac{
\la (\la - 1) \cdots (\la - i + 1)}{i!}
(\id - g)^i.
$$
In some cases the group $\wh{G}$ is generated by such powers; the
corresponding Lie algebra $\wh{\fg}$ is generated by elements $\log g = \log(\id - (\id - g)):=\sum_{i=1}^\infty (\id-g)^i$.  Note that the above formula for $g^\la$ implies,
\Beq\label{2: conjugation of g to the power la}
hg^\la h^{-1}= (hgh^{-1})^\la,
\Eeq
and
$$
g^\la g^\mu=g^{\la+\mu}
$$
for any $g,h\in G$ and $\la, \mu\in \K$.

\subsubsection{\bf Theorem \cite{Ha}}\label{2: Theorem on prounipotent}{\em If $I_\K/I_\K^2$ is a finite-dimensional vector space over $\K$, then $\wh{G}$
is a prounipotent completion of $G$.}

\mip

If $G$ is finitely generated, then the condition of this Theorem is satisfied and the above explicit construction gives us the prounipotent completion of $G$.

\sip

It is often useful for practical computations to consider an increasing filtration of the group (\ref{2: grouplike}),
 $$
\ldots \subseteq  \wh{G}_{l+1} \subseteq \wh{G}_l \subseteq \ldots \subseteq \wh{G}_1\subseteq \wh{G},
 $$
defined by
$$
\wh{G}_l =\wh{G} \cap  (\id + \wh{I_\K}^l).
$$
If the condition of Theorem {\ref{2: Theorem on prounipotent}} is satisfied, then each group,
$$
\wh{G}^l:= \wh{G}/ \wh{G}_{l+1}
$$
is a unipotent group, and one can compute the prounipotent completion of $G$ as the inverse limit
$$
\wh{G}=\lim_{\longleftarrow} \wh{G}^l.
$$

\subsubsection{\bf The case $char(\K)=0$} Consider the set of {\em primitive
elements}, $\wh{\fg}$, of the completed group algebra $\wh{\K[G]}$ defined by
$$
\wh{\fg}:=\left\{x\in \wh{I_\K}\ |\ \Delta(x)=\id\ot x+ x\ot \id\right\}.
$$
The bracket $[x,y]=xy-yx$ makes $\wh{\fg}$ into a Lie algebra. The $\wh{I_\K}$-adic topology of $\wh{\K[G]}$
induces a topology on $\wh{\fg}$ making the latter into a complete topological Lie algebra. The logarithm and exponential functions give us well-defined and mutually inverse homeomorphisms,
$$
\log: \id + \wh{I_\K} \lon \wh{I_\K}\ \ \ \  \mbox{and}\ \ \ \  \exp :
\wh{I_\K}\lon \id + \wh{I_\K}
$$
which restrict to the mutually inverse homeomorphisms,
$$
\log: \wh{G} \lon \wh{\fg}\ \ \ \  \mbox{and}\ \ \ \  \exp : \wh{\fg}\lon \wh{G}.
$$
Thus one can describe $\wh{G}$ in terms of its Lie algebra $\wh{\fg}$ which is sometimes a simpler
object to define. The latter is the inverse limit of the Lie algebras of the prounipotent groups $\wh{G}^l$
defined above.

\subsection{Examples of prounipotent completions}\label{2: Examples of prounipotent completion}
\Bi
\item[(i)] Let $G$ be a group of finite type, that is, a quotient of the free group generated by a finite set $\{g_1,\ldots, g_n\}$ by the normal subgroup generated by a finite number of relations
$$
\left\{R_i(g_1,\ldots, g_n)=1\right\}_{i\in [m]}, \ \ \  m\in \N.
 $$
Let $\K$ be a field of characteristic zero, and let $\wh{\fg}$ be the quotient of the completed free Lie algebra $\wh{\ffree}_n$ generated by symbols $\{\ga_1,\ldots, \ga_n\}$  by the ideal generated by the relations,
$$
\left\{\log R_i(e^{\ga_1},\ldots, e^{\ga_n})=0\right\}_{i\in [m]}.
$$
Then the group $\wh{G}(\K):=\exp{\wh{\fg}}$ gives us the prounipotent completion of $G$ over $\K$
\cite{Br}.

\mip

\item[(ii)] If $F_n$ is the free group in $n$ letters, $x_1,\ldots, x_n$,
then its prounipotent completion, $\wh{F}_n(\K)$, over a field $\K$ of characteristic zero is equal
to $\exp(\wh{\ffree_n})$, where $\wh{\ffree_n}$ is the completed filtered Lie algebra over $\K$
generated by the following set,
$$
\left\{\log x_i=\log(1-(1-x_i))=\sum_{k=1}^\infty \frac{1}{k}(1-x_i)^k\right\}_{i\in [n]}.
$$
Put another way, $F_n(\K)$ is generated by elements of the form $x_i^\la$, $\la\in \K$.



\Ei

The next examples are taken from \cite{Kn}.

\Bi
\item[(iii)] Let $G=\Z$ and let $\K$ be a field of characteristic zero. As $\Z$ is isomorphic to the free group generated by a symbol $t$, its prounipotent completion $\wh{\Z}(\K)$ can be computed using the construction in Example (i) above to give
$$
\wh{\Z}(\K)= \exp (\wh{\ffree_1}) \simeq \K\ \ (\mbox{as Abelian groups}).
$$
One can get the same result using Quillen's construction as follows \cite{Kn}. The group algebra $\K[\Z]$
is the ring of Laurent polynomials $\K[t, t^{-1}]$ with the augmentation ideal $I_\K$ being a principal one generated by $t-1$. The completion $\wh{\K[\Z]}$ is the formal power series $\K[[T]]$
with the inclusion $\K[\Z]\rar \wh{\K[\Z]}$ given by
$$
\Ba{ccc}
\K[t, t^{-1}] &\lon & \K[[T]]\\
 (t, t^{-1}) & \lon & (1+T, \sum_{k=1}^\infty (-1)^k T^k).
\Ea
$$
The ideal $\wh{I_\K}$ is therefore the principal one generated by the symbol $T$ (that is, the maximal ideal of $\K[[T]]$)  so that the augmentation map
$\K[[T]]\rar \K$ has the form $T\rar 0$. The coproduct is given by $\Delta(T)=1\ot T + T\ot 1+T\ot T$.

\sip

If $\K$ has characteristic zero, then the set of group-like elements in $\wh{\K[\Z]}$, i.e.\ the prounipotent completion of $\Z$ over $\K$, is given by
$$
\wh{\Z}(\K)=\left\{(1+T)^\al= \exp(\al\log (1+T))\in \K[[T]]\ | \ \al\in \K\right\} \simeq \K.
$$
Note that in this case $\wh{\Z}(\K)^l=\wh{\Z}(\K)\cap (1+ \wh{I_\K}^l)=\{1\}$ for $l\geq 2$ so that
 $\wh{\Z}(\K)_l=\wh{\Z}(\K)$ for all $l$.

 \sip

If $\K$ has characteristic $p$, then the element $1+T$ has finite order in each group
$$
\wh{\Z}(\K)_l=\wh{\Z}(\K)/ \wh{\Z}(\K)\cap (1+ \wh{I_\K}^{l+1})
$$
as
$$
(1+T)^p=1+T^p\in 1+ \wh{I_\K}^p.
$$
It follows \cite{Kn} that the prounipotent completion of $\Z$ over $\K$ is equal to $\mathbf Z_p$, the Abelian group (in fact, the ring) of $p$-adic integers.

\mip

\item[(iv)] Let $G=\Z_n=\Z/n\Z$, $n\geq 2$. The associated group algebra over a field $\K$ is $\K[\Z_n]=\K[t]/(t^n-1)$, the quotient of the polynomial algebra $\K[t]$ by the principal ideal generated by $(t^n-1)$. The augmentation ideal  $I_\K \subset \K[\Z_n]$ is the principal one generated by $t-1$. Note that $t$ is a unit with $t^{-1}=t^{n-1}$. There is a factorization in $\K[\Z_n]$,
\Beqrn
(t-1)^2 &=& t^2-2t +1\\
        &=& t^n +t^2-2t\\
        &=& t(t^{n-2} +\ldots +t +2)(t-1).
\Eeqrn
If characteristic of $\K$ is prime to $n$ or equal to zero, then  $t^{n-2} +\ldots +t +2$
is also a unit in  $\K[\Z_n]$ so that $I_\K^2=I_\K$ and hence
$$
 \wh{\K[\Z_n]}=\lim_{\longleftarrow}  \K[\Z_n] / I_\K^l =\K
$$
Hence the prounipotent completion of $\Z_n$ in this case is the trivial group, $\wh{\Z_n}(\K)=\{1\}$ (cf.\ Example (iii) above).

If $char(\K)=p$ divides $n$, then one can assume\footnote{The prounipotent completion of the cartesian product of two groups is the cartesian product of their prounipotent completions.}  that $n=p^d$ for some $d\geq 1$. As $(1-t)^{p^d}=1-t^{p^d}=1-t^n=0$, one has $I_\K^{p^d}=0$
so that $\wh{\Z_n}(\K)=\Z_n(\K)$ and hence the prounipotent completion is given by the following set of group-like elements,
$$
\wh{\Z_n}(\K)=\{t^i\ |\ 0\leq i\leq p^d\}\simeq \Z_{p^d}
$$

\item[(v)] Examples (ii) and (iii) are special cases of a more general statement \cite{Kn}: if $G$ is
a finitely generated Abelian group, then
$$
\wh{G}(\K):= \left\{\Ba{ll}
G\ot_\Z \K& \mbox{if}\ char(\K)=0\\
G\ot_\Z \mathbf Z_p & \mbox{if}\ char(\K)=p.\Ea\right.
$$
\Ei

\subsection{Profinite completions}
Let $G$ be an abstract group and $\{I\}$ a family
of all normal subgroups $I\subset G$ of finite index. This family can be made into a partially ordered set with respect to the inclusion. As $I_1 I_2 \in \{I\}$ for any $I_1, I_2\in \{I\}$, this partially ordered set is directed so that $\{I\}$ gives us an inverse system of ideals and one can  consider the inverse limit,
$$
\widetilde{G}:=\lim_{\longleftarrow} G/I,
$$
which is called the {\em profinite completion of $G$}. For example, $\widetilde{\Z}=\prod_p \mathbf{Z}_p$.

\bip

{\large
\section{\bf Monoidal categories and monoidal functors}
}

\mip

For a category $\sC$ its class of objects is denoted by $Ob(\sC)$, and the set of morphisms
from an object $A$ to an object $B$ by $\Mor_\sC(A,B)$. If  $Ob(\sC)$ is a set, then the category is called {\em small}. The composition of morphisms $\Mor(A,B)\times \Mor(B,C)\rar  \Mor(A,C)$ is denoted by $\circ$.

\sip

\subsection{Monoidal categories (see e.g.\ \cite{ES})}
A category $\sC$ is called a {\em semigroup category}\, if it is equipped
with a bifunctor
$$
\Ba{rccc}
\cdot\ \ot_\sC\, \cdot & \sC\times \sC & \lon & \sC\\
                      &  A ,B\in Ob(\sC)       & \lon & A \ot_\sC B \in Ob(\sC),
\Ea
$$
and an isomorphism of trifunctors,
$$
\Ba{rccc}
\Phi:  &  (\ \cdot\ \ot_\sC\ \cdot\ )\ot_\sC\ \cdot  & \lon &  \cdot\ \ot_\sC (\ \cdot \ \ot_\sC\ \cdot\ )\\
\Phi_{A,B,C}:&   (A\ot_\sC B)\ot_\sC C       & \lon & A \ot_\sC (B \ot_\sC C)
\Ea
$$
such that the diagram
\Beq\label{3: Pentagon}
\Ba{c}\resizebox{100mm}{!}{
\xymatrix{
& ((A\ot_\sC B)\ot_\sC C)\ot_\sC D\ar[dl]_-{\Phi_{A\ot_\sC B, C, D}}
\ar[dr]^-{\Phi_{A,B,C}\ot \mathit{Id}} &   \\
(A\ot_\sC B)\ot_\sC (C\ot_\sC D)\ar[d]_-{\Phi_{A,B,C\ot_\sC D}} &  & (A\ot_\sC (B\ot_\sC C))\ot_\sC D\ar[d]^-{\Phi_{A,B\ot_\sC C,D}}\\
A\ot_\sC (B\ot_\sC (C\ot_\sC D)) & & A\ot_\sC ((B\ot_\sC C)\ot_\sC D)\ar[ll]_-{\mathit{Id}\ot \Phi_{B,C,D}}
}}\Ea
\Eeq
commutes for any $A,B,C,D\in Ob(\sC)$.

\mip

Commutativity of the diagram (\ref{3: Pentagon}) is called the {\em pentagon axiom}. A remarkable fact is that the pentagon axiom implies commutativity of all similar diagrams in the category $\cS$. More precisely, one has

\subsubsection{\bf Mac-Lane cohenerence theorem} {\em Let $\sC$ be a semigroup category and consider an arbitrary collection $A_1,\ldots, A_n$ of objects $\sC$. Given any two complete bracketings
of the formal expression $A_1\ot_\sC\ot A_2\ot_\sC \ldots \ot_\sC A_n$, then all isomorphisms from one bracketing to another composed of the associativity isomorphisms $\Phi$ and their inverses, are equal to each other.}

\mip

A semigroup category $\sC$ is called {\em monodial}\, if it has an object $\id_\sC$
(called a {\em unit}) together with isomorphisms
$$
\la_A: \id_\sC \ot_\sC A \lon  A, \ \ \ \ \ \rho_A: A \ot_\sC \id_\sC \lon  A
$$
such that the diagram
$$
\Ba{c}\resizebox{70mm}{!}{
\xymatrix{
 (A\ot_\sC \id_\sC)\ot_\sC B  \ar[dr]_-{\rho_A\ot Id} \ar[rr]^-{\Phi_{A,\ids_\sC,B}} && A\ot_\sC (\id_\sC\ot_\sC B)
  \ar[dl]^-{Id\ot \la_A}  \\
& A\ot_\sC B &
}}\Ea
$$
commutes for any $A,B\in Ob(\sC)$.

\mip

A monoidal category is called {\em strict}\, if all the isomorphisms $\Phi$, $\la$ and $\rho$ are identities.

\subsubsection{\bf Symmetric monoidal categories}
A monoidal category is called {\em symmetric}\, if it is equipped with an
isomorphism of bifunctors,
\Beq\label{3: beta}
\Ba{rccc}
\beta: & \cdot\ \ot_\sC \,\cdot & \lon & \cdot\ \ot_\sC^{op}\, \cdot \\
\beta_{A,B}& A\ot_\sC B & \lon & B\ot_\sC A
\Ea
\Eeq
such that
\Beq\label{3: beta2 constraint}
\be_{B,A}\circ \be_{A,B}=\mathit{Id}
\Eeq
and the diagram
\Beq\label{3: hexagon for symm monoidal}
\Ba{c}\resizebox{100mm}{!}{
\xymatrix{
 & A\ot_\sC (B \ot_\sC C)\ar[r]^-{\be_{A,B\ot_\sC C}} & (B\ot_\sC C)\ot_\sC A
 \ar[dr]^{\Phi_{B,C,A}}  & \\
(A\ot_\sC B)\ot_\sC C\ar[ur]^-{\Phi_{A,B,C}}\ar[dr]^-{\be_{A,B}\ot \mathit{Id}} &&& B\ot_\sC (C\ot_\sC A) \\
 & (B\ot_\sC A) \ot_\sC C\ar[r]^-{\Phi_{B,A,C}} & B\ot_\sC (A\ot_\sC C)
 \ar[ur]^{\mathit{Id}\ot \be_{A,C}}  & \\
}}\Ea
\Eeq
commutes for any $A,B,C\in Ob(\sC)$.

\mip

Commutativity of diagram (\ref{3: hexagon for symm monoidal}) is called the {\em hexagon axiom}.

\mip

There is a symmetric monoidal analogue of the above Mac Lane coherence theorem which says, that given
 any complete bracketings
of the formal expressions $A_1\ot_\sC\ot A_2\ot_\sC \ldots \ot_\sC A_n$ and
$A_{\sigma(1)}\ot_\sC\ot A_{\sigma(2)}\ot_\sC \ldots \ot_\sC A_{\sigma(n)}$, $\sigma\in \bS_n$, then all isomorphisms from one bracketing to another composed of the associativity isomorphisms $\Phi$, symmetry isomorphisms $\beta$  and their inverses, are equal to each other.
\mip

\subsubsection{\bf Braided monoidal categories}
A monoidal category is called {\em braided}\, if it is equipped with an isomorphism of functors (\ref{3: beta})
which makes the diagram  (\ref{3: hexagon for symm monoidal}) commutative while its inverse $\beta^{-1}$ makes a version of the diagram (\ref{3: hexagon for symm monoidal}) in which the symbols $\be_{A,B\ot_\sC C}$, $\be_{A,B}$ and $\be_{A,C}$ are replaced with
$\be^{-1}_{B\ot_\sC C,A}$, $\be^{-1}_{B,A}$ and  $\be^{-1}_{C,A}$ respectively commutative as well.

\sip

Mac Lane type coherence property of a braided monoidal category  will be explained below in terms of the braid groups.

\subsection{Examples of monodial categories}

\Bi

\item[{(i)}] Let $\mathsf{Vect}_\K$  be the category whose objects are  vector spaces\footnote{The  completed filtered versions of the categories listed below will be denoted by the same kernel symbol, but equipped with a wide hat, e.g., $\wh{\mathsf{Vect}}_\K$,  $\wh{\mathsf{Ass}}_\K$
    etc.}  over a field $\K$ and whose morphisms, $f: A\lon B$, are continuous  linear maps. This is a strict symmetric monoidal category  with
\Beqrn
A\ot_{\mathsf{Vect}_\K} B &:= & A \ot_\K B,\ \ \mbox{the ordinary (completed) tensor product of vector spaces over $\K$},\\
\id_{\mathsf{Vect}_\K} &:= & \K
\Eeqrn
and the symmetry morphism is given by
$\be_{A,B}(A\ot_\K B)=B\ot_\K A$.

\mip

\item[{(ii)}] Let $\mathsf{Ass}_\K$ be the
category of  associative algebras (with unit) over a field $\K$. This is a subcategory of $\mathsf{Vect}_\K$. Moreover, $\mathsf{Ass}_\K$ is a strict symmetric monoidal category with respect to the structures induced from $\mathsf{Vect}_\K$.

\mip

\item[{(iii)}] Let $\mathsf{Hopf}_\K$ be the
category of  Hopf algebras over a field $\K$. This a subcategory of $\mathsf{Ass}_\K$ which is  is  strict symmetric monoidal  with respect to the structures induced from $\mathsf{Ass}_\K$.

\mip

\item[{(iv)}] Let $\mathsf{coAss}_\K$ be the
category of  associative coalgebras over a field $\K$. This a subcategory of $\mathsf{Vect}_\K$ which is  is  strict symmetric monoidal with respect to the structures induced from $\mathsf{Vect}_\K$.

\mip

\item[{(v)}] Let $\mathsf{Lie}_\K$ be the category whose objects are (complete Lie algebras over a field $\K$ and whose morphisms, $f: \fg_1\lon \fg_2$, are linear maps preserving Lie brackets, i.e. $f\circ [\cdot ,\cdot ]=[f(\cdot), f(\cdot)]$.
This is a strict symmetric monoidal category  with
\Beqrn
 \fg \ot_{\mathsf{Lie}_\K} \fh &:= & \fg\oplus \fh,\ \ \mbox{the direct sum of vector spaces over $\K$},\\
\id_{\mathsf{Lie}_\K} &:= & 0,
\Eeqrn
and the symmetry morphism given by
$\be_{\fg,\fh}(\fg\oplus \fh)=\fh \oplus \fg$.

\mip

\item[(i-v)$'$] There are obvious {\em dg}, that is, {\em differential graded}, versions ---
$\mathsf{dgVect}_\K$, $\mathsf{dgAss}_\K$, $\mathsf{dgHopf}_\K$, $\mathsf{dgCoAss}_\K$, $\mathsf{dgLie}_\K$ ---
 of all the above symmetric monoidal categories whose objects are {\em dg}\, vector spaces, {\em dg}\, associative algebras etc.

 \sip

 We often use the following notation. Let $V=\oplus_{i\in \Z} V^i$ be a $\Z$-graded vector space, then for any integer $k\in\Z$ the symbol
$V[k]$ stands for the $\Z$-graded  vector space with $V[k]^i:=V^{i+k}$ and
and $s^k$ for the associated isomorphism $V\rar V[k]$; for $v\in V^i$ one denotes $|v|:=i$.
For a pair of $\Z$-graded vector spaces $V_1$ and $V_2$, the symbol $\Hom_i(V_1,V_2)$ stands
for the space of homogeneous linear maps of degree $i$, and
$\Hom(V_1,V_2):=\bigoplus_{i\in \Z}\Hom_i(V_1,V_2)$; for example, $s^k\in \Hom_{-k}(V,V[k])$.
Furthermore, we use the notation $\odot^n V$ for the n-fold symmetric product of the vector space $V$.

\mip
\item[{(vi)}] Let $\mathsf{Cat}$ be the category of small  categories (with functors as morphisms). This is a strict symmetric monoidal category with
$$
Ob(C'\ot_{\mathsf{Cat}} C'')  :=  Ob(C') \times Ob(C''),
 $$
 and
$$
  \Mor_{C'\ot_{\mathsf{Cat}} C''}(X'\times X'',Y'\times Y'')=  \Mor_{C'}(X',Y')\times
 \Mor_{C''}(X'',Y'').
$$
Its full subcategory consisting of small groupoids (that is, small categories with every morphism being an isomorphism) is denoted by $\mathsf{CatG}$.
\mip

\item[{(vii)}] Let $\mathsf{Cat(Vect_\K)}$ be a category
of small $\K$-linear
categories, that is, a category of small categories which satisfy two conditions: (1) morphisms between any pair of objects form a  vector space, and (2) the composition of morphisms is a bilinear
map.  This is a subcategory of $\mathsf{Cat}$ but we make it into a  symmetric monoidal category
in a slightly different way:
$$
Ob(C'\ot_{\mathsf{Cat(Vect_\K)}} C'')  :=  Ob(C') \times Ob(C''),
 $$
 and
$$
  \Mor_{C'\ot_{\mathsf{Cat(Vect_\K)}} C''}(X'\times X'',Y'\times Y'')=  \Mor_{C'}(X',Y')\ot_\K
 \Mor_{C''}(X'',Y'').
$$
For future reference we denote by $\mathsf{Cat(coAss}_\K)$  the (strict symmetric monoidal) subcategory of  $\mathsf{Cat(Vect}_\K)$ such that the space of morphisms between
any two objects is a  $\K$-coalgebra, and
the composition is compatible with the coalgebra structures.

\mip
\item[(viii)] Let $\mathsf{Set}$ be the category of sets. This is a strict symmetric monoidal category
with
\Beqrn
I\ot_{\mathsf{Set}} J &:= & I \times J,\ \ \mbox{the Cartesian product of sets $I$ and $J$},\\
\id_{\mathsf{Set}} &:= & \text{one element set},
\Eeqrn
and the symmetry morphism given by
$\be_{A,B}(I\times J)=J\times I$.

\mip
\item[{ (ix)}] Let $\mathsf{Group}$ be the subcategory of $\mathsf{Set}$ whose objects are groups and homomorphisms as morphisms. This is a strict symmetric monoidal category with respect to the the Catersian product; the unit is given by the trivial group $1$.

\item[{ (x)}] Let $\mathsf{Top}$ be the subcategory of $\mathsf{Set}$ whose objects are Hausdorff  spaces with compactly generated topology and continuous maps as morphisms. This is a strict symmetric monoidal category with respect to the structures inherited from $\mathsf{Set}$.
\mip

\Ei
\subsection{Monoidal functors} Let $(\sC, \ot_\sC, \id_\sC)$ and $(\sD, \ot_\sD, \id_\sD)$ be monoidal categories.
A functor $F: \sC\rar \sD$ is called {\em monoidal}\, if it comes with a morphism $\id: \id_\sD \rar F(\id_\sC)$
and a natural transformation
$$
\Ba{rccc}
\phi:  &  F(\, \cdot\, ) \ot_\sD F(\,\cdot\, )  & \lon &  F(\, \cdot\, \ot_\sC \, \cdot \, )\\
\phi_{A,B}:&   F(A)\ot_\sD F( B)       & \lon &  F(A \ot_\sC B)
\Ea
$$
such that the following diagrams,
$$
\Ba{c}\resizebox{90mm}{!}{
\xymatrix{
(F(A)\ot_\sD F(B))\ot_\sD F(C)\ \ \ \ \ar[r]^-{{\Phi_{F(A), F(B), F(C)}}}\ar[d]_-{\phi_{A, B}\ot Id} &\ \ \ \ F(A)\ot_\sD (F(B))\ot_\sD F(C)) \ar[d]^-{Id \ot \phi_{B, C}} \\
F(A\ot_\sC B)\ot_\sD F(C)\ar[d]_-{\phi_{A\ot_\sC B, C}} &   F(A)\ot_\sD  F(B\ot_\sC C)\ar[d]^-{\phi_{A, B\ot_\sC C}}\\
F((A\ot_\sC B)\ot_\sC C)\ar[r]_-{\Phi_{A, B, C}} &   F(A\ot_\sC (B\ot_\sC C))
}}\Ea
$$

$$
\Ba{c}\resizebox{44mm}{!}{
\xymatrix{
\id_\sD \ot_\sD F(A)  \  \ar[r]^-{\ids\ot Id}\ar[d]_{\la_{F(A)}} &\  F(\id_\sC)\ot_\sD F(A) \ar[d]^-{\phi_{\ids_\sC, A}} \\
F(A) &   F(\id_\sC \ot_\sC A))\ar[l]_{F(\la_A)}\\
}}\Ea
\ \ \ \ \ \
\Ba{c}\resizebox{44mm}{!}{\xymatrix{
F(A)\ot_\sD \id_\sD \  \ar[r]^-{Id\ot \id}\ar[d]_{\rho_{F(A)}} &\  F(A)\ot_\sD F(\id_\sC) \ar[d]^-{\phi_{A, \ids_\sC}} \\
F(A) &   F(A\ot_\sC \id_\sC))\ar[l]_{F(\rho_A)}\\
}}\Ea
$$
commute for any $A,B,C\in Ob(\sC)$.

\mip

A functor $F$ between braided monoidal categories is called {\em braided monoidal}\, if it is monoidal and makes the diagram
$$
\Ba{c}\resizebox{48mm}{!}{
\xymatrix{
F(A)\ot_\sD F(B)\ar[r]^\be\ar[d]_-{\phi_{A,B}} & F(B)
\ot_\sD F(A)\ar[d]^-{\phi_{B,A}}\\
F(A\ot_\sC B)\ar[r]_{F(\be)} & F(B\ot_\sC A)
}}\Ea
$$
commutative for any $A,B\in Ob(\sC)$.

\mip

A braided monoidal functor between symmetric monoidal categories is called {\em symmetric monoidal}.

\subsection{Examples of monoidal  functors}\label{3: tensor functors}
\Bi
\item[(i)] The contravariant functor
$$
\Ba{rccc}
Lin_\K: & \mathsf{Set} & \lon & \mathsf{Vect}_\K\\
     &   I          & \lon & \K[I]:= \mbox{the space of functions}\ I \rar \K
\Ea
$$
is symmetric monoidal.

\mip

\item[(ii)] The (completed) universal enveloping functor,
$$
\Ba{rccc}
U: & \mathsf{Lie}_\K & \lon & \mathsf{Hopf}_\K\\
     &   \fg          & \lon &  {U}(\fg)
\Ea
$$
is symmetric monoidal because ${U}(\fg\oplus \fh)={U}(\fg)\ot_\K U(\fh)$.

\mip

\item[(iii)] Consider a functor,
$$
\Delta_\K: \mathsf{CatG}  \lon \mathsf{Cat(coAss_\K)}
$$
which is identity on objects $I\in Ob(\mathsf{CatG})$,
$$
\Delta_\K(I)=I,
$$
and on morphisms it coincides with $Lin_\K$ defined above in (i),
$$
\Delta_\K(\Mor_{\mathsf{CatG}}(I,J))=\K[\Mor_{\mathsf{CatG}}(I,J)]
$$
where $\K[\Mor_{\mathsf{Cat}}(I,J)]$ is equipped with the unique coalgebra structure in which all elements $f$ of $\Mor_{\mathsf{Cat}}(I,J)$ are group-like, $\Delta(f)= f\ot_\K f$. This functor is symmetric monoidal. Note that $\Delta_\K(\bS_n)=\K[\bS_n]$, the group algebra of $\bS_n=\Mor_{\mathsf{CatG}}([n],[n])$ equipped with its standard coalgebra structure.

\mip

\item[(iv)] Let $(\fg,\ [\ ,\ ])$ be a Lie algebra, and $CE_\bu(\fg)$ its {\em Chevalley-Eilenberg complex}\, (also known as the {\em bar constriction}) which is, by definition, a dg coalgebra
$$
CE^\bu(\fg):= \odot^\bu(\fg[1])=\bigoplus_{n=0}^\infty (\ot^n\fg[1])_{\bS_n},
$$
 with the coproduct $\Delta: CE^\bu(\fg)\rar CE^\bu(\fg)\ot_\K CE^\bu(\fg)$ given by
$$
\Delta: s(x_{1})\odot\ldots \odot s(x_{n}):= \sum_{[n]=I_1\sqcup I_2\atop
\# I_1,I_2\geq 0} (-1)^{\sigma(I_1,I_2)} x_{I_1} \ot x_{I_2},
$$
where, for a naturally ordered subset $I=\{i_1,\ldots, i_k\}$ of $[n]$, we set
$$
x_I:=s(x_{i_1})\odot \ldots \odot s(x_{i_k}),
$$
and $(-1)^{\sigma(I_1,I_2)}$ is the sign of the permutation $[n]\rar [I_1\sqcup I_2]$.
The differential, $d: CE^\bu(\fg)\rar CE^\bu(\fg)$, is
given by
$$
d s(x_{1})\odot\ldots \odot s(x_{n}):=
\sum_{
1\leq i<j\leq n}
(-1)^{i+j} s([x_i, x_j])\odot s(x_1)\odot \ldots \odot \widehat{s(x_i)} \odot \ldots \odot
 \widehat{s(x_j)}\odot\ldots \odot {s(x_n)}.
$$
As $CE^\bu(\fg\oplus \fh)=CE^\bu(\fg) \ot_\K CE^\bu(\fh)$, the functor
$$
\Ba{rccc}
CE^\bu: & \mathsf{Lie}_\K & \lon & \mathsf{dgCoAss}_\K\\
     &   \fg          & \lon &  CE^\bu(\fg)
\Ea
$$
is symmetric monoidal.
\mip

\item[(v)] An associative $\K$-algebra $A$ with unit can be interpreted as a category
with one object, $\bu$, and with $\Mor(\bu,\bu):=A$. This defines a functor
$$
\mathsf{Ass}_\K \lon \mathsf{Cat(Vect_\K)}
$$
which is symmetric monoidal.

\mip

\item[(vi)] Let $I^n:=[0,1]^n\subset \R^n$ be the $n$-cube (we set $I^0$ to be a point), and $X$ a topological space. A continuous map $f: I^n\rar X$ is called a {\em singular $n$-cube}. A singular
    cube $f$ is called {\em degenerate}\, if the function $f$ does not depend on (at least one) one of the $n$ canonical coordinates, $(t^1,t^2,\ldots,t^n)$ on $I^n$. Define a $\Z$-graded vector space of (cubical) chains,
$$
Chains_{\bu}(X)=\bigoplus_{n=0}^\infty Chains_{-n}(X),
$$
by setting $Chains_{n}(X)$ to be the quotient of the free Abelian group consisting of formal linear combinations,
$$
\sum_i n_i f_i, \ \ n_i\in \Z,
$$
of  non-degenerate singular cubes $f_i: I^n\rar X$ modulo the relation,
$$
f= (-1)^\sigma f\circ \sigma,
$$
for any permutation $\sigma\in \bS_n$ acting on the cube as follows,
$$
\Ba{rccc}
\sigma: & I^n   &\lon &  I^n\\
        & (t_1,t_2,\ldots, t_n) & \lon & (t_{\sigma(1)},t_{\sigma(2)},\ldots, t_{\sigma(n)})
        \Ea
$$
The differential, $\p$, on  $Chains_{\bu}(X)$ is defined in the usual way,
$$
\Ba{rccl}
\p: & Chains_{n}(X)  &\lon &  Chains_{n-1}(X)\\
        & f(t_1,\ldots, t_n) & \lon & (\p f)(t_1,\ldots,t_{n-1}):=\sum_{i=1}^n(-1)^i
        f(t_1,\ldots,t_{i-1},0, t_i,\ldots, t_n)-\\
        &&&
        \hspace{35mm} - \sum_{i=1}^n(-1)^i f(t_1,\ldots,t_{i-1},1, t_i,\ldots, t_n)
        \Ea
$$
The homology of the complex $(Chains_{\bu}(X), \p)$ is denoted by $H_\bu(X,\Z)$ and is called the {\em singular (cubical) homology}\, of the topological space $X$. This homology is equal to the singular simplicial homology of $X$  (so that at the homology level the adjective {\em cubical}\, can be omitted). For a field $\K$ we denote $Chains_{\bu}(X,\K):= Chains_{\bu}(X)\ot_\Z \K$ and its homology by $H_\bu(X,\K)$.

\sip

\noindent For any finite collection $\{X_i\}_{i\in I}$ of topological
spaces,  there is a product
map
$$
\bigotimes_{i\in I}  Chains_\bu\left(X_i\right) \lon  Chains_\bu\left(\prod_{i\in I} X_i\right),
$$
which is a {\em natural}\,  homomorphism of complexes.

\sip

\noindent The cubical  chains functor
$$
\Ba{ccc}
\mathsf{Top} &\lon & \mathsf{dgVect}_\K\\
X& \lon & Chains_{\bu}(X,\K)
\Ea
$$ and the associated
homology functor
$$
\Ba{ccc}
\mathsf{Top} &\lon & \mathsf{dgVect}_\K\\
X& \lon & H_{\bu}(X,\K)
\Ea
$$
 are symmetric monoidal.

\mip

\item[(vii)] A {\em path}\, in a topological space $X$ from a point $x_0\in X$ to a point  $x_1 \in  X$ is a continuous map $u: [0,1]\rar X$ of the unit
interval into $X$ with $u(0) = x_0$ and $u_(1) = x_1$. Let $\pi(X)(x_0,x_1)$  stand for the set of of all homotopy classes of paths in $X$ from $x_0$ to $x_1$. There is a natural composition map,
$$
\pi(X)(x_0,x_1)\times \pi(X)(x_1,x_2) \lon \pi(X)(x_0,x_2),
$$
which makes the set of points of $X$ into a category, $\pi(X)$,
in which the set of morphisms from an object $x_0\in X$ to an object $x_1\in X$ is identified with $\pi(X)(x_0,x_1)$. Every such a morphisms
is an isomorphism so that $\pi(X)$ is in fact a groupoid called the {\em fundamental groupoid}\, of $X$. This construction gives us a monoidal functor,
\Beq\label{3: fundam groupoid monoidal functor}
\Ba{rccc}
\pi: & \mathsf{Top} &\lon & \mathsf{CatG}\\
& X& \lon & \pi(X).
\Ea
\Eeq
Note that for any subset $A\subset X$ one can define a full subcategory
of $\pi(X)$ whose objects are points of $A$; it is denoted by $\pi(X,A)$.

\mip

The set $\pi(X)(x_0,x_0)=:\pi_1(X,x_0)$ is a group called the
{\em fundamental group of $X$  base at $x_0$}. If $\mathsf{Top}_*$ stands for the
category of {\em based topological spaces}\, then $\pi_1$ gives us a
monoidal functor
\Beq\label{3: fundam group monoidal functor}
\Ba{rccc}
\pi_1:& \mathsf{Top}_* &\lon & \mathsf{Group}\\
& (X,x_0)& \lon & \pi_1(X,x_0)
\Ea
\Eeq
We shall use below both these functors to construct operads in the categories
$\mathsf{CatG}$ and $\mathsf{Group}$ out of operads in the categories
in $\mathsf{Top}$ and $\mathsf{Top}_*$ respectively.
\mip

\item[(viii)] Let $G$ be a (discrete) group. For non-empty subsets $A,B\subset G$, let
$$
[A, B]=\left\{[a,b]:= aba^{-1}b^{-1}\ |\ a\in A,\ b\in B\right\}
$$
stand for the subgroup of $G$ generated by all  commutators of elements of $A$ with elements of $B$.
Define inductively  the {\em descending central series} of $G$ as the following
series of normal subgroups,
$$
G_1:=G \ \supseteq \ G_2=[G,G]\ \supseteq\  G_3=[G_2,G]\ \supseteq \ldots  \supseteq
\ G_n=[G_{n-1},G]  \supseteq \
G_{n+1}=[G_{n},G]\ \supseteq \dots
$$
and let
$$
\fg\fr_n(G):= G_n/G_{n+1}
$$
 be the  associated $n$-th quotient. This is an Abelian group, i.e.\ a $\Z$-module, so that
 it makes sense to define the following positively graded  vector space over a field $\K$,
 $$
 \fg\fr_\K(G) := \bigoplus_{n=1}^\infty\fg\fr_n(G)\ot_\Z \K,
 $$
The commutator map
$$
\Ba{rccc}
[\ ,\ ]:& G\times G &\lon & G\\
        &  (a,b)     & \lon&  aba^{-1}b^{-1},
\Ea
$$
induces on $\fg\fr(G)_\K$ a linear skew-symmetric map,
$$
[\ ,\ ]: \wedge^2 \fg\fr_\K(G) \lon  \fg\fr_\K(G)
$$
 which
satisfies the Jacobi identity. Hence this construction (which is natural in $G$) gives us a functor
$$
\Ba{rccc}
\wh{\fg}\fr_\K: & \mathsf{Group} & \lon & \mathsf{Lie}_\K\\
& G &\lon & \wh{\fg\fr}_\K(G)
\Ea
$$
which is obviously symmetric monoidal. Here $\wh{\fg\fr}_\K(G)$ is the completion,
$\prod_{n=1}^\infty\fg\fr_n(G)\ot_\Z \K$,
 of the filtered Lie algebra  $\fg\fr_\K(G)$. It is called the (completed)
 {\em graded Lie algebra of $G$ over $\K$}.

\Ei

\subsubsection{\bf A useful fact}\label{3: fact on mono of gr(G)} Let $f: G\rar H$ be a morphism of groups, and assume that  $G$ is {\em residually nilpotent}, that is, its
descending central series  satisfies the condition
$$
\bigcap_{i=1}^\infty G_n=1.
$$
If the associated morphism of graded Lie algebras,
$$
\fg\fr^\Q(f): \fg\fr^{\Q}(G) \lon   \fg\fr^{\Q}(H)
$$
is a monomorphism, then $f$ is a monomorphism \cite{CW}.

\subsection{Closed symmetric monoidal categories} A symmetric monoidal category $\sC$ is called {\em closed} if for any object $A\in Ob(\sC)$  the functor
$$
\Ba{rccc}
\cdot\ \ot_\sC A & \sC & \lon & \sC \\
                 &  B & \lon & B\ot_\sC A
\Ea
$$
has a {\em right adjoint}\, functor,
$$
\Ba{rccc}
\Hom_\sC (A,\ \cdot\ ) & \sC & \lon & \sC \\
                 &  B & \lon &  \Hom_\sC (A,B)
\Ea
$$
which, by definition,  is a functor satisfying the condition,
$$
 \Mor_\sC(B\ot_\sC A, C) \cong  \Mor_\sC(B,  \Hom_\sC (A,C)),
$$
for any $A,B,C\in Ob(\sC)$. One can show  using Yoneda lemma  that this condition implies
\Beq\label{3: Int hom equation}
\Hom_\sC(B\ot_\sC A, C) \cong  \Hom_\sC(B,  \Hom_\sC (A,C)).
\Eeq
The object $ \Hom_\sC (A,B)\in Ob(\sC)$  is called the {\em internal hom}\, of $A$ and $B$.

\mip

The categories  $\mathsf{Cat}$, $\mathsf{Set}$ and the category of  finite-dimensional vector spaces are closed with
  $ \Hom_\sC (A,B)=\Mor_\sC(A,B)$.


\bip


{\large
\section{\bf Operads}
}

\mip

We shall mostly use in the subsequent sections the notion of operad at a rather elementary level --- at the level of its basic properties which follow from its definition almost immediately. However, the definition itself (which is due to P.\ May \cite{May}) is not, perhaps, very elementary, and we recommend a newcomer to pay attention not only to the formal definition(s) of operad given below, but also to examples illustrating that definition.  Comprehensive expositions of the theory of operads can be found  in the books \cite{LV,MSS}.

\sip

From now on $\sS$ stands for the groupoid  of finite sets and their bijections. This is a subcategory of
$\mathsf{Set}$.

\mip

\subsection{Basic axioms for operadic compositions}\label{4: subsec on def of operad}
An {\em operad}\, $\f$ in a symmetric monoidal category $\sC$ is a
collection of the following data \cite{Ta2}:

\mip

1) a functor $\f: \sS \rar \sC$ (which is often called an $\sS$-{\em module});

\mip

2) for any finite sets $I$, $J$ and an element $i\in I$  a morphism
\Beq\label{4: operadic composition}
\circ_i^{I,J}: \f(I) \ot_\sC \f(J) \lon \f(I\setminus\{i\}\, \sqcup J ),
\Eeq
\ \ natural in $i$, $I$, $J$ (``the insertion of $\f(J)$ into the $i$-th slot of $\f(I)$"), such that
the following

\ \ associativity conditions hold:
\Bi
\item[(a)] for any finite sets $I$, $J$, $K$ and any $i_1,i_2\in I$, $i_1\neq i_2$,  the diagram

$$
\Ba{c}\resizebox{120mm}{!}{
\xymatrix{
\f(I)\ot_\sC \f(J)\ot_\sC \f(K) \ \ \ar @{->}[r]^-{\circ_{i_1}^{I,J}\ot \mathit{Id}}
  \ar[d]_{\mathit{Id}\ot \beta_{\f(J)\f(K)}} & \ \ \f(I\setminus\{i_1\}\, \sqcup J )\ot_\sC \f(K)\ \ \
\ar[r]^-{\circ_{i_2}^{I\setminus\{i_1\}\sqcup J,K}} & \ \ \   \f(I\setminus \{i_1,i_2\}\sqcup J\sqcup K)  \ar @{=}[d] \\
\f(I)\ot_\sC \f(K)\ot_\sC \f(J) \ \  \ar[r]^-{\circ_{i_2}^{I,K}\ot \mathit{Id}} & \ \
\f(I\setminus\{i_2\}\, \sqcup K )\ot_\sC \f(J)\ \ \ \ar[r]^-{\circ_{i_1}^{I\setminus\{i_2\}\sqcup K,J}} & \ \ \ \f(I\setminus \{i_1,i_2\}\sqcup K\sqcup J)
}}\Ea
$$
commutes;
\item[(b)] for any finite sets $I$, $J$, $K$ and any $i\in I$ and $j\in J$ the diagram
$$
\Ba{c}\resizebox{114mm}{!}{
\xymatrix{
 & \f(I)\ot _\sC\f(J\setminus\{j\}\, \sqcup K)
\ar[dr]^-{\circ_{i}^{I, J\setminus\{j\}\sqcup K}} &    \\
\f(I)\ot_\sC \f(J)\ot_\sC \f(K)  \ar[ur]^-{\mathit{Id}\ot\circ_{j}^{J,K}}
  \ar[dr]_{\circ_{i}^{I,J}\ot \mathit{Id}} &&  \f(I\setminus \{i\}\sqcup J\setminus \{j\} \sqcup K)\\
& \ \
\f(I\setminus\{i\}\, \sqcup J )\ot_\sC \f(K)\ \ \ \ar[ur]^-{\circ_{j}^{I\setminus\{i\}\sqcup J,K}} &
}}\Ea
$$
commutes.
\Ei

Sometimes we abbreviate compositions $\circ_i^{I,J}$ to $\circ_i$.

\subsection{Axioms for a unit}
An {\em operad with unit}\, is, by definition, an operad $\f$ equipped with
a morphism,
$$
e: \id_\sC \lon \f(\bu),
$$
for any one-element set $\{\bu\}$ which is natural in $\{\bu\}$ and makes the following two compositions,
$$
\f(I)\lon \f(I)\ot_\sC\id_\sC \xrightarrow{\mathit{Id}\ot e} \f(I)\ot_\sC \f(\bu) \xrightarrow{\circ_i^{I,\bu}} \f(I),
$$
$$
\f(I)\lon \id_\sC \ot_\sC \f(I) \xrightarrow{e\ot \mathit{Id}} \f(\bu) \ot_\sC \f(I)  \xrightarrow{\circ_\bu^{\bu,I}}\f(I),
$$
into identity maps for any set $I$ and any $i\in I$.

\subsubsection{\bf Remark}  The skeleton of the groupoid $\sS$ consists of the sets
$[n]$, $n\geq 0$ (with $[0]:=\emptyset$). If the category $\sC$ admits small colimits
(and all categories we work with in this paper do have this property), then a functor
$\f: \sS \rar \sC$ is equivalent to a collection of $\bS_n$-modules,
$\{\f(n):=\f([n]\}_{n\geq 0}$, which is often abbreviated to an $\bS$-module in the category $\sC$.
One can reformulate axioms of an operad in terms of  the
$\bS$-module $\{\f(n):=\f([n]\}_{n\geq 0}$ and the operadic compositions,
$$
\circ_i: \f(n)\ot_\sC \f(m)\lon \f([n]\setminus i\sqcup [m])=\f(n+m-1),\ \ \forall\ m,n\in \N, i\in[n],
$$
where one first identifies
$$
[n]\setminus i\sqcup [m]=\{1,2,i-1,i+1,\ldots, n\}\sqcup\{1',2',\ldots, m'\}
$$
 with the linearly ordered  set,
$$
\{1,2,\ldots, i-1,1',2',\ldots, m',i+1,\ldots, n\},
$$
 and then identifies the latter with
$[n+m-1]=\{1,2,\ldots,i-1,i,i+1,\ldots, m+n-1\}$ as linearly ordered sets.
In this way one recovers the original definition
of an operad given in \cite{May}. We shall often use this freedom below to
describe an operad either as a functor $\f: \sS \rar \sC$  or simply as an $\bS$-module
 $\{\f(n)\in Ob(\sC)\}_{n\in \N}$.

 \subsubsection{\bf Non-$\bS$ operads} As noted just above an operad
 $\f$
 in a monoidal category $\sC$ can be defined as a collection of objects,
 $\f=\{\f(n)\in Ob(\sC)\}_{n\geq 0}$, such that  each objects $\f(n)$ carries a representation
 of $\bS_n$, and  there are equivariant compositions,
 $$
 \circ_{i}: \f(n)\ot_\sC \f(m)\lon \f(n+m-1)\ \ \ \forall \ i\in[n],
 $$
which satisfy axioms (a) and (b) in \S{\ref{4: subsec on def of operad}}.

\sip

If we forget in the definition above $\bS_n$ actions on $\f(n)$, $n\geq 0$, and, correspondingly,
omit the equivariance condition on the compositions
$\circ_i$, then we get a notion of a {\em non-$\bS$ operad}. More precisely, a {\em non-$\bS$ operad}\,
in a (semigroup) category $\sC$ is an operad $\f=\{\f(n)\}$ in $\sC$ such that the action of $\bS_n$ on each
object $\f(n)$ is trivial.

\subsubsection{\bf Exercise} Let $\f$ be an operad in a symmetric monoidal
category $\sC$ and $F: \sC \rar \sD$ a symmetric monoidal functor to some other
category. Show that the  data
\Bi
\item[(i)]  the $\sS$-module structure, $F\f: \sS \rar \sD$, given by the composition $\sS\xrightarrow{\f} \sC\xrightarrow{F} \sD$, and
\item[(ii)] the operadic ``insertions",
$$
\bar{\circ}_i^{I,J}: F\f(I) \ot_\sC F\f(J) \lon F\f(I\setminus\{i\}\, \sqcup J ),
$$
given by the compositions
$$
 F\f(I) \ot_\sC F\f(J) \xrightarrow{\phi_{\f(I),\f(J)}} F( \f(I) \ot_\sC \f(J))
 \xrightarrow{F( \circ_i^{I,J})} F(\f(I\setminus\{i\}\, \sqcup J ))=
 F\f(I\setminus\{i\}\, \sqcup J ),
$$
\Ei
give us an operad $F\f$ in the symmetric monoidal category $\sD$. This fact is
of an extreme importance in applications --- starting with a ``geometric" operad in the category, say, of topological spaces, and applying the chain or homology functor one arrives to an operad in the category of vector spaces. This particular property of operads is another manifestation of the {\em  amazing unity of mathematics}.

\subsubsection{{\bf Exercise}} Introduce the notion of
 a {\em morphism},
 $
 \rho: \f_1 \lon \f_2,
 $
 of operads in a symmetric monoidal category $\sC$.

 \subsubsection{{\bf Exercise}} Introduce the notion of an {\em ideal}\, $\cI$, of an operad $\f$ in the category $\mathsf{Vect}_\K$ and construct the {\em quotient}\, operad
 $\f/\cI$.

\subsection{Basic examples}\label{4: basic examples of operads} Here we define a few operads which will be used in some constructions below.

\subsubsection{\bf Endomorphism operad} Let $\sC$ be a closed symmetric monoidal category. An arbitrary object $A\in Ob(\sC)$ gives rise to an operad $\cE nd_A$ with the $\bS$-module,
$$
    \left\{\cE nd_A(n):= \Hom_\sC(A^{\ot n}, A) \right\}
$$
and with the compositions $\circ_i^{[n],[m]}$ given by the isomorphisms (use \ref{3: Int hom equation})),
$$
\Hom_\cC\left(A^{\ot (i-1)}\ot_\sC \Hom_\sC(A^{\ot m},A)\ot_\sC A^{\ot(n-i+1)}, A\right)
\cong \Hom_\sC(A^{\ot(n+m+1)}, A),
$$
i.e.\ literally by the substitution of the object $\Hom_\sC(A^{\ot m},A)$ into the $i$-th input slot of $\Hom_\sC(A^{\ot n},A)$.

\sip

Let $\f$ be an operad in $\sC$. An $\f$-{\em algebra}\, structure on an
object $A$ is, by definition, a morphism of operads,
$$
\rho: \f \lon \cE nd_A,
$$
which is also often called a {\em representation}\, of $\f$ in $A$.

\mip

\subsubsection{\bf Operad of parenthesized permutations}\label{4: operad Pa}  Consider an $\sS$-module in the category $\mathsf{Set}$,
$$
\Ba{rccl}
\cP a: & \sS & \lon & \mathsf{Set}\\
    &  I   & \lon & \cP a(I):= \mbox{the set of all parenthesized permutations of}\ I
\Ea
$$
One can equivalently define $\cP a(I)$ as
\Bi
\item[-] the set of all monomials
built from elements of $I$ using a non-commutative and non-associative product (using each element of $I$ once),
\item[-] the set of all planar binary trees whose legs are labelled
by elements of $I$ (using each element once).
\Ei
For example, for $I=\{a,b,c,d,e\}$,
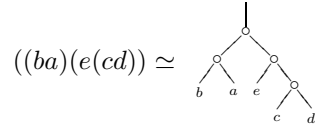
$$
((ba)(e(cd)) \simeq
\Ba{c}\resizebox{16mm}{!}{\xy
(-8,-3.7)*{_b},
(-2,-3.7)*{_a},
(2,-3.7)*{_e},
(5.5,-8.2)*{_c},
(11.5,-8.2)*{_d},
(0,12)*{}="0",
 (0,7)*{\circ}="a",
(-5,2)*{\circ}="b_1",
(5,2)*{\circ}="b_2",
(-8,-2)*{}="c_1",
(-2,-2)*{}="c_2",
(2,-2)*{}="c_3",
(8.5,-2.5)*{\circ}="c_4",
(5.5,-6.5)*{}="d_1",
(11.5,-6.5)*{}="d_2",
\ar @{-} "a";"0" <0pt>
\ar @{-} "a";"b_1" <0pt>
\ar @{-} "a";"b_2" <0pt>
\ar @{-} "b_1";"c_1" <0pt>
\ar @{-} "b_1";"c_2" <0pt>
\ar @{-} "b_2";"c_3" <0pt>
\ar @{-} "b_2";"c_4" <0pt>
\ar @{-} "c_4";"d_1" <0pt>
\ar @{-} "c_4";"d_2" <0pt>
\endxy}
\Ea
$$
 is an element of $\cP a(I)$. The operadic compositions $\circ_i^{I,J}=\circ_i$ are given by the substitution of a monomial from $\cP a(J)$ into the $i$-th letter of a monomial $\cP a(I)$, e.g.\
for  $I=\{a,b,c,d,e\}$ and $J=\{x,y,z\}$,
$$
    ((ba)(e(cd)) \   \circ_c \   ((xy)z)   =   ((ba)(e( ((xy)z)d)) \in \cP
     a\left(\{I\setminus c\}\sqcup J\right)
      $$
or simply by grafting the root (which in our pictorial notation grows upward) of a planar binary tree from   $\cP a(J)$ into the $i$-labelled input leg of a tree from     $\cP a(I)$,
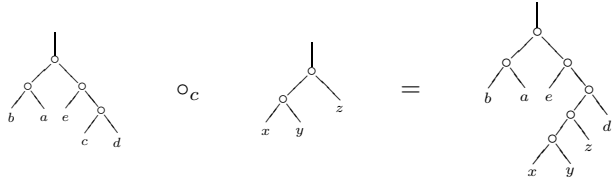
 $$
\Ba{c}\resizebox{16mm}{!}{
\xy
(-8,-3.7)*{_b},
(-2,-3.7)*{_a},
(2,-3.7)*{_e},
(5.5,-8.2)*{_c},
(11.5,-8.2)*{_d},
(0,12)*{}="0",
 (0,7)*{\circ}="a",
(-5,2)*{\circ}="b_1",
(5,2)*{\circ}="b_2",
(-8,-2)*{}="c_1",
(-2,-2)*{}="c_2",
(2,-2)*{}="c_3",
(8.5,-2.5)*{\circ}="c_4",
(5.5,-6.5)*{}="d_1",
(11.5,-6.5)*{}="d_2",
\ar @{-} "a";"0" <0pt>
\ar @{-} "a";"b_1" <0pt>
\ar @{-} "a";"b_2" <0pt>
\ar @{-} "b_1";"c_1" <0pt>
\ar @{-} "b_1";"c_2" <0pt>
\ar @{-} "b_2";"c_3" <0pt>
\ar @{-} "b_2";"c_4" <0pt>
\ar @{-} "c_4";"d_1" <0pt>
\ar @{-} "c_4";"d_2" <0pt>
\endxy}
\Ea
     \ \ \ \   \circ_c \ \ \ \
\Ba{c}\resizebox{12mm}{!}{
 \xy
(-8,-3.7)*{_x},
(-2,-3.7)*{_y},
(5,0.3)*{_z},
(0,12)*{}="0",
 (0,7)*{\circ}="a",
(-5,2)*{\circ}="b_1",
(5,2)*{}="b_2",
(-8,-2)*{}="c_1",
(-2,-2)*{}="c_2",
\ar @{-} "a";"0" <0pt>
\ar @{-} "a";"b_1" <0pt>
\ar @{-} "a";"b_2" <0pt>
\ar @{-} "b_1";"c_1" <0pt>
\ar @{-} "b_1";"c_2" <0pt>
\endxy}\Ea
\ \ \ \
= \ \ \ \
\Ba{c}\resizebox{18mm}{!}{
\xy
(-8,-3.7)*{_b},
(-2,-3.7)*{_a},
(2,-3.7)*{_e},
(8.5,-11.9)*{_z},
(11.5,-8.2)*{_d},
(-0.7,-15.9)*{_x},
(5.5,-15.9)*{_y},
(0,12)*{}="0",
 (0,7)*{\circ}="a",
(-5,2)*{\circ}="b_1",
(5,2)*{\circ}="b_2",
(-8,-2)*{}="c_1",
(-2,-2)*{}="c_2",
(2,-2)*{}="c_3",
(8.5,-2.5)*{\circ}="c_4",
(5.5,-6.5)*{\circ}="d_1",
(11.5,-6.5)*{}="d_2",
(2.5,-10.5)*{\circ}="e_1",
(8.5,-10.5)*{}="e_2",
(-0.7,-14.5)*{}="f_1",
(5.5,-14.5)*{}="f_2",
\ar @{-} "a";"0" <0pt>
\ar @{-} "a";"b_1" <0pt>
\ar @{-} "a";"b_2" <0pt>
\ar @{-} "b_1";"c_1" <0pt>
\ar @{-} "b_1";"c_2" <0pt>
\ar @{-} "b_2";"c_3" <0pt>
\ar @{-} "b_2";"c_4" <0pt>
\ar @{-} "c_4";"d_1" <0pt>
\ar @{-} "c_4";"d_2" <0pt>
\ar @{-} "d_1";"e_1" <0pt>
\ar @{-} "d_1";"e_2" <0pt>
\ar @{-} "e_1";"f_1" <0pt>
\ar @{-} "e_1";"f_2" <0pt>
\endxy}
\Ea
 $$
 \subsubsection{{\bf Exercises}}\label{4: Exercise on operads 1} (i)  Show that there is a one-to-one correspondence between
  $\Mor_{\mathsf{Set}}(X\times X, X)$ and $\cP a$-algebra structures on $X \in Ob(\mathsf{Set})$,
  that is,  representations, $\rho:\cP a \rar \cE nd_X$,
  of the operad $\cP a$ in a set $X$.

\sip

(ii) Describe representations of the operad $Lin_\K(\cP a)$ in a vector space $V$, where
 $Lin_\K$ is a monoidal functor $\mathsf{Set}\rar \mathsf{Vect}_\K$ defined
in \S {\ref{3: tensor functors}}(i).

\sip

(iii) Let $\f_1$ and $\f_2$ be operads in a symmetric monoidal category $\sC$. Show that
the $\sS$-module defined by,
$$
\Ba{rccc}
\f_1\ot_\sC \f_2:& \sS & \lon & \sC\\
                & I& \lon & \f_1(I)\ot_\sC \f_2(I),
\Ea
$$
has an induced operadic structure. The resulting operad $\f_1\ot_\sC \f_2$ is called the {\em tensor product}\, of operads $\f_1$ and $\f_2$.

\subsubsection{\bf Graphs and rooted trees}\label{4: subsec on graphs}
A {\em graph with legs (or hairs)}\,  is a triple $\Ga=(H(\Ga), \sqcup, \tau )$, where
\Bi
\item[-] $H(\Ga)$ is a finite set whose elements
are called {\em half-edges}\ (or {\em flags}),
\item[-] $\sqcup$ is a partition of $H(\Ga)$ into a disjoint union of subsets,
$$
H(\Ga)=\coprod_{v\in V(\Ga)} H(v),
$$
parameterized by a set $V(\Ga)$  which is called the {\em set of vertices of $\Ga$}; the subset $H(v)\subset H(\Ga)$ is called the {\em set of half-edges attached to the vertex $v$};
its cardinality, $\# H(v)$, is called the {\em valency}\, of $v$,
\item[-] $\tau: H(\Ga)\rar H(\Ga)$ is an involution, that is, a map satisfying the condition $\tau^2=Id$. This map defines a new partition of $H(\Ga)$ into orbits $(h,\tau(h))\subset H(\Ga)$ which in general can have cardinality  two ($h\neq \tau(h)$) or one ($h= \tau(h)$). Orbits of cardinality $2$ are called  {\em internal edges}\, or simply {\em edges}\, of the graph $\Ga$; the set of (internal) edges is denoted by  $E(\Ga)$. Orbits of cardinality $1$ are called {\em legs}\, (or {\em hairs}); the set of legs is denoted by $L(G)$.
\Ei
Note that $L(\Ga)=\emptyset$ if and only the involution $\tau$ has no fixed points. In this case $\Ga=(H(\Ga,\sqcup, \tau)$ is called simply a {\em graph}. We shall study operads of such graphs in \S 7 below.

\sip

It is convenient to think of a graph (with or without legs) $\Ga$  in terms of its {\em geometric realization}\, which is a topological space constructed as follows: (i) for each vertex $v$ take $\# H(v)$ copies of the interval $[0,1]$
labelled by elements of $H(v)$ and glue these intervals at the end-point $0$, the result is a topological space (equipped with the quotient topology from $[0,1]^{\# H(v)}$) which is
called the {\em corolla of $v$}; (ii) then consider a union of all corollas and, for each internal edge $(h,\tau(h))$, identify the end-points $1$ of the intervals $[0,1]$ labelled by $h$ and $\tau(h)$.

\sip

An {\em isomorphism}\, of graphs with legs $i: \Ga_1\rar \Ga_2$ is a bijection
$i: H(\Ga_1)\rar H(\Ga_2)$ which preserves partitions $\sqcup_1$ and $\sqcup_2$ (and hence induces
a bijection of the sets of vertices $V(\Ga_1)\rar V(
\Ga_2)$) and commutes with involutions, $i\circ \tau_1= \tau_2\circ i$
(and hence induces bijections $E(\Ga_1)\rar E(\Ga_2)$ and $L(\Ga_1)\rar L(\Ga_2)$). The associated groupoid of graphs and their isomorphisms is denoted by $\mathsf{IsoGraph}$.

\sip

A graph (with legs) is called connected (resp., simply connected) if its geometric realization is
a connected (resp., simply connected) topological space. A connected and simply connected
graph is called a {\em tree}. A {\em rooted tree}\, is a tree $T$ with $L(T)\neq \emptyset$ and with one one of the legs,
say $r\in L(T)$,  marked and called {\em the root}. Elements in $L(T)\setminus r$ are called
{\em  input legs}\, or {\em leaves} of the rooted tree $T$; the set of leaves is denoted by
$Leaf(T)$. Each vertex on a geometric realization of a rooted tree is connected to a root
by a unique path; this defines a flow on the tree and hence a partition,
$$
H(v)=In(v)\sqcup r_v.
$$
 of the set
of half-edges at each vertex $v$ into a disjoint union of a unique element $r_v$ which lies on the path from $v$ to the root and the subset of the remaining half-edges, $In(v)$, which are called
{\em input (half) edges}. By abuse of language, the cardinality of $In(v)$ is called the {\em valency}
\, of a vertex $v$ of a {\em rooted}\, tree $T$ (which is equal to the valency of $v$ in $T$ viewed as an unrooted tree minus one). The corolla of a vertex in a rooted tree can be visualized  as follows (with flow directed from bottom to the top)
$$
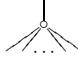
\Ba{c}\resizebox{10mm}{!}{\xy
(0,0)*{\circ}="a",
(0,5)*{}="0",
(0,-5)*{\ldots},
(4,-5)*{}="b_1",
(7,-5)*{}="b_2",
(-4,-5)*{}="b_3",
(-7,-5)*{}="b_4",
\ar @{-} "a";"0" <0pt>
\ar @{-} "a";"b_1" <0pt>
\ar @{-} "a";"b_2" <0pt>
\ar @{-} "a";"b_3" <0pt>
\ar @{-} "a";"b_4" <0pt>
 \endxy}\Ea
$$

Let $T_1$ and $T_2$ be rooted trees. An isomorphism of graphs, $i: T_1\rar T_2$, which sends the root of
$T_1$ into the root of $T_2$ is called an {\em isomorphism of rooted trees}. The associated groupoid of rooted trees and their isomorphisms is denoted by  $\mathsf{IsoTree}$.

\subsubsection{\bf Free operads} Let $\mathsf{IsoTree}(I)$ be a category (in fact, a groupid) whose objects are
rooted trees $T$ equipped with an isomorphism $l_T: Leaf(T) \rar I$ (``labeling of input
legs with elements of $I$") and morphisms, $f: T_1\rar T_2$, are isomorphisms of rooted trees
which respect the labelings, i.e.\ satisfy the condition $l_{T_1}= l_{T_2}\circ f$.
An example of such a morphism is given by the following picture

\Beq\label{4: morphism of labeled trees}
\Ba{c}\resizebox{18mm}{!}{
\xy
(-12,-11.7)*{_1},
(-2,-11.7)*{_2},
(2,-11.7)*{_3},
(12,-11.7)*{_4},
(-2,10)*{_{v_1}},
(-10,0)*{_{v_2}},
(10,0)*{_{v_3}},
(-5.6,5)*{_{e_{12}}},
(5.9,5)*{_{e_{13}}},
(-8.5,-6)*{_{k}},
(-5,-6)*{_{l}},
(9,-6)*{_{n}},
(5.8,-6)*{_{m}},
(0,20)*{}="0",
 (0,10)*{\circ}="a",
(-7,0)*{\circ}="b_1",
(7,0)*{\circ}="b_2",
(-12,-10)*{}="c_1",
(-2,-10)*{}="c_2",
(2,-10)*{}="c_3",
(12,-10)*{}="c_4",
\ar @{-} "a";"0" <0pt>
\ar @{-} "a";"b_1" <0pt>
\ar @{-} "a";"b_2" <0pt>
\ar @{-} "b_1";"c_1" <0pt>
\ar @{-} "b_1";"c_2" <0pt>
\ar @{-} "b_2";"c_3" <0pt>
\ar @{-} "b_2";"c_4" <0pt>
\endxy}
\Ea
\underset{{(v_1,v_2,v_3)\rar (v_1, v_2,v_3)\atop
(e_{12}, e_{13})\rar (e_{13},e_{12})}\atop
(k,l,m,n)\rar (n,m,k,l)
 }{\xrightarrow{\ \ \ \ \ \ \ \ f \ \  \ \ \ \ }}
\Ba{c}\resizebox{18mm}{!}{
\xy
(-12,-11.7)*{_3},
(-2,-11.7)*{_4},
(2,-11.7)*{_2},
(12,-11.7)*{_1},
(-2,10)*{_{v_1}},
(-10,0)*{_{v_2}},
(10,0)*{_{v_3}},
(-5.6,5)*{_{e_{12}}},
(5.9,5)*{_{e_{13}}},
(-8.5,-6)*{_{k}},
(-5,-6)*{_{l}},
(9,-6)*{_{n}},
(5.8,-6)*{_{m}},
(0,20)*{}="0",
 (0,10)*{\circ}="a",
(-7,0)*{\circ}="b_1",
(7,0)*{\circ}="b_2",
(-12,-10)*{}="c_1",
(-2,-10)*{}="c_2",
(2,-10)*{}="c_3",
(12,-10)*{}="c_4",
\ar @{-} "a";"0" <0pt>
\ar @{-} "a";"b_1" <0pt>
\ar @{-} "a";"b_2" <0pt>
\ar @{-} "b_1";"c_1" <0pt>
\ar @{-} "b_1";"c_2" <0pt>
\ar @{-} "b_2";"c_3" <0pt>
\ar @{-} "b_2";"c_4" <0pt>
\endxy}
\Ea
\Eeq
where $V(T)=\{v_1,v_2,v_3\}$, $E(T)=\{e_{12}, e_{21}\}$ and $Leaf(T)=\{k,l,m,n\}$.

\sip

 Given any $ \bS$-module $\cE :\sS \rar \sC$ in a symmetric monoidal category $\sC$ (with small colimits) there is an associated {\em free}\, operad
$\cF ree\langle \cE \rangle$ which has an obvious universal property. For simplicity, we shall explain the construction of $\cF ree\langle \cE \rangle$ in the case $\sC=\mathsf{Vect}_\K$, the general case being essentially the same.

\sip

Given an  $ \bS$-module,  $\cE :\sS \rar \mathsf{Vect}_\K$,
to any $I$-labelled rooted tree $T \in Ob(\mathsf{IsoTree}(I))$ with, say, $n$ vertices, one can associate a vector space (``unordered tensor product over the set of vertices of $T$"),
\[
T\langle \cE \rangle=
\bigotimes_{v\in V(T)} \cE(In_v):=\left(
 \bigoplus_{p: [n]\rar V(T)} \cE\left(In_{p(1)}\right)\ot\ldots\ot
\cE\left(In_{p(n)}\right)\right)_{\bS_n}
\]
(in a general symmetric monoidal category one has to take the equalizer over the set
of isomorphisms $[n]\rar V(T)$). It is not hard to see directly from this definition that the association $T \rar T\langle \cE \rangle$ defines a functor from the category $\mathsf{IsoTree}(I)$ to the category $\mathsf{Vect}_K$, and hence it makes sense to consider its colimit
$$
\underset{T\in\mathsf{IsoTree}(I)}{\mathrm{colim}} T\langle \cE\rangle,
$$
which exists in $\mathsf{Vect}_K$. As a vector space this colimit is (non-canonically) isomorphic to the direct sum,
$$
\underset{T\in\mathsf{IsoTree}(I)}{\mathrm{colim}} T\langle \cE\rangle \simeq \bigoplus_{T \in Ob(\mathsf{IsoTree}(I))/\sim} T\langle \cE \rangle
$$
where the summation runs over representatives $T$ of isomorphism classes of $I$-labelled rooted trees from the category $\mathsf{IsoTree}(I)$.
We use this colimit to define an $\sS$-module,
$$
\Ba{rccc}
\cF ree\langle \cE \rangle: & \sS & \lon &  \mathsf{Vect}_\K\\
&                             I   & \lon &  \cF ree\langle \cE \rangle(I):= \underset{T\in
\mathsf{IsoTree}(I)}{\mathrm{colim}} T\langle \cE\rangle,
\Ea
$$
and notice that it has  a natural structure of the operad with respect to graftings of trees and unordered tensor products. This operad is called the {\em free operad generated by the $\sS$-module $\sE$}.

\sip

Elements of  $\cF ree\langle \cE \rangle(I)$ can be visualized as linear combinations
 of {\em $\cE$-decorated $I$-labelled
rooted trees}. One can represent pictorially such a decorated tree as follows,
$$
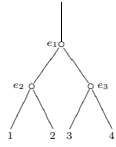
\Ba{c}\resizebox{15mm}{!}{
\xy
(-12,-11.7)*{_1},
(-2,-11.7)*{_2},
(2,-11.7)*{_3},
(12,-11.7)*{_4},
(-2,10)*{_{e_1}},
(-10,0)*{_{e_2}},
(10,0)*{_{e_3}},
(0,20)*{}="0",
 (0,10)*{\circ}="a",
(-7,0)*{\circ}="b_1",
(7,0)*{\circ}="b_2",
(-12,-10)*{}="c_1",
(-2,-10)*{}="c_2",
(2,-10)*{}="c_3",
(12,-10)*{}="c_4",
\ar @{-} "a";"0" <0pt>
\ar @{-} "a";"b_1" <0pt>
\ar @{-} "a";"b_2" <0pt>
\ar @{-} "b_1";"c_1" <0pt>
\ar @{-} "b_1";"c_2" <0pt>
\ar @{-} "b_2";"c_3" <0pt>
\ar @{-} "b_2";"c_4" <0pt>
\endxy}
\Ea
$$
where $e_i$ are vectors from $\cE(In_{v_i})$, $i=1,2,3$. Such a representation is not unique; for example,  isomorphism (\ref{4: morphism of labeled trees}) leads to the identification
$$
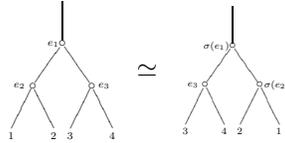
\Ba{c}\resizebox{15mm}{!}{
\xy
(-12,-11.7)*{_1},
(-2,-11.7)*{_2},
(2,-11.7)*{_3},
(12,-11.7)*{_4},
(-2,10)*{_{e_1}},
(-10,0)*{_{e_2}},
(10,0)*{_{e_3}},
(0,20)*{}="0",
 (0,10)*{\circ}="a",
(-7,0)*{\circ}="b_1",
(7,0)*{\circ}="b_2",
(-12,-10)*{}="c_1",
(-2,-10)*{}="c_2",
(2,-10)*{}="c_3",
(12,-10)*{}="c_4",
\ar @{-} "a";"0" <0pt>
\ar @{-} "a";"b_1" <0pt>
\ar @{-} "a";"b_2" <0pt>
\ar @{-} "b_1";"c_1" <0pt>
\ar @{-} "b_1";"c_2" <0pt>
\ar @{-} "b_2";"c_3" <0pt>
\ar @{-} "b_2";"c_4" <0pt>
\endxy}
\Ea
\simeq
\Ba{c}\resizebox{15mm}{!}{
\xy
(-12,-11.7)*{_3},
(-2,-11.7)*{_4},
(2,-11.7)*{_2},
(12,-11.7)*{_1},
(-4,10)*{_{\sigma(e_1)}},
(-10,0)*{_{e_3}},
(11.5,0)*{_{\sigma(e_2)}},
(0,20)*{}="0",
 (0,10)*{\circ}="a",
(-7,0)*{\circ}="b_1",
(7,0)*{\circ}="b_2",
(-12,-10)*{}="c_1",
(-2,-10)*{}="c_2",
(2,-10)*{}="c_3",
(12,-10)*{}="c_4",
\ar @{-} "a";"0" <0pt>
\ar @{-} "a";"b_1" <0pt>
\ar @{-} "a";"b_2" <0pt>
\ar @{-} "b_1";"c_1" <0pt>
\ar @{-} "b_1";"c_2" <0pt>
\ar @{-} "b_2";"c_3" <0pt>
\ar @{-} "b_2";"c_4" <0pt>
\endxy}
\Ea
$$
\sip
\noindent
where $\sigma$ is the unique non-trivial automorphism of a two element set $I_0$ (say $I_0=\{e_{12},e_{13}\}$, $\{k,l\}$ or $\{m,n\}$, see (\ref{4: morphism of labeled trees})), and  $\sigma(e_i)$
is the result of an action of this automorphism on the corresponding elements of $\cE(I_0)$.

\sip

Every element in   $\cF ree\langle \cE \rangle(I)$ is an iterated operadic composition of (isomorphism classes of) the following basic decorated graphs,
$$
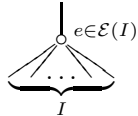
\underbrace{
\xy
(0,0)*{\circ}="a",
(0,5)*{}="0",
(0,-5)*{\ldots},
(4,-5)*{}="b_1",
(7,-5)*{}="b_2",
(-4,-5)*{}="b_3",
(-7,-5)*{}="b_4",
\ar @{-} "a";"0" <0pt>
\ar @{-} "a";"b_1" <0pt>
\ar @{-} "a";"b_2" <0pt>
\ar @{-} "a";"b_3" <0pt>
\ar @{-} "a";"b_4" <0pt>
 \endxy}_I
 \hspace{-6mm}{^{e\in \cE(I)}}
$$
called the {\em generating corollas}.

\subsubsection{\bf Exercises} (i) Show that the operad of parenthesized permutations, $\cP a$,
is the free operad in the category $\mathsf{Set}$ generated by the following $\sS$-module,
$$
\cE(I):=\left\{\Ba{ll}
Bij(I\rar [\# I]) & \mbox{if}\ \# I =2\\
0  &\mbox{otherwise}.
\Ea\right.
$$
where $Bij(I\rar [\# I])$ is the set of all bijections $I\rar  [\# I]$.

\mip

(ii)
Show that the operad $Lin_\K(\cP a)$ (see \S~{\ref{4: Exercise on operads 1}}ii) is the free operad in the category $\mathsf{Vect}_K$ generated by the following $\bS$-module,
$$
\cE(n):=\left\{\Ba{ll}
\K[\bS_2] & \mbox{if}\ n=2\\
0  &\mbox{otherwise}.
\Ea\right.
$$

\subsubsection{\bf Operad of little disks}\label{4: subsect on little disks}  Let
$$
D_{a,\la}:=\{z\in \C\ | \ |z-a|\leq \la\}
$$
be the closed disk in $\C$ of radius $\la$ and with center $a\in \C$.
The {\em operad of little disks}, $\caD=\{\caD(n)\}_{n\geq 1}$, is an operad in the category, $\mathsf{Top}_{pc}$,
of path-connected topological spaces given by the following data,
\Bi
\item[-] $\caD(1)$ is the point (the unit in $\caD$);
\item[-] $\caD(n)$ for $n\geq 2$
 is the space of configurations of $n$ disjoint closed disks, $\{D_{a_i,\la_i}\}_{i\in [n]}$, lying
 inside the standard unital disk, $D_{0,1}$.
\item[-] the operadic composition, for any $k\in [n]$, is given by
$$
\Ba{rccc}
\circ_k: & \caD(n)\times \caD(m) & \lon & \caD(n+m-1)\\
&    \{D_{a_i,\la_i}\}_{i\in [n]} \times \{D_{b_j,\nu_i}\}_{i\in [m]} &\lon&
\displaystyle \left\{ \{D_{a_i,\la_i}\}_{i\in [n]\setminus k} \coprod \{D_{\la_k b_j+a_k ,\la_k \nu_j}\}_{j\in [m]}\right\}.
\Ea
$$
Put another way, we replace the closed $k$-th disk $D_{a_k,\la_k}$ in a configuration from  $\caD(n)$ with the
$\la_k$-rescaled and $a_k$-translated configuration from  $\caD(m)$.
\footnote{A pictorial description of this composition can be found at http://en.wikipedia.org/wiki/Operad\_theory}
\Ei

The homology  monoidal functor (see \S{\ref{3: tensor functors}}vi) sends the topological operad $\caD$ into an operad,
$H_\bu(\caD)$, in the category $\mathsf{dgVect}_\K$. It was proven by F.\ Cohen that $H_\bu(\caD)$ is precisely the operad
of {\em Gerstenhaber algebras}, i.e.\ its
representations $H_\bu(\caD)\rar \cE nd_V$ in a graded vector space $V$ are equivalent
to a {\em Gerstenhaber algebra}\, \footnote{See
http://en.wikipedia.org/wiki/Gerstenhaber\_algebra for the definition of this notion.} structure on $V$.

\sip

The {\em fundamental groupoid}\, monoidal functor $\pi: \mathsf{Top}_{*}\rar \mathsf{CatG}$
(see \S{\ref{3: tensor functors}}vii) sends $\caD$ into an operad
$\pi(\caD)$ in the category of groupoids.




\subsection{Degree shift functor}\label{4: subsec on degree shoft of operads}
 Let $\f$ be an operad in the category
$\mathsf{dgVect}_\K$. For any integer $m$ one can uniquely associate to $\f$ an
operad $\f\{m\}$ with the property that there is a 1-1 correspondence,
$$
\left\{\Ba{c}\mbox{reprersentations of}\ \f\\
\mbox{in a graded vector space}\ V \Ea\right\} \stackrel{1:1}{\longleftrightarrow}
\left\{\Ba{c}\mbox{reprersentations of}\ \f\{m\}\\
\mbox{in a graded vector space}\ V[m] \Ea\right\}.
$$
From this property one easily reads the structure of $\f\{m\}=\{\f\{m\}(n)\}_{n\geq 1}$
as an $\bS$-module,
$$
\f\{m\}(n)=\f(n)[m(1-n)]\ot sgn_n^{\ot m}.
$$
where $sgn_n$ is the one-dimensional sign representation $\bS_n$.

\sip

Put another way, $\f\{m\}$ is the tensor product of the operad $\f$ and the endomorphism
operad $\cE nd_{\K[-m]}$ of the 1-dimensional vector space concentrated in degree $m$.

\subsection{Cosimplicial structures on operads} We shall often use below the following observation.

\subsubsection{\bf Lemma}\label{4: Lemma on cosimpl operad structures} {\em
 Let $\f=\{\f(n)\}_{n\geq 1}$ be an operad in a symmetric monoidal category $\sC$. Assume there is an element $e\in \f(2)$ satisfying the condition
 \Beq\label{4: ass cond for e}
 e\circ_1 e=e\circ_2 e.
 \Eeq
Then the family of maps,
\Beq\label{4: formulae for d_i}
\Ba{rccc}
d_i: &  \f(n) & \lon & \f(n+1)\\
     &  f & \lon & d_i f:=\left\{\Ba{cc} e \circ_2 f & \mbox{for}\ i=0\\
     f\circ_i e & \mbox{for}\ i\in [n]\\
     e\circ_1 f & \mbox{for}\ i=n+1\\
     \Ea\right.
\Ea
\Eeq
satisfies the equations
\Beq\label{4: cosimplisial relations in an operad}
d_jd_i = d_id_{j-1} \ \ \ \mbox{for any}\ \ i < j,
\Eeq
and hence makes the collection  $\f=\{\f(n)\}_{n\geq 1}$ into a cosimiplicial object
in the category $\sC$.
}
\begin{proof}
Using the second definition of operads in terms of decorated trees, we shall identify
elements $e\in \f(2)$ and $f\in \f(n)$ with decorated corollas,
$$
e=
\Ba{c}\resizebox{6mm}{!}{\xy
(0,0)*{\circ}="a",
(0,4)*{}="0",
(-2,0.5)*{_e},
(-3,-4.5)*{_1},
(3,-4.5)*{_2},
(3,-3)*{}="b_1",
(-3,-3)*{}="b_2",
\ar @{-} "a";"0" <0pt>
\ar @{-} "a";"b_1" <0pt>
\ar @{-} "a";"b_2" <0pt>
 \endxy}\Ea=:
 \Ba{c}\resizebox{6mm}{!}{ \xy
(0,0)*{\bu}="a",
(0,4)*{}="0",
(-3,-4.5)*{_1},
(3,-4.5)*{_2},
(3,-3)*{}="b_1",
(-3,-3)*{}="b_2",
\ar @{-} "a";"0" <0pt>
\ar @{-} "a";"b_1" <0pt>
\ar @{-} "a";"b_2" <0pt>
 \endxy}\Ea, \ \ \ \ \ \ \
f=\Ba{c}\resizebox{11mm}{!}{\xy
(0,0)*{\circ}="a",
(-2,0.5)*{_f},
(-6,-5.5)*{_1},
(-3,-5.5)*{_2},
(6,-5.5)*{_n},
(0,4)*{}="0",
(0,-4)*{...},
(3,-4)*{}="b_1",
(6,-4)*{}="b_2",
(-3,-4)*{}="b_3",
(-6,-4)*{}="b_4",
\ar @{-} "a";"0" <0pt>
\ar @{-} "a";"b_1" <0pt>
\ar @{-} "a";"b_2" <0pt>
\ar @{-} "a";"b_3" <0pt>
\ar @{-} "a";"b_4" <0pt>
 \endxy}\Ea=:
 \Ba{c}\resizebox{11mm}{!}{\xy
(0,0)*{\blacktriangle}="a",
(-6,-5.5)*{_1},
(-3,-5.5)*{_2},
(6,-5.5)*{_n},
(0,4)*{}="0",
(0,-4)*{...},
(3,-4)*{}="b_1",
(6,-4)*{}="b_2",
(-3,-4)*{}="b_3",
(-6,-4)*{}="b_4",
\ar @{-} "a";"0" <0pt>
\ar @{-} "a";"b_1" <0pt>
\ar @{-} "a";"b_2" <0pt>
\ar @{-} "a";"b_3" <0pt>
\ar @{-} "a";"b_4" <0pt>
 \endxy}\Ea
$$
Then
$$
d_1d_0f=d_1\left(
 \Ba{c}\resizebox{12mm}{!}{\xy
(0,0)*{\bu}="a",
(0,4)*{}="0",
(-3,-4.5)*{_1},
(3,-3)*{}="a_1",
(-3,-3)*{}="a_2",
(3,-4)*{\blacktriangle}="x",
(-1,-10.5)*{_2},
(9,-10.5)*{_{n+1}},
(3,-9)*{...},
(5,-9)*{}="b_1",
(7,-9)*{}="b_2",
(1,-9)*{}="b_3",
(-1,-9)*{}="b_4",
\ar @{-} "x";"b_1" <0pt>
\ar @{-} "x";"b_2" <0pt>
\ar @{-} "x";"b_3" <0pt>
\ar @{-} "x";"b_4" <0pt>
\ar @{-} "a";"0" <0pt>
\ar @{-} "a";"a_1" <0pt>
\ar @{-} "a";"a_2" <0pt>
 \endxy}
\Ea
\right)=
 \Ba{c}\resizebox{15mm}{!}{
\xy
(0,0)*{\bu}="a",
(0,4)*{}="0",
(-6,-10.5)*{_1},
(-2,-10.5)*{_2},
(0.5,-10.5)*{_3},
(3,-3)*{}="a_1",
(-3,-3)*{}="a_2",
(-4,-4)*{\bu}="y",
(-6,-9)*{}="y_1",
(-2,-9)*{}="y_2",
(4,-4)*{\blacktriangle}="x",
(10,-10.5)*{_{n+2}},
(4,-9)*{...},
(6,-9)*{}="b_1",
(8,-9)*{}="b_2",
(2,-9)*{}="b_3",
(0,-9)*{}="b_4",
\ar @{-} "x";"b_1" <0pt>
\ar @{-} "x";"b_2" <0pt>
\ar @{-} "x";"b_3" <0pt>
\ar @{-} "x";"b_4" <0pt>
\ar @{-} "a";"0" <0pt>
\ar @{-} "a";"a_1" <0pt>
\ar @{-} "a";"a_2" <0pt>
\ar @{-} "y";"y_1" <0pt>
\ar @{-} "y";"y_2" <0pt>
 \endxy}
 \Ea
$$
On the other hand,
$$
d_0d_0f=d_0\left(
 \Ba{c}\resizebox{12mm}{!}{
\xy
(0,0)*{\bu}="a",
(0,4)*{}="0",
(-3,-4.5)*{_1},
(3,-3)*{}="a_1",
(-3,-3)*{}="a_2",
(3,-4)*{\blacktriangle}="x",
(-1,-10.5)*{_2},
(9,-10.5)*{_{n+1}},
(3,-9)*{...},
(5,-9)*{}="b_1",
(7,-9)*{}="b_2",
(1,-9)*{}="b_3",
(-1,-9)*{}="b_4",
\ar @{-} "x";"b_1" <0pt>
\ar @{-} "x";"b_2" <0pt>
\ar @{-} "x";"b_3" <0pt>
\ar @{-} "x";"b_4" <0pt>
\ar @{-} "a";"0" <0pt>
\ar @{-} "a";"a_1" <0pt>
\ar @{-} "a";"a_2" <0pt>
 \endxy}
\Ea
\right)= \Ba{c}\resizebox{15mm}{!}{
\xy
(-3,3)*{\bu}="o",
(-6,0)*{}="o_1",
(0,0)*{\bu}="a",
(-3,7)*{}="0",
(-3,-4.5)*{_2},
(-6,-1.5)*{_1},
(3,-3)*{}="a_1",
(-3,-3)*{}="a_2",
(3,-4)*{\blacktriangle}="x",
(-1,-10.5)*{_3},
(9,-10.5)*{_{n+2}},
(3,-9)*{...},
(5,-9)*{}="b_1",
(7,-9)*{}="b_2",
(1,-9)*{}="b_3",
(-1,-9)*{}="b_4",
\ar @{-} "x";"b_1" <0pt>
\ar @{-} "x";"b_2" <0pt>
\ar @{-} "x";"b_3" <0pt>
\ar @{-} "x";"b_4" <0pt>
\ar @{-} "o";"0" <0pt>
\ar @{-} "a";"a_1" <0pt>
\ar @{-} "a";"a_2" <0pt>
\ar @{-} "a";"o" <0pt>
\ar @{-} "o_1";"o" <0pt>
 \endxy}
 \Ea
$$
The condition (\ref{4: ass cond for e}) reads
$
 \Ba{c}\resizebox{8mm}{!}{
\xy
(-3,3)*{\bu}="o",
(0,0)*{}="o_1",
(-6,0)*{\bu}="a",
(-3,7)*{}="0",
(-3,-4.5)*{_2},
(-9,-4.5)*{_1},
(-3,-3)*{}="a_1",
(-9,-3)*{}="a_2",
(0,-1.5)*{_3},
\ar @{-} "o";"0" <0pt>
\ar @{-} "a";"a_1" <0pt>
\ar @{-} "a";"a_2" <0pt>
\ar @{-} "a";"o" <0pt>
\ar @{-} "o_1";"o" <0pt>
 \endxy}
 \Ea
 =
  \Ba{c}\resizebox{8mm}{!}{
\xy
(-3,3)*{\bu}="o",
(-6,0)*{}="o_1",
(0,0)*{\bu}="a",
(-3,7)*{}="0",
(-3,-4.5)*{_2},
(-6,-1.5)*{_1},
(3,-3)*{}="a_1",
(-3,-3)*{}="a_2",
(3,-4.5)*{_3},
\ar @{-} "o";"0" <0pt>
\ar @{-} "a";"a_1" <0pt>
\ar @{-} "a";"a_2" <0pt>
\ar @{-} "a";"o" <0pt>
\ar @{-} "o_1";"o" <0pt>
 \endxy}
 \Ea
$
so that we get $d_1d_0f=d_0d_0f$ for any $f\in \f(n)$. Similar pictures establish
identities  (\ref{4: cosimplisial relations in an operad}) for all other values of $i\in[n+1]$ and $j\in [n+2]$ with $i<j$. We leave it to the reader to draw the details.
\end{proof}

If $\f$ is an operad in the category $\mathsf{Vect}_\K$, then under the conditions of the above Lemma one can associate to $\f$ a dg vector space,
$$
\mathsf{Simp}^\bu(\f):=\bigoplus_{n\geq 1}\f(n)[-n],
$$
equipped with the differential
$$
\Ba{rccc}
d& \f(n) & \lon & \f(n+1)\\
 &   f & \lon & df:= \sum_{i=0}^{n+1}(-1)^i d_i f,
\Ea
$$
whose  cohomology is denoted by $H\mathsf{Simp}^\bu(\f)$ and is called the {\em cosimplicial}\, cohomology of $\f$.  We shall see below that, for example, the cosimplicial cohomology of the operad of infinitesimal braids has much to do with the main topic of these lectures,
 the Grothendieck-Teichm\"uller group.

\subsection{Lie algebras associated with operads}\label{4: subset on Kapr-Manin} Let $\f=\{\f(n)\}_{n\geq 1}$ be an operad in the category $\mathsf{dgVect}_\K$ with operadic compositions $\circ_i: \f(n)\ot \f(m) \rar \f(m+n-1)$, $1\leq i\leq n$. Consider a vector space
 $$
\f^{tot}:= \bigoplus_{n\geq 1}\f(n)
$$
Then the map
$$
\Ba{rccc}
[\ ,\ ]:&  \f^{tot} \ot \f^{tot} & \lon & \f^{tot}\\
& (a\in \f(n), b\in \f(m)) & \lon &
[a, b]:= \sum_{i=1}^n a\circ_i b - (-1)^{|a||b|}\sum_{i=1}^m b\circ_i a
\Ea
$$
makes $
\f^{tot}$
into a dg Lie algebra \cite{KM}. Moreover, the bracket descends to the space
of coinvariants
$
(\f^{tot})_\bS:=  \bigoplus_{n\geq 1}\f(n)_{\bS_n}$ making  this vector spaces into a dg Lie algebra as well. If $char(\K)=0$, the space of coinvariants $(\f^{tot})_\bS$ is isomorphic
 via the symmetrization map to the space of  $\bS$-invariants,
$(\f^{tot})^\bS:=  \bigoplus_{n\geq 1}\f(n)^{\bS_n}$, so that one has an induced
dg Lie algebra structure on  $(\f^{tot})^\bS$ with Lie bracket
given by the above formula followed with the symmetrization.


\mip

{\large
\section{\bf Operad of parenthesized  braids and $\wh{GT}$}
}

\sip
\subsection{Geometric definition of the braid group} Let $C_n(\C)$ be the configuration space
of pairwise distinct points in the plane $\C$,
$$
C_n(\C)=\left\{(z_1,z_2,\ldots, z_n)\in \C\ |\ z_i\neq z_j\ \mbox{for}\ i\neq j\right\}.
$$
The group $\bS_n$ acts on $C_n(\C)$ by permuting the points,
$$
\sigma: (z_1,z_2,\ldots, z_n)  \lon \left(z_{\sigma(1)},z_{\sigma(2)},\ldots, z_{\sigma(n)}\right), \ \ \ \sigma\in \bS_n.
$$
Let $C_n(\C)/\bS_n$ stand for the associated set of orbits equipped with the quotient topology.
Let $p_n=(z_1^0,\ldots, z_n^0)$  be any point in $C_n(\C)$ (say $p_n=\{(1,2,\ldots, n)\in \C\}$) and let $\bar{p}_n$ be
its image under the projection $C_n(\C)\rar C_n(\C)/\bS_n$.

\sip

\subsubsection{\bf Definitions} (i) The fundamental group $\B_n:=\pi_1(C_n(\C)/\bS_n, \bar{p}_n)$
is called the {\em group of braids}.

  (ii)  The fundamental group  $\P\B_n:= \pi_1(C_n(\C), p_n)$) is called
  the {\em group of pure braids}.

  \mip
\subsubsection{\bf Remark} It is often useful not to distinguish isomorphic groups.  As
topological spaces $C_n(\C)$ and $C_n(\C)/\bS_n$ are path-connected, the
groups $ \pi_1(C_n(\C), p_n)$ and $\pi_1(C_n(\C)/\bS_n, \bar{p}_n)$ are isomorphic to each for
different choices of base points. Therefore one should view  $\B_n$ and $\P\B_n$ as
isomorphism classes, $\pi_1(C_n(\C)/\bS_n)$ and, respectively,
 $\pi_1(C_n(\C))$, of the fundamental groups $\pi_1(C_n(\C)/\bS_n, \bar{p}_n)$ and, respectively,
 $\pi_1(C_n(\C), p_n)$; this point of view makes  the choice of the base points
 irrelevant.
  \mip

The projection $C_n(\C)\rar C_n(\C)/\bS_n$ is a regular $n!$-sheeted covering with $\bS_n$ acting as covering transformations. Thus the subgroup $\P\B_n\subset \B_n$ has index $n!$, and we have a short exact sequence of groups,
$$
1\lon
\P\B_n \lon \B_n  \stackrel{p}{\lon}\bS_n\lon  1.
$$

\mip

 Let $b$ be an element in $\B_n$ and $\sigma_b:=p(b)\in \bS_n$ the associated permutation.
One can visualize an element  $b\in \B_n$  as  (an isotopy classes of)
 $n$ disjoint continuous curves (strands),
 $$
 s_i: [0,1] \lon \C\times [0,1], \ \ \ i=1,2,\ldots, n,
 $$
 such that
$$
  s_i(0)=z_i^0, \ \ \
 s_i(1)=z^0_{\sigma_b(i)},
$$
and  the composition
$$
 [0,1] \stackrel{s_i}\lon \C\times [0,1] \stackrel{proj}{\lon} [0,1]
$$
is a continuous monotonous function (with flow along the strands running up). For example,
$$
b_1=
\Ba{c}\resizebox{8mm}{!}{
\xy
(-5,-10)*{_1},
(-5,10)*{_1},
(5,-10)*{_2},
(5,10)*{_2},
\vtwist~{(-5,8)}{(5,8)}{(-5,-8)}{(5,-8)};
\endxy}\Ea
\ \ \ \ \ \mbox{and}\ \ \ \ \
b'=
\Ba{c}\resizebox{12mm}{!}{\xy
(-10,-10)*{_1},
(-10,10)*{_1},
(1,-10)*{_2},
(-2,10)*{_2},
(7,-10)*{_3},
(6,10)*{_3},
(5,8)*{}; (-10,-8)*{} **\crv{(6,0)&(-10,3)}
\POS?(.4)*+{}="x" \POS?(.6)*+{}="y";
(-3,8)*{}; "x" **\crv{};
(-10,8)*{}; "y" **\crv{};
"y"+(.5,-1); (.5,-8)*{} **\dir{-};
"x"+(.5,-1); (6,-8)*{} **\dir{-};
\endxy}\Ea
$$
represent braids $b_1\in \B_2$, $b'\in \B_3$ with $\sigma_{b_1}=(12)$ and  $\sigma_{b'}=\left(\Ba{ccc} 1&2&3\\3&1&2 \Ea\right)$, while
$$
b''=
\Ba{c}\resizebox{8mm}{!}{
\xy
(-5,-17)*{_1},
(-5,7)*{_1},
(5,-17)*{_2},
(5,7)*{_2},
\vtwist~{(-5,5)}{(5,5)}{(-5,-5)}{(5,-5)};
\vtwist~{(-5,-5)}{(5,-5)}{(-5,-15)}{(5,-15)};
\endxy}
\Ea
$$
represents a pure braid from $\P\B_2$. Multiplication of braids is represented by the concatenation of strands,
$$
xy=
\Ba{c}\resizebox{11mm}{!}{
\xy
(7.5,12.2)*{y},
(0,10); (15,10)*{} **\dir{-};
(0,15); (15,15)*{} **\dir{-};
(0,10); (0,15)*{} **\dir{-};
(15,10); (15,15)*{} **\dir{-};
(7.5,2.2)*{x},
(0,0); (15,0)*{} **\dir{-};
(0,5); (15,5)*{} **\dir{-};
(0,0); (0,5)*{} **\dir{-};
(15,0); (15,5)*{} **\dir{-};
(2,5); (2,10)*{} **\dir{-};
(4,5); (4,10)*{} **\dir{-};
(13,5); (13,10)*{} **\dir{-};
(11,5); (11,10)*{} **\dir{-};
(7.5,7.25)*{\ldots},
(2,0); (2,-3)*{} **\dir{-};
(4,0); (4,-3)*{} **\dir{-};
(13,0); (13,-3)*{} **\dir{-};
(11,0); (11,-3)*{} **\dir{-};
(7.5,-2.25)*{\ldots},
(2,15); (2,18)*{} **\dir{-};
(4,15); (4,18)*{} **\dir{-};
(13,15); (13,18)*{} **\dir{-};
(11,15); (11,18)*{} **\dir{-};
(7.5,17.25)*{\ldots},
\endxy}
\Ea
\ \ \ \ \
 \forall x,y\in \B_n
$$
For example,
$
b''= (b_1)^2.
$

\subsubsection{\bf Remark}\label{5: Remark on B_2} $\B_2$ is the free group on one generator
$b_1$, i.e.\ there is an isomorphism of groups $\Z\rar \B_2$ which sends $n$ into $b_1^n$.
Similarly, $\P\B_2$ is the free group on one generator $b_1^2$, i.e.\ its every element
is of the form $b_1^{2n}$ for some uniquely defined $n\in \Z$. If $\K$ is a field of
characteristic zero, then the prounipotent completions, $\wh{\B}_2(\K)$ and
$\wh{\P\B}_2(\K)$, of these groups coincide with
$\{b_1^{2\mu}\ |\ \mu\in \K\}\simeq \K=\wh{F}_1(\K)$, the prounipotent completion of the
free group in one generator.

\subsection{Algebraic definition of the (pure) braid group}
\label{5: Subsection on generators of (pure) braids} According to M.\ Artin, the braid
group $\B_n$ can be identified with a group of finite type with generators,
\[
b_i:=
\Ba{c}\resizebox{23mm}{!}{\xy
(-16,-10)*{_1},
(-16,10)*{_1},
(-13.2,-10)*{_2},
(-13.2,10)*{_2},
(16.3,-10)*{_n},
(16.3,10)*{_n},
(-10,0)*{...},
(11,0)*{...},
(-4,-10)*{_i},
(-4,10)*{_i},
(4,-10)*{_{i+1}},
(4,10)*{_{i+1}},
\vtwist~{(-4,8)}{(4,8)}{(-4,-8)}{(4,-8)};
(-16,-8)*{}; (-16,8)*{} **\dir{-};
(-14,-8)*{}; (-14,8)*{} **\dir{-};
(-7,-8)*{}; (-7,8)*{} **\dir{-};
(14,-8)*{}; (14,8)*{} **\dir{-};
(7.5,-8)*{}; (7.5,8)*{} **\dir{-};
(16,-8)*{}; (16,8)*{} **\dir{-};
\endxy}\Ea \]
and (so called {\em braid}) relations,
\Beqrn
b_ib_j&=&b_jb_i\ \ \ \mbox{for}\ |i-j|>1,\\
b_ib_{i+1}b_i&=&b_{i+1}b_ib_{i+1}.
\Eeqrn
Hence the prounipotent completion $\wh{\B}_n(\K)$ of $\B_n$ over a field $\K$ can be
identified (see \S {\ref{2: examples of completed filtered}}i) with
$\exp{\wh{\mathfrak{b}}}$, where $\wh{\mathfrak{b}}$ is the quotient of the completed free
Lie algebra $\wh{\ffree}_n$ on letters $\{\ga_i\}_{1\leq i\leq n}$ with respect to the Lie
ideal generated by the following formal power series of Lie words,
\Beqrn
\log(e^{\ga_i}e^{\ga_j}e^{-\ga_i}e^{-\ga_j})&=&0\ \ \ \mbox{for}\ |i-j|>1,\\
\log(e^{\ga_i}e^{\ga_{i+1}}e^{\ga_i}e^{-\ga_{i+1}} e^{-\ga_{i}} e^{-\ga_{i+1}}) &=&0.
\Eeqrn

\mip

Note that the permutation group $\bS_n$ is a group of finite type with generators
$\sigma_i=(i(i+1))$ and relations
\Beqrn
\sigma_i^2&=& \id\\
\sigma_i\sigma_j&=&\sigma_j\sigma_i\ \ \ \mbox{for}\ |i-j|>1,\\
\sigma_i\sigma_{i+1}\sigma_i&=&\sigma_{i+1}\sigma_i\sigma_{i+1}.
\Eeqrn
Hence the projection $p: \B_n\rar \bS_n$ is given on generators by $b_i\rar \sigma_i$.

\mip

According to W.\ Burau and A.\ A.\ Markov, the pure braid group $\P\B_n$ can be identified
with a group of finite type with generators, $\{x_{ij}\}_{1\leq i< j \leq n}$, where
$$
x_{ij}=(b_{j-1}b_{j-2}\cdots b_{i+1})b_i^2 (b_{j-1}b_{j-2}\cdots b_{i+1})^{-1},
$$
and relations,
\Beqrn
x_{ij}x_{kl} &=& x_{kl}x_{ij}\ \ \mbox{for}\ \ \ \ \ \ i<j<k<l\ \ \mbox{and}\ \ i<k<l<j,\\
x_{ij}x_{ik}x_{jk} &=& x_{ik}x_{jk}x_{ij} \ \ \ \mbox{for}\ \ i<j<k\\
x_{ik}x_{jk}x_{ij} &=& x_{jk}x_{ij}x_{ik} \ \ \ \mbox{for}\ \ i<j<k\\
x_{ik}x_{jk}x_{jl}x_{jk}^{-1} &=& x_{jk}x_{jl}x_{jk}^{-1}x_{ik} \ \ \mbox{for}\ \ i<j<k<l.
\Eeqrn
Hence its prounipotent completion $\wh{\P\B}_n$ can be described via the construction in
\S{\ref{2: examples of completed filtered}}(i)
(see Exercise {\ref{5: Exercise on graded of pb}} below).

\subsubsection{\bf Exercise}\label{5: Exercise on b_2b_1 squared b_2 inverse} (i)
Show that $(b_1b_2)^3=(b_2b_1)^3$.

\sip

(ii) Show that the element $(b_1b_2)^3$ lies in $\P\B_3\subset \B_3$ and is central both in $\P\B_3$ and in  $\B_3$.

\sip

(iii) Show that
$$
b_2 b_1^2 b_2^{-1}= (b_1b_2)^3 b_1^{-2}b_2^{-2}, \ \ \ (b_2b_1)b_2^2(b_2b_1)^{-1}= b_1^2,\
\ \ (b_1b_2b_1)^2=(b_2b_1b_2)^2=(b_1b_2)^3
$$
in $\B_3$.

\subsection{Pure braids and semidirect products of free groups} For $n\geq 1$, the
``forgetful" map
$$
\Ba{rccc}
f : & C_{n+1}(\C) & \lon & C_{n}(\C)\\
& (z_1, \ldots, z_{n+1}) & \lon & (z_1, \ldots , z_n)
\Ea
$$
is a locally trivial fiber bundle \cite{FV}. Its fiber is isomorphic to
$\C\setminus \{z_1, \ldots , z_n\}$ whose first homotopy group is the free group $F_n$ on $n$ generators
while the second homotopy group is trivial, $\pi_2(\C\setminus \{z_1, \ldots , z_n\}, *)=1$. We can choose base points $p_n$ in $C_n(\C)$, $n\geq 1$, such that $f(p_{n+1})=p_n$. Hence the associated long exact sequence of homotopy groups reads,
$$
1\lon \pi_2(C_{n+1}(\C),p_{n+1}) \lon  \pi_2(C_{n}(\C),p_n)  \lon F_n \lon
\pi_1(C_{n+1}(\C),p_{n+1}) \lon
 \pi_1(C_{n}(\C),p_n) \lon 1,\ \ \ \ n\geq 1.
$$
As $\pi_2(C_{2}(\C),p_2)=\pi_2(\C,p_1)=1$, one obtains by induction that
$\pi_2(C_{n}(\C),p_n)=1$ for all $n$, and  hence one gets an exact sequence of groups,
$$
1\lon F_n\lon \P\B_{n+1}\lon \P\B_{n} \lon 1.
$$
This exact sequence splits as the fiber bundle $f : C_{n+1}(\C)\rar C_{n}(\C)$ admits a
cross-section,
$$
\Ba{rccc}
s : & C_{n}(\C) & \lon & C_{n+1}(\C)\\
 & (z_1, \ldots , z_n) &\lon &  (z_1, \ldots , z_n, |z_1|+\ldots + |z_n|+1).
\Ea
$$
Hence we proved the following

\subsubsection{\bf Theorem} {\em For any $n\geq 1$ one has}
\Beqrn
\P\B_{n+1} &= &F_n\rtimes \P\B_n \\
&=& F_n\rtimes \left(F_{n-1} \rtimes \ldots \rtimes \left( F_3\rtimes
\left( F_2 \rtimes F_1\right)\right)\right)
\Eeqrn

\subsubsection{\bf Lemma}\label{5: Lemma on PB_3}{\em Every element of $\P\B_3$ can be
uniquely represented as $f(b_1^2,b_2^2) (b_1b_2)^{3n}$ for some  $n\in \Z$ and some
element $f(x,y)$ in the free group $F_2$ on the generators $x$ and $y$.}
\begin{proof}
The group  $\P\B_3$ is generated by $x_{12}=b_1^2$, $x_{23}=b_2^2$ and
$x_{13}=b_2b_1^2b_2^{-1}$, or, equivalently, by $b_1^2$, $b_2^2$ and
$$
x_{13}x_{23}x_{12}=b_2b_1^2b_2 b_1^2=(b_2b_1)^3=(b_1b_2)^3,
$$
where we used the braid relations $b_1b_2b_1=b_2b_1b_2$ in $\B_3$. As
$(b_1b_2)^3$ is central in $\P\B_3$ and  there are no relations between the generators
$b_1^2$ and $b_2^2$, the result follows.
\end{proof}

\subsubsection{\bf Corollary}\label{5: Corollary on PB_3}  $\P\B_3=F_2\times \Z$.
\mip

This corollary can be seen also from the fact that the fibration $\C_3(\C)\rar \C_2(\C)$
is trivial.

\mip

If $\K$ is a field of characteristic zero, then $\wh{\P\B}_3(\K)= \wh{F_2}(\K)\times \K$,
with $\wh{F_2}(\K)$  generated by $[b_1^{2}]^{\la_1}$, $[b_2^{2}]^{\la_2}$, and the factor
$\K$ corresponding to  $[(b_1b_2)^{3}]^{\la_3}$, $\la_i\in \K$,
$i=1,2,3$. On the contrary, $\wh{\B}_3(\K)$ is a much simpler group.

\subsubsection{\bf Proposition} {\em If $char(\K)=0$, then}\, $
\wh{B}_3(\K)\simeq \wh{F}_1(\K)\simeq \K$.

\begin{proof} Set $c:= (b_1b_2)^3$. As it is a central element in $\B_3$,
$c^{\la}$ commutes with any element in $\wh{\B}_3(\K)$ for any $\la\in \K$. Using
the last equality in \S{\ref{5: Exercise on b_2b_1 squared b_2 inverse}}(iii),
we obtain
$$
b_1b_2b_1=b_2b_1b_2=c^{\frac{1}{2}}\ \ \ \mbox{in}\ \ \ \wh{\B}_3(\K)
$$
so that $b_1b_2=b_2^{-1}c^{\frac{1}{2}}$, $b_2^{-1}c^{\frac{1}{2}}b_1=c^{\frac{1}{2}}$,
and hence
$
b_1=b_2=c^{\frac{1}{6}} \ \ \mbox{in}\ \ \ \wh{\B}_3(\K).
$
\end{proof}

\subsubsection{\bf Corollary}\label{5: corollary on sfb_1 sfb_2}{\em Let $char(\K)=0$.
The injection $\P\B_3\rar \B_3$
induces a surjection $\wh{\P\B}_3(\K)\rar \wh{\B}_3(\K)$ whose kernel is the prounipotent
completion of the free group generated by
$\sfb_1= c^{-1/3}b_1^2$ and $\sfb_2:= c^{-1/3}b_2^2$.}

\subsection{Non-$\bS$ operad of pure braids}\label{5: subsect on non sigma operad of pure braids}
Consider a collection of groups,
$$
\cP\cB:=\{\pi_1(C_n(\C))\simeq \P\B_n\}.
$$
Recall that $\pi_1(C_n(\C))$ stands for the isomorphism class of the fundamental group $\pi_1(C_n(\C), p_n)$
defined for a particular choice of the base point $p_n$; as $C_n(\C)$ is path connected, the  class  $\pi_1(C_n(\C))$
does not depend on such a choice.
The collection $\cP\cB$  can be made into a non-$\bS$ operad with respect to the operadic
composition,
$$
\Ba{rccc}
\circ_i :& \P\B_n \times \P\B_n &\lon & \P\B_{n+m-1}\ \ \ \ \  \forall n,m\in \N,i\in [n]\\
      &    (b',b'') & \lon b'\circ_i b''
      \Ea
$$
which is
defined by
replacing the $i$-labelled strand in $b'$  by the braid $b''$ made very thin.

\subsubsection{\bf Cosimplicial structure on $\cP\cB$} Let $e$ be the identity element in the group
$\cP\cB(2)=\P\B_2$. Then it satisfies the condition
$$
e\circ_1 e=e \circ_2 e
$$
and hence  gives rise to a pre-cosimplicial structure,
$$
d_i: \P\B_n \lon \P\B_{n+1}, \ \ \ \ i=0,1,\ldots, n+1,
$$
 on the collection $\cP\cB=\{\P\B_n\}_{n\geq 1}$
given by the explicit formulae in (\ref{4: formulae for d_i}). For example, one has in $\cP\cB_3$,
\Beq\label{5: d_i on pure braids}
d_0(b_1^2)= (b_2)^2, \ \ \ d_1(b_1^2)=(b_2b_1)^2,\ \ \   d_2(b_1^2)=(b_1b_2)^2,\ \ \
d_3(b_1^2)= (b_1)^2.
\Eeq
  In fact, this structure can be completed to a {\em cosimplicial}\, structure on $\cP\cB$ by defining the operators
 $$
 s_i: \P\B_{n+1} \lon \P\B_{n}, \ \ \ \ \ i=1,2,\ldots,n+1,
 $$
 where $s_i$ erases the strand labelled by $i$.
\subsection{An operad of parenthesized braids} If we were able to make the operad of little disks
$\caD$ into an operad in the category of {\em based}\, topological spaces, then, by applying
the monoidal functor (\ref{3: fundam group monoidal functor}), we would have obtained an operad $\pi_1(\caD)$
in the category of groups. However this is impossible as there is {\em no}\, $\bS_n$ invariant
configuration of little disks in $\caD(n)$ for every $n\geq 2$.

\sip

The fundamental groupoid functor
(\ref{3: fundam groupoid monoidal functor}) applied to $\caD$  gives us an operad in the
category of groupoids, $\mathsf{CatG}$,  with nontrivial action of permutations groups
as endofunctors, $\sigma: \pi(\caD)(n) \rar  \pi(\caD)(n)$, $\sigma\in \bS_n$, $n\geq 1$.
However this operad is too big for our purposes, the category $\pi(\caD)(n)$ has
too many objects (which are points in $\caD(n)$). Let us construct by induction a
 ``smallest possible" non-trivial suboperad, $\cP a \cB$,  of the operad
 $\pi(\caD)$ such that $\cP a\cB(n)$ is a {\em full}\, subcategory of $\pi(\caD)(n)$
 for each $n\geq 1$:
 \Bi
 \item[$n=1$]: Set $\cP a\cB(1)$ to a full subcategory of $\pi(\caD)(1)$ which has only
 one object, the configuration $D_{0,1}\in \pi(\caD)(1)$ (which is the unit in the operad
 $\caD$).
 \item[$n=2$]: Let $\mathbf{1}$ be any configuration of little disks in $\caD(2)$, that is, a pair of disjoint disks inside the unit disk labeled by 1 and 2, and
 set $\mathbf{2}:=\sigma(\mathbf 1)$, where $\sigma=(12)\in \bS_2$. Set   $\cP a\cB(2)$
 be the full subcategory of $\pi(\caD)(2)$ whose only objects are
 $\mathbf{1}$ and $\mathbf{2}$; the $\bS_2$ action on  $\pi(\caD)(2)$
 leaves the subcategory $\cP a\cB(2)$ invariant. It is useful to identify
 $Ob(\cP a\cB(2))$ with  $\cP a(2)$, the set of two planar binary corollas,
 $$
 \mathbf{1}=\xy
(0,0)*{\circ}="a",
(0,3)*{}="0",
(-3,-4.5)*{_1},
(3,-4.5)*{_2},
(-3,-3)*{}="b_1",
(3,-3)*{}="b_2",
\ar @{-} "a";"0" <0pt>
\ar @{-} "a";"b_1" <0pt>
\ar @{-} "a";"b_2" <0pt>
 \endxy
 \ \ \ \ \ \ \ \
 \mathbf{2}=\xy
(0,0)*{\circ}="a",
(0,3)*{}="0",
(-3,-4.5)*{_2},
(3,-4.5)*{_1},
(-3,-3)*{}="b_1",
(3,-3)*{}="b_2",
\ar @{-} "a";"0" <0pt>
\ar @{-} "a";"b_1" <0pt>
\ar @{-} "a";"b_2" <0pt>
 \endxy
 $$

 \item[$n=3$]: Let $\circ_i: \caD(2)\times \caD(2)\rar \caD(3)$, $i=\{1,2\}$,
  be the operadic
 compositions in $\caD(n)$. Set   $\cP a\cB(3)$
 be the full subcategory of $\pi(\caD)(3)$ whose only objects are
 $\mathbf{p}\circ_i \mathbf{q}$, where $\mathbf{p},\mathbf{q}\in \{\mathbf{1},\mathbf{2}\}$. It is clear
 that the set $Ob( \cP a\cB(3))$ can be identified with $\cP a(3)$, the set
 of parenthesized permutations of $\mathbf[3]=\{ \mathbf{1}, \mathbf{2}, \mathbf{3}\}$
 (see \S {\ref{4: operad Pa}}), or, equivalently,  the set of $[\mathbf 3]$-labelled
 planar binary trees,
 $$
 \mathbf{1}\circ_1 \mathbf{1}= (\mathbf{1}\mathbf{2})\mathbf{3}, \ \ \
 \mathbf{1}\circ_2 \mathbf{1}= \mathbf{1}(\mathbf{2}\mathbf{3}),\ \ \
 \mathbf{1}\circ_1 \mathbf{2}= (\mathbf{2}\mathbf{1})\mathbf{3}, \ \ \
 \mathbf{2}\circ_2 \mathbf{2}= \mathbf{3}(\mathbf{2}\mathbf{1}),\ \ \  etc.
 $$
 \item[$n\geq 3$]: Let $Ob(\cP a\cB)$ be the suboperad of $Ob(\pi(\caD))$
 (in the category of sets) generated by the $\cS$-module $E=\{E(n)\}$, where
 $$
E(n):=\left\{\Ba{ll} Ob(\cP a\cB(1))=D_{0,1} & \mbox{for}\ n=1,\\
Ob(\cP a \cB(2))=\{\mathbf{1},\mathbf{2}\} & \mbox{for}\ n=2,\\
\emptyset & \mbox{for}\ n\geq 2
\Ea
\right.
$$
It can be identified with the operad $\cP a$.
 \Ei
 Let now $\cP a \cB(n)$ be the full subcategory of the category $\pi(\caD(n))$.
 By construction, the collection of groupoids,
$$
\cP a \cB:=\{\cP a \cB(n)\}_{n\geq 1},
$$
is a suboperad of $\pi(\caD)$. It was first introduced by D.\ Bar-Natan \cite{BN}.
The idea to use the language of operads in D.\ Bar-Natan's approach to Drinfeld associators
and the Grothendieck-Teichm\"uller group is due to Tamarkin \cite{Ta2}
(see also the works  of P.\ Severa and T.\ Willwacher  \cite{Se,SW}).
This idea makes D.\ Bar-Natan's story short and transparent.

\mip

One can equivalently describe the  underlying S-module of the operad $\cP a \cB$ as functor,
$$
\Ba{rccc}
\cP a \cB: & \sS &\lon  & \mathsf{CatG}\\
           &  I  & \lon & \cP a \cB(I)= \left\{\Ba{c} \mbox{objects},\ \ Ob(\cP a \cB(I)) \\
           \mbox{sets of morphisms},\ \{ \Mor(A_I,B_I)\}_ {A_I, B_I\in Ob(\cP a \cB(I))}
           \Ea
           \right.
\Ea
$$
where $\cP a \cB(I)$ is the category given by the following data:
\Bi
\item[(i)] the objects  are  parenthesized permutations of the set $I$, or,
equivalently, planar $I$-labelled binary trees equipped with the operadic structure explained in  \S{\ref{4: operad Pa}};
put another way, the ``object" part of the operadic structure is given by the functor
$$
\Ba{rccc}
 \cP a: & \sS &\lon  & \mathsf{Set}\\
   & I  & \lon & \cP a(I)=:Ob(\cP a \cB(I))
\Ea
$$

\item[(ii)] for any two objects $A_I$ and $B_I$ in  $Ob(\cP a \cB(I))$, that is, for any two parenthesized permutations
of the set $I$, the associated set of morphisms,  $\Mor(A(I),B(I))$, is defined to be
the set of braids whose strands  connect the same elements
of $I$ (and hence have a non-ambiguous  $I$-labelling). For example,
$$
\beta_{1,2}:=
\Ba{c}\resizebox{7mm}{!}{\xy
(-5,-10)*{_{1}},
(-5,10)*{_{2}},
(5,-10)*{_{2}},
(5,10)*{_{1}},
\vtwist~{(-5,8)}{(5,8)}{(-5,-8)}{(5,-8)};
\endxy}\Ea \in \Mor\left(12, 21\right),
\ \ \ \ \ \id_{1,2,3}:=
\Ba{c}\resizebox{12mm}{!}{\xy
(-8,-10)*{_{(1}},
(-7,10)*{_{1}},
(0,-10)*{_{2)}},
(0,10)*{_{(2}},
(7,-10)*{_{3}},
(8,10)*{_{3)}},
(-7,-8)*{}; (-7,8)*{} **\dir{-};
(0,-8)*{}; (0,8)*{} **\dir{-};
(7,-8)*{}; (7,8)*{} **\dir{-};
\endxy}\Ea\in \Mor\left((12)3), 1(23)\right)\ \ \ \
$$
The operadic composition $b'\circ_i b''$ on morphisms
is defined by
replacing the $i$-labelled strand of the braid $b'$ with the braid $b''$ made thin.
\Ei


\subsubsection{\bf Exercise} Show that the operad $\cP a \cB$ is generated (via operadic and categorical compositions) by the elements $\beta_{1,2}$, $\id_{1,2,3}$ and their inverses.

Find an explicit expression for the element
$\Ba{c}\resizebox{11mm}{!}{\xy
(-7,-8)*{_{1}},
(-7.4,8)*{_{(2}},
(1,-8)*{_{(2}},
(-2,8)*{_{3)}},
(7,-8)*{_{3)}},
(6,8)*{_{1}},
(5,6)*{}; (-6,-6)*{} **\crv{(4,0)&(-4,3)}
\POS?(.4)*+{}="x" \POS?(.6)*+{}="y";
(-2,6)*{}; "x" **\crv{};
(-6,6)*{}; "y" **\crv{};
"y"+(.5,-1); (.5,-6)*{} **\dir{-};
"x"+(.5,-1); (6,-6)*{} **\dir{-};
\endxy}\Ea\in \Mor\left(1(23), (23)1\right)$
in terms of the generators $\be_{\bu\bu }$ and  $\id_{\bu,\bu,\bu}$.

\subsubsection{\bf A family of functors on $\cP a \cB$}\label{5: cosimiplic struc on PaB} Let $\id_{12}$ be the identity element in $\Mor(12, 12)$,
$$
 \id_{12}:=
\Ba{c}\resizebox{4mm}{!}{\xy
(-2,8)*{_{1}},
(-2,-8)*{_{1}},
(2,-8)*{_{2}},
(2,8)*{_{2}},
(-2,-6)*{}; (-2,6)*{} **\dir{-};
(2,-6)*{}; (2,6)*{} **\dir{-};
\endxy}\Ea
$$
and, for any $n\in \N_{\geq 1}$, define a family of functors,
$$
d_i: \cP a \cB(n) \lon \cP a \cB(n+1), \ \ \ \ i=0,1,\ldots, n, n+1,
$$
as follows:
\Bi
\item[(i)] on objects (i.e.\ on parenthesized permutations of $[n]$), one has
$$
\Ba{rccc}
d_i: &  Ob(\cP a \cB(n)) & \lon &  Ob(\cP a \cB(n+1))\\
     &  A_n & \lon & d_i A_n:=\left\{\Ba{cc} (12)\circ_2 A_n & \mbox{for}\ i=0\\
     A_n\circ_i (12) & \mbox{for}\ i\in [n]\\
     (12)\circ_1 A_n & \mbox{for}\ i=n+1\\
     \Ea\right.
\Ea
$$
\item[(ii)] on morphisms (i.e.\ on braids connecting the same elements in a pair of parenthesized permutations), one has
$$
\Ba{rccc}
d_i: &  \Mor(A_n,B_n) & \lon &  \Mor(d_i A_n, d_iB_n)\\
     &  f_n & \lon & d_i f_n:=\left\{\Ba{cc} \id_{1,2}\circ_2 f_n & \mbox{for}\ i=0\\
     f_n\circ_i \id_{1,2} &  \mbox{for}\ i\in [n]\\
     \id_{1,2}\circ_1 f_n & \mbox{for}\ i=n+1\\
     \Ea\right.
\Ea
$$
\Ei
Note that in the operad $Ob(\cP a \cB)$,
$$
(12)\circ_1 (12)\neq (12)\circ_2 (12)
$$
so that the functors $d_i$ do {\em not}\,
make the collection of groupoids $\{\cP a\cB(n)\}_{n\geq 1}$ into
a (pre)cosimplicial object in the category of small categories.
However these functors will be very useful for us below for the the following reason.
Let $\f$ be an operad in the category of small categories with same set
of objects as in $\cP a \cB$. Let $\id_{1,2}$ be the identity morphism from $(12)$ to $(12)$ in the category $\f(2)$. Then the formulae formally identical to the ones above give
us a family of functors, $d_i:\f(n)\rar \f(n+1)$; moreover, any morphism
of operads $F:\cP a \cB\rar \f$ {\em respects}\, these functors, i.e.
$$
d_i\circ F(n)= F(n+1)\circ d_i,
$$
for any $n\in \N$ and any $i\in \{0,1,\ldots n+1\}$.

\sip

\subsubsection{\bf Remark} The non-$\bS$ operad $\cP\cB$ also has a presimplicial structure given formally by the same formulae as above.

\subsubsection{\bf Pentagon equation} Consider the elements
$$
 d_0\id_{1,2,3}=\Ba{c}\resizebox{12mm}{!}{ \xy
(0,6)*{_{1}},
(0,-6)*{_{1}},
(4,-6)*{_{((2}},
(4,6)*{_{(2}},
(8.5,-6)*{_{3)}},
(8.0,6)*{_{((3}},
(12.7,-6)*{_{4)}},
(12.8,6)*{_{4))}},
(0,-4)*{}; (0,4)*{} **\dir{-};
(4,-4)*{}; (4,4)*{} **\dir{-};
(8,-4)*{}; (8,4)*{} **\dir{-};
(12,-4)*{}; (12,4)*{} **\dir{-};
\endxy}\Ea, \ \ \ \
d_1\id_{1,2,3}=
\Ba{c}\resizebox{12mm}{!}{\xy
(0,6)*{_{(1}},
(0,-6)*{_{((1}},
(4.4,-6)*{_{2)}},
(4.4,6)*{_{2)}},
(8.5,-6)*{_{3)}},
(8.0,6)*{_{(3}},
(12.7,-6)*{_{4}},
(12.8,6)*{_{4)}},
(0,-4)*{}; (0,4)*{} **\dir{-};
(4,-4)*{}; (4,4)*{} **\dir{-};
(8,-4)*{}; (8,4)*{} **\dir{-};
(12,-4)*{}; (12,4)*{} **\dir{-};
\endxy}\Ea, \ \ \ \
 d_2\id_{1,2,3}=
 \Ba{c}\resizebox{12mm}{!}{\xy
(0,6)*{_{1}},
(0,-6)*{_{(1}},
(4.0,-6)*{_{(2}},
(4.0,6)*{_{((2}},
(8.8,-6)*{_{3))}},
(8.5,6)*{_{3)}},
(12.7,-6)*{_{4}},
(12.8,6)*{_{4)}},
(0,-4)*{}; (0,4)*{} **\dir{-};
(4,-4)*{}; (4,4)*{} **\dir{-};
(8,-4)*{}; (8,4)*{} **\dir{-};
(12,-4)*{}; (12,4)*{} **\dir{-};
\endxy}\Ea,
$$

$$
 d_3\id_{1,2,3}=\Ba{c}\resizebox{12mm}{!}{\xy
(0,6)*{_{1}},
(0,-6)*{_{(1}},
(4.4,-6)*{_{2)}},
(3.8,6)*{_{((2}},
(8.5,-6)*{_{(3}},
(8.2,6)*{_{3)}},
(12.7,-6)*{_{4)}},
(12.8,6)*{_{4)}},
(0,-4)*{}; (0,4)*{} **\dir{-};
(4,-4)*{}; (4,4)*{} **\dir{-};
(8,-4)*{}; (8,4)*{} **\dir{-};
(12,-4)*{}; (12,4)*{} **\dir{-};
\endxy}\Ea, \
\ \ \ \
 d_4\id_{1,2,3}=  \Ba{c}\resizebox{12mm}{!}{\xy
(0,6)*{_{(1}},
(0,-6)*{_{((1}},
(4.4,-6)*{_{2)}},
(4.0,6)*{_{(2}},
(8.8,-6)*{_{3)}},
(8.5,6)*{_{3)}},
(12.7,-6)*{_{4}},
(12.8,6)*{_{4}},
(0,-4)*{}; (0,4)*{} **\dir{-};
(4,-4)*{}; (4,4)*{} **\dir{-};
(8,-4)*{}; (8,4)*{} **\dir{-};
(12,-4)*{}; (12,4)*{} **\dir{-};
\endxy}\Ea,
$$
and check that they satisfy the pentagon equation (cf.\ (\ref{3: Pentagon}))
\Beq\label{5: penragon in PaB}
d_1\id_{1,2,3} \circ d_3\id_{1,2,3}=d_4\id_{1,2,3}\circ d_2\id_{1,2,3}\circ d_0\id_{1,2,3}
\Eeq
Here $\circ$ denotes the categorical (rather than operadic) composition of morphisms.

\sip

Following a long tradition (see, e.g., \cite{Dr,ES}), the elements  $d_0\id_{1,2,3}$,
 $d_1\id_{1,2,3}$,  $d_1\id_{1,2,3}$,  $d_3\id_{1,2,3}$ and  $d_4\id_{1,2,3}$ should have been
 denoted, respectively, by $\id_{2,3,4}$, $\id_{12,3,4}$, $\id_{1,23,4}$, $\id_{1,2,34}$
 and  $\id_{1,2,3}$ so that the pentagon equation takes a familiar form (cf.\ equation (17.2) in the book \cite{ES})
 $$
  \id_{12,3,4}\circ\id_{1,23,4} = \id_{1,2,3}\circ \id_{1,23,4} \circ  \id_{2,3,4}.
 $$
However, we shall mostly stick to the cosimplicial notation in this survey.
\subsubsection{\bf Hexagon equations} One has in $\cP a \cB(3)$,
$$
d_0\beta_{1,2}= \Ba{c}\resizebox{11mm}{!}{\xy
(-5,-10)*{_{(2}},
(-5,10)*{_{(3}},
(5,-10)*{_{3)}},
(5,10)*{_{2)}},
(-10,-10)*{_{1}},
(-10,10)*{_{1}},
\vtwist~{(-5,8)}{(5,8)}{(-5,-8)}{(5,-8)};
(-10,-8)*{}; (-10,8)*{} **\dir{-};
\endxy}\Ea\ \ \ , \ \ \ \
d_3\beta_{1,2}= \Ba{c}\resizebox{11mm}{!}{
\xy
(-5,-10)*{_{(1}},
(-5,10)*{_{(2}},
(5,-10)*{_{2)}},
(5,10)*{_{1)}},
(10,-10)*{_{3}},
(10,10)*{_{3}},
\vtwist~{(-5,8)}{(5,8)}{(-5,-8)}{(5,-8)};
(10,-8)*{}; (10,8)*{} **\dir{-};
\endxy}\Ea\ \ \ , \ \ \ \ \ \
d_1 \beta_{1,2}=\Ba{c}\resizebox{13mm}{!}{
\xy
(-10,-10)*{_{(1}},
(-5,-10)*{_{2)}},
(9,-10)*{_{3}},
(-10,10)*{_{3}},
(5,10)*{_{(1}},
(10,10)*{_{2)}},
(5,8)*{}; (-10,-8)*{} **\crv{(6,0)&(-10,3)}
\POS?(.4)*+{}="x" \POS?(.5)*+{}="y";
(9,8)*{}; (-6,-8)*{} **\crv{(7,0)&(-8,3)}
\POS?(.4)*+{}="u" \POS?(.54)*+{}="v";
(-10,8)*{}; "y" **\crv{};
"v"+(.5,-1); (8,-8)*{} **\dir{-};
\endxy}\Ea$$
The permutation group $\bS_n$ acts on $\cP a \cB(n)$ by relabelling of strands, for example,
$$
\mbox{(\scriptsize 23)}\ccdot \id_{1,2,3}= \id_{1,3,2}=
\Ba{c}\resizebox{11mm}{!}{
\xy
(-6,-10)*{_{(1}},
(-5,10)*{_{1}},
(0,-10)*{_{3)}},
(0,10)*{_{(3}},
(5,-10)*{_{2}},
(6,10)*{_{2)}},
(-5,-8)*{}; (-5,8)*{} **\dir{-};
(0,-8)*{}; (0,8)*{} **\dir{-};
(5,-8)*{}; (5,8)*{} **\dir{-};
\endxy}\Ea, \ \ \
\mbox{(\scriptsize 123)}\ccdot \id_{1,2,3}= \id_{3,1,2}=
\Ba{c}\resizebox{11mm}{!}{
\xy
(-6,-10)*{_{(3}},
(-5,10)*{_{3}},
(0,-10)*{_{1)}},
(0,10)*{_{(1}},
(5,-10)*{_{2}},
(6,10)*{_{2)}},
(-6,-8)*{}; (-6,8)*{} **\dir{-};
(0,-8)*{}; (0,8)*{} **\dir{-};
(5,-8)*{}; (5,8)*{} **\dir{-};
\endxy}\Ea, \ \ \ \
\mbox{(\scriptsize 23)}\ccdot (d_3\beta_{1,2})=
\Ba{c}\resizebox{13mm}{!}{
 \xy
(-5,-10)*{_{(1}},
(-5,10)*{_{(3}},
(5,-10)*{_{3)}},
(5,10)*{_{1)}},
(10,-10)*{_{2}},
(10,10)*{_{2}},
\vtwist~{(-5,8)}{(5,8)}{(-5,-8)}{(5,-8)};
(10,-8)*{}; (10,8)*{} **\dir{-};
\endxy}\Ea
$$
It is easy to see that the generators satisfy the following hexagon equation in $\cP a \cB (3)$,
\Beq\label{5: first hehagon in PaB(3)}
d_1 \beta_{1,2}= \id_{1,2,3}\  \circ\  d_0\beta_{1,2}\ \circ\ \left(\mbox{(\scriptsize 23)}\ccdot \id_{1,2,3}^{-1}\right)\ \circ \
\mbox{(\scriptsize 23)} \ccdot d_3\beta_{1,2}\  \circ\ \mbox{(\scriptsize 123)}\ccdot \id_{1,2,3}.
\Eeq
Studying analogously
$$
\be^{-1}_{2,1}=\mbox{(\scriptsize 12)}\ccdot \be^{-1}_{1,2}= = \left( \mbox{(\scriptsize 12)}\cdot \beta_{1,2}\right)^{-1}=
\Ba{c}\resizebox{8mm}{!}{
\xy
(-5,-10)*{_{1}},
(-5,10)*{_{2}},
(5,-10)*{_{2}},
(5,10)*{_{1}},
\vcross~{(-5,8)}{(5,8)}{(-5,-8)}{(5,-8)};
\endxy}\Ea \in \Mor\left(12, 21\right)
$$
we get a second hexagon equation,
\Beq\label{5: second hehagon in PaB(3)}
d_1 \be_{2,1}^{-1}= \id_{1,2,3}\  \circ\  d_0\be^{-1}_{2,1}\ \circ\ \left(\mbox{(\scriptsize 23)}\ccdot \id_{1,2,3}^{-1}\right)\ \circ \
\mbox{(\scriptsize 23)} \ccdot d_3\be^{-1}_{2,1}\  \circ\ \mbox{(\scriptsize 123)}\ccdot \id_{1,2,3}.
\Eeq
for the generators $\beta_{1,2}$ and $\id_{1,2,3}$ of the operad $\cP a \cB$.

\subsubsection{\bf Theorem}\label{5: Theorem on morhisms from PaB} {\em Let $\f$ be an operad in the small category of groupoids with same set
of objects as in $\cP a \cB$.
    There is a one-to-one correspondence between morphisms of operads,
$$
F: \cP a \cB \lon \f,
$$
which are identical on objects, and elements
$F(\beta_{1,2})\in \Mor_\f(12,21)$ and {\em $F(\id_{1,2,3})\in \Mor_\f((12)3,1(23))$}, which satisfy
the pentagon equation,}
\Beq\label{5: pentagon in PaB--f}
d_1F(\id_{1,2,3}) \circ d_3F(\id_{1,2,3})
=d_4F(\id_{1,2,3})\circ d_2F(\id_{1,2,3})\circ d_0F(\id_{1,2,3})
\Eeq
{\em and the hexagon equations},
\Beq\label{5: first hexagon in PaB--f}
d_1 F(\beta_{1,2})= F(\id_{1,2,3})\  \circ\  d_0F(\beta_{1,2})\ \circ\ \left(\mbox{\scriptsize (23)}\cdot F(\id_{1,2,3})^{-1}\right)\ \circ \
\mbox{\scriptsize (23)} \ccdot d_3F(\beta_{1,2})\  \circ\ \mbox{\scriptsize (321)}\ccdot F(\id_{1,2,3}),
\Eeq
\Beq\label{5: second hehagon in PaB--f}
d_1 F(\be_{2,1})^{-1}= F(\id_{1,2,3})\  \circ\  d_0(F(\be_{2,1})^{-1})\ \circ\ \left(\mbox{\scriptsize (23)}\ccdot F(\id_{1,2,3})^{-1}\right)\ \circ \
\mbox{\scriptsize (23)} \ccdot d_3(F(\be_{2,1})^{-1})\  \circ\ \mbox{\scriptsize (321)}\ccdot F(\id_{1,2,3}).
\Eeq
\begin{proof}
As $\cP a \cB$ is generated by $\beta_{1,2}$ and $\id_{1,2,3}$, any morphism of operads
$F: \cP a \cB\rar \f$ is uniquely determined by its values on the generators which must satisfy
the equations obtained by applying $F$ to (\ref{5: penragon in PaB}), (\ref{5: first hehagon in PaB(3)}) and (\ref{5: second hehagon in PaB(3)}).

\sip

The reverse statement follows from the Mac Lane coherence theorem for braided monoidal categories, i.e.\ from the fact that all the relations between  iterations (via operadic and categorial compositions) of the generators follow from  (\ref{5: penragon in PaB}), (\ref{5: first hehagon in PaB(3)}) and (\ref{5: second hehagon in PaB(3)}).
\end{proof}



\subsection{D.\ Bar-Natan's operad $\wh{\cP a \cB}_\K$}\label{5: subsection on wh{PaB}_K}
Applying the functor $\Delta_\K$ (see \S{\ref{3: tensor functors}}(iii)) to the operad $\cP a \cB$
we get an operad,
$$
\cP a \cB_\K:=\Delta_\K(\cP a \cB),
$$
 in the symmetric monoidal category $\mathsf{Cat(coAss_\K)}$. For any finite set $I$ the objects, $A_I$, of the category
 $\cP a \cB_\K(I)$ are the same as objects of $\cP a \cB$, that is, the parenthesized permutations of the set $I$, while morphisms are formal $\K$-linear combinations of braids,
 $$
\Mor_{\cP a \cB_\K}(A_I,B_I)= \left\{\sum_{i} \la_i f_i\ | \la_i\in \C, f_i\in \Mor_{\cP a \cB}(A_I,B_I)\right\},
 $$
 which connect the same letters in the words $A_I$ and $B_I$. The coproduct structure on morphisms is uniquely defined by the condition that all $f_i$ are group-like.

 \sip

 Every operad $\f$  in $\mathsf{Cat(coAss_\K)}$, in particular  $\cP a \cB_\K$,   comes equipped with a {\em diagonal}\, morphism of operads,
\Beq\label{6: diagonal in Cat(coAss)}
\Delta: \f \lon \f \ot \f,
\Eeq
which  is given on objects $A_I\in Ob(\f(I))$ by
$$
\Delta(A_I)=A_I\times A_I,
$$
and on morphisms $\Delta$ is given by the above coproduct in $\Mor_{\f}(A_I,B_I)$.

\mip

Let us define one more operad in this category,
$$
\cP a_\K: \sS \lon \mathsf{Cat(coAss_\K)},$$
whose objects are, by definition, identical to objects of the operad $\cP a$ (and hence of  $\cP a \cB_\K$)
and whose morphisms spaces, $\Mor_{\cP a_\K}(A_I,B_I)$, are identified with $\K$ equipped with the coproduct $\Delta(1):=1\ot 1$. Speaking plainly, $\cP a_\K$ is the operad $\cP a\cB_\K$
with ``braiding" on morphisms forgotten.

\mip

There is a augmentation morphism of operads (cf.\ \S{\ref{2: Quillen's construction}}),
$$
\var: \cP a \cB_\K \lon \cP a_\K,
$$
which is identical on objects, and is given on morphisms by the formula,
$$
\Ba{rccc}
\var: & \Mor_{\cP a \cB_\K}(A_I,B_I) &\lon & \Mor_{\cP a_\K}(A_I,B_I)\\
     &  \sum_{i} \la_i f_i &\lon & \sum_{i} \la_i
\Ea
$$
The kernel, $\cI_\K:= \Ker\, \var$ is an operadic  ideal in $\cP a \cB_\K$.  One can define its power, $\cI^m_\K$, $m\geq 1$, which is a suboperad of  $\cI_\K$ with the same objects as
in $\cI_\K$ but with morphisms given by the categorical composition   of at least $m$ morphisms from $\cI_\K$.
 In this way we get a family of ideals, $\{\cI^m_\K\}_{m\geq 1}$, of the operad $\cP a \cB_\K$, and hence an associated inverse system of operads in the category  $\mathsf{Cat(coAss_\K)}$,
$$
\ldots
\lon
\cP a\cB_\K^{(m)}
\lon
\cP a\cB_\K^{(m-1)}
\lon
\ldots
\lon
\cP a\cB_\K^{(1)}
\lon \cP a \cB_\K^{(0)},
$$
where $\cP a\cB_\K^{(m)}$ is the quotient operad $\cP a \cB_\K/{\cI_\K^{m+1}}$, $m\geq 0$.
The inverse limit,
$$
\wh{\cP a \cB}_\K:= \lim_{\longleftarrow} \cP a\cB_\K^{(m)} .
$$
is an operad in the category $\wh{\mathsf{Cat}}(\mathsf{coAss}_\K)$.

\subsection{Grothendieck-Teichm\"uller group $\wh{GT}(\K)$} The group of those automorphisms of the operad
$\wh{\cP a \cB}_\K$ which are identical on objects  is called the {\em Grothendieck-Teichm\"uller group}\, and is denoted by $\wh{GT}(\K)$ or simply by $\wh{GT}$ or even $GT$.\footnote{Since we understand $\wh{\cP a \cB}_\K$ as an operad in $\wh{\mathsf{Cat}}(\mathsf{coAss}_\K)$, it is tacitly assumed  that its every  automorphism respects the diagonal $\Delta$.}

 A similar group of automorphisms of  $\cP a \cB_\K^{(m)}$
is denoted by $GT^{(m)}(\K)$. The group  $\wh{GT}(\K)$ coincides with the inverse limit $\displaystyle\lim_{\longleftarrow} GT^{(m)}(\K)$.

\sip

Let us next characterize  elements of $\wh{GT}(\K)$ as solutions of certain explicit algebraic equations. This characterization proves that the above operadic definition of  $\wh{GT}(\K)$
coincides with the one given by V.\ Drinfeld in his original paper \cite{Dr}.

\subsubsection{{\bf Theorem} \cite{Dr}}\label{5: Drinfeld Theorem} {\em Let $\wh{F}_2(\K)$ be the pro-unipotent completion of the free group in two variables $x$ and $y$.
 There is a one-to-one correspondence between elements of
$\wh{GT}(\K)$ and pairs,
$$
\left(f\in  \wh{F}_2(\K), \la\in \K^*\right),
$$
 which satisfy the following equations,
\Beq\label{5: Main Th 1st cond}
f(x,y)=f(y,x)^{-1},
\Eeq
\Beq\label{5: Main Th 2nd cond}
f(x_3,x_1)x_3^{\mu} f(x_2,x_3)x_2^{\mu} f(x_1,x_2)x_1^{\mu}=1
\ \ \ \ \mbox{for}\ \ x_1x_2x_3=1 \ \ \mbox{and}\ \ \mu:=\frac{\la-1}{2} ,
\Eeq
and
\Beq\label{5: Main Th 3rd cond}
f(x_{12}, x_{23}x_{24})f(x_{13}x_{23}, x_{34}) = f(x_{23}, x_{34})f(x_{12}x_{13}, x_{23}x_{34})f(x_{12}, x_{23})\ \  \ \mbox{in}\ \ \ \wh{\P\B}_4(\K).
\Eeq
The multiplication law in $\wh{GT}(\K)$ is given by,
$$
\left(\la_1, f_1\right)\left(\la_2, f_2\right) = \left(\la_1\la_2, \ \ f_2(f_1(x,y)x^{\la_1}f_1(x,y)^{-1}, y^{\la_1} )
f_1(x,y)\right)
$$
}
\begin{proof}
By Theorem~{\ref{5: Theorem on morhisms from PaB}}, every automorphism $F$  of $\wh{\cP a \cB}$ is uniquely determined by its values on the generators $\beta_{1,2}$ and $\id_{1,2,3}$ so that we have only to check that these values satisfy pentagon equation (\ref{5: pentagon in PaB--f}), two hexagon equations (\ref{5: first hexagon in PaB--f}) and (\ref{5: second hehagon in PaB--f}), and that $F$ sends group-like elements in  $\wh{\cP a \cB}$  to group-like elements with respect to the diagonal $\Delta$. The latter condition is equivalent to saying that $F(\beta_{1,2})$ and
$F(\id_{1,2,3})$ are group-like in $\Mor_{\wh{\cP a \cB}}(\mbox{\scriptsize 12,21})$ and  $\Mor_{\wh{\cP a \cB}}( \mbox{\scriptsize (12)3,1(23)})$. Hence
$$
F(\beta_{1,2})\circ \beta_{1,2}^{-1}\ \ \mbox{is grouplike in the coalgebra}\ \  \Mor_{\wh{\cP a \cB}}(\mbox{\scriptsize 12,12})
$$
and
$$
F(\id_{1,2,3})\circ \id_{1,2,3}^{-1} \ \mbox{is grouplike in the coalgebra}\ \
\Mor_{\wh{\cP a \cB}}(\mbox{\scriptsize (12)3,1(23)})
$$
which in turn imply (see \S\S
{\ref{5: Corollary on PB_3}}
and {\ref{5: Remark on B_2}}) that
$$
F(\beta_{1,2})\circ \beta_{1,2}^{-1}= [\beta_{1,2}^{2}]^\mu\in \wh{\P\B_2}(\K)\
$$
for some $\mu\in \K$ with $\mu\neq -\frac{1}{2}$ (as otherwise $F$ would not be invertible),
and
$$
F(\id_{1,2,3})\circ \id_{1,2,3}^{-1}=f(\sfb_2^2,\sfb_1^2)[(b_1b_2)^{3}]^\eta\in \wh{\P\B_3}(\K)
$$
for some $f(x,y)\in \wh{F_2}(\K)$ and some $\eta\in \K$. Here  $\sfb_1^2$ and $\sfb_2^2$ stand for the
renormalized generators (cf.\ \S{\ref{5: corollary on sfb_1 sfb_2}}),
$$
\sfb_1^2:= c^{-1/3} b_1^2, \ \ \ \ \  \sfb_2^2:= c^{-1/3} b_2^2,
$$
and  $b_1^2$,  $b_2^2$ and $c:=(b_1b_2)^3$ are the grouplike elements of $\Mor_{\wh{\cP a \cB}}(\mbox{\scriptsize (12)3,(12)3})$
obtained from the corresponding braids in $\P\B_3$
by attaching strands' inputs and outputs to one and the same object $(12)3$, for example,
$$
b_1^2=\hspace{-2mm}\Ba{c}
\Ba{c}\resizebox{11mm}{!}{\xy
(-5,-17)*{_{(1}},
(-5,7)*{_{(1}},
(5,-17)*{_{2)}},
(5,7)*{_{2)}},
(9,7)*{_{3}},
(9,-17)*{_{3}},
\vtwist~{(-5,5)}{(5,5)}{(-5,-5)}{(5,-5)};
\vtwist~{(-5,-5)}{(5,-5)}{(-5,-15)}{(5,-15)};
(9,-15)*{}; (9,5)*{} **\dir{-};
\endxy}\Ea
\Ea= d_3(\beta_{1,2}^2)  , \ \ \ \
b_2^2:=\hspace{-1mm}
\Ba{c}\resizebox{11mm}{!}{
\xy
(-5,-17)*{_{2)}},
(-5,7)*{_{2)}},
(5,-17)*{_{3}},
(5,7)*{_{3}},
(-8,7)*{_{(1}},
(-8,-17)*{_{(1}},
\vtwist~{(-5,5)}{(5,5)}{(-5,-5)}{(5,-5)};
\vtwist~{(-5,-5)}{(5,-5)}{(-5,-15)}{(5,-15)};
(-8,-15)*{}; (-8,5)*{} **\dir{-};
\endxy}\Ea
=\id_{1,2,3}\circ
\Ba{c}\resizebox{11mm}{!}{
\xy
(-5,-17)*{_{(2}},
(-5,7)*{_{(2}},
(5,-17)*{_{3)}},
(5,7)*{_{3)}},
(-8,7)*{_{1}},
(-8,-17)*{_{1}},
\vtwist~{(-5,5)}{(5,5)}{(-5,-5)}{(5,-5)};
\vtwist~{(-5,-5)}{(5,-5)}{(-5,-15)}{(5,-15)};
(-8,-15)*{}; (-8,5)*{} **\dir{-};
\endxy}
\Ea\circ \id_{1,2,3}^{-1}
=
\id_{1,2,3}\circ d_0(\beta_{1,2}^2) \circ \id_{1,2,3}^{-1}
$$

Using (\ref{2: conjugation of g to the power la}) we obtain
\Beqrn
d_0F(\beta_{1,2}) &=& d_0([\beta_{1,2}^{2}]^\mu\circ\beta_{1,2}) \\
&=&  d_0([\beta_{1,2}^{2}]^\mu) \circ d_0 (\beta_{1,2})  \\
&=&\id_{1,2,3}^{-1}\circ [b_2^{2}]^\mu\circ \id_{1,2,3}\circ d_0\beta_{1,2}
\Eeqrn
so that
$$
F(\id_{1,2,3})\  \circ\  d_0F(\beta_{1,2})= \left[f(\sfb_2^2,\sfb_1^2) (\sfb_2^2)^\mu c^{\eta+ \frac{1}{3}\mu} \right]
\circ\id_{1,2,3}\circ d_0\beta_{1,2}.
$$
Denoting
$$
g=\id_{1,2,3}\circ d_0\beta_{1,2} \circ \mbox{(\scriptsize 23)}\id_{1,2,3}^{-1}=
\Ba{c}\resizebox{10mm}{!}{\xy
(-5.3,-8)*{_{(1}},
(-5,8)*{_{1}},
(0,-8)*{_{2)}},
(0,8)*{_{(2}},
(5,-8)*{_{3}},
(6,8)*{_{3)}},
(-5,-6)*{}; (-5,6)*{} **\dir{-};
(0,-6)*{}; (0,6)*{} **\dir{-};
(5,-6)*{}; (5,6)*{} **\dir{-};
\endxy}\Ea
\ \
\circ
\ \
\Ba{c}\resizebox{10mm}{!}{
\xy
(-4,-8)*{_{(2}},
(-4,8)*{_{(3}},
(4,-8)*{_{3)}},
(4,8)*{_{2)}},
(-8,-8)*{_{1}},
(-8,8)*{_{1}},
\vtwist~{(-4,6)}{(4,6)}{(-4,-6)}{(4,-6)};
(-8,-6)*{}; (-8,6)*{} **\dir{-};
\endxy}\Ea
\ \
\circ
\ \
\Ba{c}\resizebox{10mm}{!}{
\xy
(-5,-8)*{_{1}},
(-5,8)*{_{(1}},
(0,-8)*{_{(3}},
(0,8)*{_{3)}},
(5,-8)*{_{2)}},
(6,8)*{_{2}},
(-5,-6)*{}; (-5,6)*{} **\dir{-};
(0,-6)*{}; (0,6)*{} **\dir{-};
(5,-6)*{}; (5,6)*{} **\dir{-};
\endxy}\Ea\  \ =\ \ \Ba{c}\resizebox{10mm}{!}{\xy
(-3.5,-8)*{_{2)}},
(-3.5,8)*{_{3)}},
(4,-8)*{_{3}},
(4,8)*{_{2}},
(-8,-8)*{_{(1}},
(-8,8)*{_{(1}},
\vtwist~{(-4,6)}{(4,6)}{(-4,-6)}{(4,-6)};
(-8,-6)*{}; (-8,6)*{} **\dir{-};
\endxy}\Ea   
$$
we obtain, using the results of \S{\ref{5: Exercise on b_2b_1 squared b_2 inverse}},
\Beqrn
g\circ \left[ \mbox{(\scriptsize 23)}\hspace{-1mm} \cdot\hspace{-1mm}  f(\sfb_2^2,\sfb_1^2)\right]^{-1} \circ g^{-1}  &=&
 f\left(g \circ  \mbox{(\scriptsize 23)}
 \hspace{-1mm} \cdot\hspace{-1mm} \sfb_2^2 \circ g^{-1},    g \circ  \mbox{(\scriptsize 23)}
 \hspace{-1mm} \cdot\hspace{-1mm} \sfb_1^2 \circ g^{-1}  \right)^{-1}\\
 &=& f(\sfb_2^2, \sfb_1^{-2}\sfb_2^{-2} )^{-1}
\Eeqrn
so that the product of the first three terms on the r.h.s.\ of  (\ref{5: first hexagon in PaB--f})
now reads,
$$
F(\id_{1,2,3})\  \circ\  d_0F(\beta_{1,2}) \circ \ \left(\mbox{(\scriptsize 23)}\ccdot F(\id_{1,2,3})^{-1}\right)=
\left[f(\sfb_2^2,\sfb_1^2) (\sfb_2^2)^\mu  f(\sfb_2^2,\sfb_1^{-2}\sfb_2^{-2})^{-1} c^{\frac{1}{3}\mu}\right] \circ g.
$$
Next we have,
\Beqrn
g\circ \mbox{(\scriptsize 23)} \ccdot d_3F(\beta_{1,2})\circ g^{-1} &=&
g\circ \mbox{(\scriptsize 23)} \ccdot [b_1^{2}]^\mu \circ g^{-1}\circ g\circ \mbox{(\scriptsize 23)} \ccdot  d_3\be_{1,2} \circ g^{-1}\\
&=& \left[(\sfb_1^{-2}\sfb_2^{-2})^\mu c^{\frac{1}{3}\mu}\right]\circ  g\circ \mbox{(\scriptsize 23)} \ccdot  d_3\be_{1,2} \circ g^{-1},
\Eeqrn
and hence
$$
F(\id_{1,2,3})\  \circ\  d_0F(\beta_{1,2}) \circ \ \left(\mbox{(\scriptsize 23)}\ccdot F(\id_{1,2,3})^{-1}\right)\circ
\mbox{(\scriptsize 23)} \ccdot d_3F(\beta_{1,2})=\hspace{60mm}
$$
$$
\hspace{60mm}=\left[f(\sfb_2^2,\sfb_1^2) (\sfb_2^2)^\mu  f(\sfb_2^2,\sfb_1^{-2}\sfb_2^{-2} )^{-1}(\sfb_1^{-2}\sfb_2^{-2})^\mu
c^{\frac{2}{3}\mu}\right] \circ h.
$$
where
$$
h:= g\circ \mbox{(\scriptsize 23)} \ccdot  d_3\be_{1,2}=
\Ba{c}\resizebox{12mm}{!}{
\xy
(-3.5,-8)*{_{2)}},
(-3.5,20)*{_{1)}},
(4,-8)*{_{3}},
(4,20)*{_{2}},
(-12,-8)*{_{(1}},
(-12,20)*{_{(3}},
\vtwist~{(-4,6)}{(4,6)}{(-4,-6)}{(4,-6)};
\vtwist~{(-12,18)}{(-4,18)}{(-12,6)}{(-4,6)};
(-12,-6)*{}; (-12,6)*{} **\dir{-};
(4,6)*{}; (4,18)*{} **\dir{-};
\endxy}\Ea  =d_1\be_{1,2}\circ \mbox{(\scriptsize 123)}\ccdot \id_{1,2,3}^{-1}.
$$
As $(b_2b_1)(\sfb_1^2)(b_2b_1)^{-1}= \sfb_1^{-2}\sfb_2^{-2}$ and
$(b_2b_1)(\sfb_2^2)(b_2b_1)^{-1}=\sfb_1^2$ in $\wh{\B}_3(\K)$, we have
$$
h\circ \mbox{(\scriptsize 123)}\ccdot F(\id_{1,2,3}) \circ h^{-1}=
[f( \sfb_1^2, \sfb_1^{-2}\sfb_2^{-2})]\circ h\circ  \mbox{(\scriptsize 123)}\ccdot \id_{1,2,3}\circ h^{-1}
$$
and hence
$$
\mbox{r.h.s.\ of}\  (\ref{5: first hexagon in PaB--f})=
\left[f(\sfb_2^2,\sfb_1^2) (\sfb_2^2)^\mu  f(\sfb_2^2,\sfb_1^{-2}\sfb_2^{-2} )^{-1}(\sfb_1^{-2}\sfb_2^{-2})^\mu f(\sfb_1^2,\sfb_1^{-2}\sfb_2^{-2} ) c^{\eta+\frac{2}{3}\mu}\right]\circ h\circ  \mbox{(\scriptsize 123)}\ccdot \id_{1,2,3}.
$$
On the other hand, one has $d_1 b_1^2=b_2b_1^2b_2=(b_1b_2)^3b_1^{-2}=c^{2/3}\sfb_1^{-2}$ in $\wh{\B}_3(\K)$, and hence
$$
\mbox{l.h.s.\ of}\  (\ref{5: first hexagon in PaB--f})=
d_1 F(\beta_{1,2})=[c^{\frac{2}{3}\mu}(\sfb_1^{-2})^\mu]\circ h \circ \mbox{(\scriptsize 123)}\ccdot \id_{1,2,3}.
$$
Therefore, the first hexagon equation (\ref{5: first hexagon in PaB--f})
is equivalent to the following equation in $\wh{\P\B}_3(\K)$,
$$
f(\sfb_2^2,\sfb_1^2) (\sfb_2^2)^\mu  f(\sfb_2^2,\sfb_1^{-2}\sfb_2^{-2} )^{-1}(\sfb_1^{-2}\sfb_2^{-2})^\mu
f(\sfb_1^2,\sfb_1^{-2}\sfb_2^{-2} )
(\sfb_1^{2})^\mu c^{\eta}=1
$$
Hence $\eta=0$ and the element $f(x,y)\in \wh{F_2}(\K)$ satisfies the
equation
\Beq\label{5: first heganon for f(x,y)}
f(x_3, x_1)x_3^\mu f(x_3,x_2)^{-1} x_2^\mu f(x_1,x_2)x_1^\mu=1
\Eeq
for $x_1x_2x_3=1$.
Studying similarly the second  hexagon equation (\ref{5: second hehagon in PaB--f}), one obtains
(using identities $b_2^{-1}\sfb_1^2 b_2= \sfb_2^{-2}\sfb_1^{-2}$,
 $(b_1b_2)^{-1}(\sfb_1^2)(b_1b_2)= \sfb_2^{-2}\sfb_1^{-2}$ and
$(b_1b_2)^{-1}(\sfb_2^2)(b_1b_2)=\sfb_1^2$) the following equation,
$$
(\sfb_1^{2})^{-\mu}f(\sfb_2^2,\sfb_1^2) (\sfb_2^{2})^{-\mu}
f(\sfb_2^2,\sfb_2^{-2}\sfb_1^{-2} )^{-1}(\sfb_2^{-2}\sfb_1^{-2})^{-\mu} f(\sfb_1^2,\sfb_2^{-2}\sfb_1^{-2} )
 =1
$$
which says that $f(x,y)\in \wh{F_2}(\K)$ satisfies
\Beq\label{5: second hexagon for f(x,y)}
f(x_3, x_2)^{-1}x_2^\mu f(x_1,x_2) x_1^\mu f(x_1,x_3)^{-1}x_3^\mu=1.
\Eeq
for $x_1x_2x_3=1$.

\sip

Equations (\ref{5: first heganon for f(x,y)}) and (\ref{5: second hexagon for f(x,y)}) imply the
first two claims, (\ref{5: Main Th 1st cond})
and (\ref{5: Main Th 2nd cond}), of the Theorem.

\sip

Note that relation between the generators $x_{ij}$ of $\P\B_n$ are invariant under
their renormalizations, $x_{ij}\rar c^{\la_{ij}} x_{ij}$,  by powers of the central element
$c=(b_1b_2\cdots b_n)^n$. Therefore, when computing  $d_i(f(\sfb_2^2,\sfb_1^2))$, $i=0,1,2,3$,
 in (\ref{5: pentagon in PaB--f}) we can replace $\sfb_1^2$ and $\sfb_2^2$ with $b_2^2$ and $b_2^2$.
 As one has in $\P\B_4$
$$
d_1(b_2^2)= x_{34}, \ \ d_1(b_1^2)= x_{13}x_{23},\ \ d_3(b_2^2)= x_{23}x_{24}, \ \ d_3(b_1^2)= x_{12},
$$
the l.h.s.\ of the ``pentagon" equation (\ref{5: pentagon in PaB--f}) reads,
$$
d_1F(\id_{1,2,3}) \circ d_3F(\id_{1,2,3})=f(x_{34},x_{13}x_{23})f(x_{23}x_{24}, x_{12}).
$$
Continuing in the same way we conclude that   (\ref{5: pentagon in PaB--f}) is equivalent to
the following equality in $\wh{\P\B}_4(\K)$,
$$
f(x_{34},x_{13}x_{23})f(x_{23}x_{24}, x_{12})=f(x_{23}, x_{12})f(x_{24}x_{34}, x_{12}x_{13})
f(x_{34}, x_{23}).
$$
Taking inverses of both sides and using (\ref{5: Main Th 1st cond}) we obtain from the latter
equation the third claim,
equation (\ref{5: Main Th 3rd cond}), of the Theorem.

\sip

The final claim about the formula for the group multiplication in $\wh{GT}(\K)$ is left as an exercise.
\end{proof}

\subsubsection{\bf Definition} There is a surjection of groups,
$$
\Ba{ccc}
\wh{GT}(\K) &\lon & \K^*\\
(f,\la) &\lon & \la.
\Ea
$$
Its kernel is a subgroup of $\wh{GT}(\K)$ often denoted by $\wh{GT}_1(\K)$.

\subsection{H.\ Furusho's Theorem} A remarkable result of H.\ Furusho \cite{Fu} gives us a much shorter
algebraic characterization of elements of the Grothendieck-Teichm\"uller group than the one given
in the above Theorem {\ref{5: Drinfeld Theorem}}: the pentagon equation (\ref{5: Main Th 3rd cond}) often implies in the hexagon equations (\ref{5: Main Th 1st cond}) and (\ref{5: Main Th 2nd cond}).
 More precisely, one has

\subsubsection{{\bf Theorem} \cite{Fu}} {\em Let $\K$ be a field of characteristic $0$ and
$\overline{\K}$ be its algebraic closure. Suppose that an  element $f = f(x, y)\in \wh{F}_2(\K)$  satisfies Drinfeld's pentagon equation (\ref{5: Main Th 3rd cond}). Then there exists an element
 $\la\in \overline{\K}$ (which is unique up to a sign) such that
the pair $(\la, f)$ satisfies hexagon equations (\ref{5: Main Th 1st cond}) and
(\ref{5: Main Th 2nd cond}). Moreover, this $\la$ is equal to $\pm\sqrt{24c_2(f) + 1}$,
where $c_2(f)$ stands for the coefficient of $XY$ in the formal power series $f(e^X, e^Y)$.}

\mip

Therefore, the group $\wh{GT}(\overline{\K})$ can be identified with the set of pairs
$(f = f(x, y)\in \wh{F}_2(\overline{\K}), \pm\sqrt{24c_2(f) + 1})\in \overline{\K}\setminus 0) $, where $f$ is a solution of Drinfeld's pentagon equation with
$c_2(f)\neq -\frac{1}{24}$.

\bip

{\large
\section{\bf Infinitesimal braids, $GRT_1$  and associators}

\sip
\subsection{Graded Lie algebra of $\P\B_n$} Let $\fg\fr_\K$ be the functor associating to a discrete group
the graded Lie algebra obtained from the descending central series of that group
(see \S{\ref{3: tensor functors}}(viii) for details). It was proven in \cite{Koh}
that $\fg\fr_\K(\P\B_n)$ can be identified with the quotient of the free Lie algebra
$\ffree_{n(n-1)/2}$  generated by the symbols $\{t_{ij}=t_{ji}\}_{1\leq i< j\leq n}$ modulo the
ideal generated by relations,
\Beqrn
[t_{ij},t_{kl}] &= & 0\ \ \ \mbox{if}\ \ \#\{i,j,k,l\}=4  \\
\mbox{$ [t_{ij}, t_{ik} + t_{kj}]$} &=& 0 \ \ \ \mbox{if}\ \ \#\{i,j,k\}=3
\Eeqrn

The Lie algebra $\ft_n=\fg\fr_\K(\P\B_n)$ is called the {\em Lie algebra
of infinitesimal braids}. This is a filtered Lie algebra whose completion is denoted by $\wh{\ft}_n$ (which is equal, of course, to $\wh{\fg\fr}_\K(\P\B_n)$).

\subsubsection{\bf Lemma} {\em The element $c_n:=\sum_{i< j}t_{ij}$ is central in $\ft_n$}.
\begin{proof}
$\displaystyle [t_{ij}, c_n]=\sum_{k\neq i,j}[t_{ij}, t_{ik} + t_{jk}]=0$.
\end{proof}

\subsubsection{\bf Lemma} (i) {\em  The Lie subalgebra of\, $\ft_n$ generated by
$t_{ij}$ with $i,j\in [n-1]$ can be identified with $\ft_{n-1}$.}

(ii) {\em  The Lie subalgebra of\, $\ft_n$ generated by
$t_{1,n},t_{1,n}\ldots,t_{n-1,n}$ is a Lie ideal and can be identified with the free Lie algebra
$\ffree_n$.}

(iii) {\em There exists an isomorphism of Lie algebras},
$$
\ft_n\simeq \ft_{n-1}\ \oplus \ \ffree_n.
$$

The statement (i) is obvious. The statement (ii) follows from the fact that, for any $i,j,k\in [n-1]$, $i\neq j$,  one has
$$
[t_{ij}, t_{kn}]=\left\{\Ba{ll}
0 & \mbox{if}\ k\neq i, k\neq j \\

[t_{in},t_{jn}] & \mbox{if}\ k=i.
\Ea \right.
$$
We skip the proof of (iii).

\subsubsection{\bf Corollary}\label{6: corollary on t_3} {\em There is an isomorphism of Lie algebras,
$$
\ft_3=\K c_3 \oplus \ffree_2,
$$
 where $\K c_3$ is the Abelian Lie algebra generated by $c_3=t_{12}+t_{13}+t_{23}$, and  $\ffree_2$ is the free Lie algebra generated by $t_{13}$ and $t_{23}$ (or, equivalently, by $t_{12}$ and $t_{23}$)}.

 \subsubsection{\bf Exercise}\label{5: Exercise on graded of pb} As $\P\B_n$ is a group of finite type, its prounipotent completion $\wh{\P\B}_n(\K)$ over a field $\K$
 can be identified with the group $\exp(\wh{\fp\fb}_n(\K))$, where the complete Lie algebra
 $\wh{\fp\fb}_n(\K)$ is the quotient of the free Lie algebra $\wh{\ffree}_{n(n-1)/2}$  on letters $\{\ga_{ij}\}_{1\leq i< j\leq n}$ modulo the ideal generated by following formal power series of Lie words,
\Beqrn
\log\left(e^{\ga_{ij}}e^{\ga_{kl}}e^{-\ga_{ij}}e^{-\ga_{kl}} \right) &=& 0\ \ \mbox{for}\ \ \ \ \ \ i<j<k<l\ \ \mbox{and}\ \ i<k<l<j,\\
\log\left(e^{\ga_{ij}}e^{\ga_{ik}}e^{\ga_{jk}}e^{-\ga_{ik}}e^{-\ga_{jk}}e^{-\ga_{ij}}\right) &=& 0 \ \ \ \mbox{for}\ \ i<j<k\\
\log\left(e^{\ga_{ik}}e^{\ga_{jk}}e^{\ga_{ij}}e^{-\ga_{jk}}e^{-\ga_{ij}}\ga^{-\ga_{ik}}\right)&=& 0 \ \ \ \mbox{for}\ \ i<j<k\\
\log\left(e^{\ga_{ik}}e^{\ga_{jk}}e^{\ga_{jl}}e^{-\ga_{jk}} e^{-\ga_{jk}}e^{-\ga_{jl}}
e^{\ga_{jk}}e^{-\ga_{ik}}\right)&=& 0 \ \ \mbox{for}\ \ i<j<k<l.
\Eeqrn
The Lie algebra $\wh{\fp\fb}_n(\K)$ has a natural filtration,
$$
\wh{\fp\fb}_n(\K)=F_1(\wh{\fp\fb}_n(\K))\supset F_2(\wh{\fp\fb}_n(\K)) \supset \ldots
\supset F_p(\wh{\fp\fb}_n(\K)) \supset F_{p+1}(\wh{\fp\fb}_n(\K)) \supset \ldots,
$$
 where
$F_p(\wh{\fp\fb}_n(\K))$ is the Lie subalgebra generated by Lie words with at at least $p$ letters. Let
$$
gr(\wh{\fp\fb}_n(\K)):=\prod_{i=1}^\infty \frac{F_{p+1}(\wh{\fp\fb}_n(\K))}{F_p(\wh{\fp\fb}_n(\K))}
$$
be the associated graded Lie algebra, and let $[\ga_{ij}]$ be the image of the generators
under the projection $gr(\wh{\fp\fb}_n(\K))\rar \wh{\fp\fb}_n(\K)$.

\sip

Show that the association $t_{ij}\lon [\ga_{ij}]$ induces an isomorphism of Lie algebras,
$$
\al:  \wh{\ft}_n \lon  gr(\wh{\fp\fb}_n(\K))
$$

\subsection{Operadic structure on Lie algebras of infinitesimal braids \cite{Ta2}} Consider
a functor from the groupoid of finite sets to the category of Lie algebras,
$$
\Ba{rccc}
\ft: & \sS & \lon & \mathsf{Lie}_\K\\
     & I   & \lon & \ft(I)
     \Ea
$$
where $t(I)$ is the Lie algebra of infinitesimal braids on the set $I$, that is,
 the quotient of the free Lie algebra generated by symbols $t_{ij}=t_{ji}$, $i,j\in I$,
$i\neq j$, modulo the Lie ideal generated by the relations $[t_{ij}, t_{kl}]=0$ whenever
$i,j,k,l$ are four different elements of $I$, and $[t_{ij}, t_{ik}+ t_{jk}]=0$ whenever
$i,j,k$ are three different elements. We set $\ft(I)=0$ if cardinality of $I$ is $1$.

\subsubsection{\bf Exercises} (i) Show that for any {\em injective}\, morphism
of sets $f:I\rar J$ there is an associated {\em push-forward}\,   homomorphism
of Lie algebras
 $$
\Ba{rccc}
f_*: & \ft(I) & \lon & \ft(J)\\
     & t_{ij}   & \lon & f(t_{ij}):=t_{f(i)f(j)}.
     \Ea
$$

(ii) Show that for any morphism
of sets $f:I\rar J$ there is an associated {\em pull-back}\, Lie algebra  homomorphism
 $$
\Ba{rccc}
f^*: & \ft(J) & \lon & \ft(I)\vspace{3mm}\\
     & t_{\al\be}   & \lon & \displaystyle f^*(t_{\al\be}):=\hspace{-3mm}\sum_{i,j\in I\atop
     f(i)=\al, f(j)=\be} t_{ij}.
     \Ea
$$

(iii) Let
$$
I\xrightarrow{f} J \xrightarrow{g} K
$$
be a sequence of sets and their morphisms such that $\Img f\circ g$ consists of at
most one element. Show that the Lie subalgebras $f_*(I), g^*(K)\subset \ft(J)$ satisfy
$$
[f_*(I), g^*(K)]=0.
$$

\subsubsection{\bf Lemma} {\em The $\sS$-module $\ft$ is an operad
in the category $\mathsf{Lie}_\K$ with respect to the  composition which
is given on the generators as follows,}
 $$
\Ba{rccll}
\circ^{I,J}_k: & \ft(I)\oplus \ft(J) & \lon & \ft(I\setminus k \sqcup J) & \forall I,J,\ k\in I
\vspace{3mm}\\
& (0, t_{\al\be}) & \lon & t_{\al\be} & \forall\ \al,\be\in J\vspace{3mm} \\
& (t_{ij}, 0)
& \lon & \displaystyle \left\{\Ba{ll}
     t_{ij}  &\mbox{if}\ \  i,j\in I\setminus k\vspace{3mm} \\
     \sum_{p\in J} t_{jp}  &\mbox{if}\ \ i= k.\\
     \Ea\right. &
     \Ea
$$
\begin{proof} 
Consider the maps of sets,
$$
\hspace{-19mm}
\Ba{rccc}
f: & J &\lon & I\setminus k\sqcup J\\
   & p & \lon & f(p):=p
\Ea
$$
and
$$
\hspace{30mm}
\Ba{rccl}
g: & I\setminus k \sqcup J & \lon & I\\
   &        x   & \lon & g(x):=\displaystyle \left\{\Ba{ll}
     i  & \mbox{if}\ \ x=i \in I\setminus k\vspace{3mm} \\
     k  & \mbox{if}\ \ x=p\in J.\\
     \Ea\right.
   \Ea
$$
with $f$ being an injection, and $\Img f\circ g$ being an subset of $I$ of cardinality $1$.
Then the operadic composition $\circ^{I,J}_k$ can be identified with the map $f_*\oplus g^*$.
Hence the Exercises
just above imply that  $\circ_k^{I,J}$ is a Lie algebra homomorphism.
We leave it to the reader to check that the homomorphisms $\circ_k^{I,J}$ satisfy the
axioms of a (non-unital) operad.
\end{proof}

\sip

The $\sS$-module $\wh{\ft}$ of {\em completed}\, Lie algebras of infinitesimal braids is,
of course, also an operad on the category $\mathsf{Lie}_\K$.

\subsubsection{\bf Exercise} Check that the operadic composition
$$
\Ba{rccc}
\circ_2: & \ft(3) \oplus  \ft(3)\lon & \ft(5)\\
        & \displaystyle \left([t_{12},t_{23}], [t_{13},t_{23}]\right)
        &\lon & [t_{12},t_{23}] \circ_2  [t_{13},t_{23}]
\Ea
$$
gives
$$
[t_{12},t_{23}] \circ_2  [t_{13},t_{23}]=[t_{12}+ t_{13}+t_{14},t_{25}+ t_{35}+ t_{45}] + [t_{24},t_{34}]
$$

 \subsection{Cosimplicial complex of the operad of infinitesimal braids}\label{6: simp t}
 Let $e:=0$ be the zero element
 in $\ft(2)$. It obviously satisfies the equation (\ref{4: ass cond for e}) and hence makes by Lemma \S{\ref{4: Lemma on cosimpl operad structures}} the collection $\ft=\{\ft(n)\}_{\geq 1}$
into the cosimplicial object in the category $\mathsf{Lie}_\K$. The cohomology
of the associated cosimplicial dg Lie algebra,
$$
\mathsf{Simp}^\bu(\ft)=\left(\bigoplus_{n\geq 1}\ft_n[2-n], d=\sum_{i}(-1)^id_i\right),
$$
was studied by Thomas Willwacher in his breakthrough paper \cite{Wi1}. Its relation to the Grothendieck-Teichm\"uller group is discussed below.

\subsection{Operad of chord diagrams}  Applying the universal enveloping monoidal
functor $U: \mathsf{Lie}_\K \rar \mathsf{Hopf}_\K$ we get an operad $$
\cC\caD_\K:=U(\ft)=\{\cC\caD_\K(n):=U(\ft(n))\}
$$
called in \cite{BN} the operad of {\em chord diagrams} because its elements, say, elements of $\cC\caD_\K(n)$, can be pictorially presented as a collection of chords on $n$ vertical strands (labelled by $1,2,\ldots, n$ from left to right) e.g.
$$
\cC\caD_\K(3) \ni t_{23}t_{12} \longleftrightarrow
\xy
(-5,-4)*{}="1",
(-5,6)*{}="1'",
(0,-4)*{}="2",
(0,6)*{}="2'",
(5,-4)*{}="3",
(5,6)*{}="3'",
(0,-2)*{}="a",
(5,-2)*{}="a'",
(0,2)*{}="b",
(-5,2)*{}="b'",
\ar @{->} "1";"1'" <0pt>
\ar @{->} "2";"2'" <0pt>
\ar @{->} "3";"3'" <0pt>
\ar @{-} "a";"a'" <0pt>
\ar @{-} "b";"b'" <0pt>
 \endxy
 $$
As an associative algebra $\cC\caD_\K(n)$, $n\geq 2$, is the quotient of the free non-commutative algebra generated by symbols $t_{ij}=t_{ji}$, $i\neq j$ $i,j\in [n]$, modulo the ideal generated by the relations
$[t_{ij},t_{kl}]=  0$ if $\#\{i,j,k,l\}=4$ and  $[t_{ij}, t_{ik} + t_{kj}] =0$ if $ \#\{i,j,k\}=3$.
The coproduct in $\cC\caD_\K(n)$ is given on the generators as follows,
$$
\Delta(t_{ij})=1\ot t_{ij} + t_{ij}\ot 1.
$$
The operadic structure in $\cC\caD_\K$ is completely determined by the one in $\ft$. For example, one has
$$
\Ba{rccc}
\circ_2: & \cC\caD_\K(3) \ot_\K  \cC\caD_\K(3)\lon & \cC\caD_\K(5)\\
        & \displaystyle\left(t_{12}t_{23}, t_{13}t_{23}\right)
        &\lon & (t_{12}+ t_{13}+t_{14})(t_{25}+ t_{35}+ t_{45})t_{24}t_{34}.
\Ea
$$

The $I$-adic completion of the associative algebra $\cC\caD_\K(n)$ with respect to the maximal ideal $I$
is denoted by $\wh{\cC\caD}_\K(n)$. The collection $\wh{\cC\caD}_\K=\{\wh{\cC\caD}_\K(n)\}_{n\geq 1}$
 is an operad in the category of complete filtered Hopf algebras.
For our purposes it is useful to view $\wh{\cC\caD}_\K$ instead as an operad in the category $\mathsf{Cat(Coass_\K)}$
as follows: each competed filtered  Hopf algebra $\wh{\cC\caD}_\K(n)=\wh{U}(\ft_n)$ is understood from now on as a category with one object $\bu$ and with
$Mor(\bu,\bu):=\wh{\cC\caD}_\K(n)$ equipped with its standard coalgebra structure.

\subsection{Cosimplicial structure on $\wh{\cC\caD}_\K$} The element $e:=1$ in $\cC\caD(2)$ satisfies equations (\ref{4: ass cond for e}) so that formulae (\ref{4: formulae for d_i})
make the collection $\wh{\cC\caD}_\K=\{\wh{U}(\ft_n)\}_{n\geq 1}$ into a cosimplicial object in the category $\mathsf{dgAss}_\K$.

\subsubsection{\bf Exercise} Let $t_{12}^2t_{23}\in \wh{U}(\ft_3)$. Check that
$$
d_0\left(t_{12}^2t_{23}\right)= t_{23}^2t_{34}, \ \ \ d_1\left(t_{12}^2t_{23}\right)= \left(t_{12}+t_{13}\right)^2t_{34},\ \ \ \
d_2\left(t_{12}^2t_{23}\right)=\left(t_{12}+t_{13}\right)^2\left(t_{24}+t_{34}\right).
$$


\subsection{D.\ Bar-Natan's operad of parenthesized chord diagrams} We have two operads in the symmetric monoidal category $\mathsf{Cat(Coass_\K)}$, $\cP a_\K$ defined in \S{\ref{5: subsection on wh{PaB}_K}}
and $\wh{\cC\caD}_\K$ defined just above. Hence
$$
\wh{\cP a\cC\caD}_\K:= \cP a_\K \ot_{\mathsf{Cat(Coass_\K)}} \wh{\cC\caD}_\K
$$
is an operad in the category $\mathsf{Cat(Coass_\K)}$ called in \cite{BN} the operad of
{\em parenthesized chord diagrams}. Elements of  $\cP a_\K$ can be pictorially represented as elements
of $\cP a \cB_\K$ but with braiding of strands forgotten so that we represent morphisms in the category
$\wh{\cP a\cC\caD}_\K(n)$ as linear combinations, with coefficients in the algebra  $\wh{\cC\caD}_\K(n)$, of non-braided strands connecting identical symbols in a pair of parenthesized permutations of $[n]$,
for example
$$
\mathbb{H}_{1,2}:=t_{12}\cdot
\Ba{c}\resizebox{5mm}{!}{
\xy
(-2,8)*{_{1}},
(-2,-8)*{_{1}},
(2,-8)*{_{2}},
(2,8)*{_{2}},
(-2,-6)*{}; (-2,6)*{} **\dir{-};
(2,-6)*{}; (2,6)*{} **\dir{-};
\endxy}\Ea\in \Mor_{\wh{\cP a\cC\caD}_\K(2)}(12,12), \ \ \
\mathbb{X}_{1,2}=1\cdot
\Ba{c}\resizebox{5mm}{!}{
\xy
(-2,8)*{_{2}},
(-2,-8)*{_{1}},
(2,-8)*{_{2}},
(2,8)*{_{1}},
(-2,-6)*{}; (2,6)*{} **\dir{-};
(2,-6)*{}; (-2,6)*{} **\dir{-};
\endxy}\Ea \in \Mor_{\wh{\cP a\cC\caD}_\K(2)}(12,21)
$$
$$
\id_{1,2,3}:=
\Ba{c}\resizebox{9mm}{!}{
\xy
(-5,-8)*{_{(1}},
(-5,8)*{_{1}},
(0,-8)*{_{2)}},
(0,8)*{_{(2}},
(5,-8)*{_{3}},
(6,8)*{_{3)}},
(-5,-6)*{}; (-5,6)*{} **\dir{-};
(0,-6)*{}; (0,6)*{} **\dir{-};
(5,-6)*{}; (5,6)*{} **\dir{-};
\endxy}\Ea\in \Mor_{\wh{\cP a\cC\caD}_\K(3)}\left((12)3), 1(23)\right)\ \ \ \
$$

\subsubsection{\bf Exercise} Show that the operad $\wh{\cP a \cC\caD}_\K$ is generated
by $\mathbb{H}_{1,2}$, $\mathbb{X}_{1,2}$ and $\id_{1,2,3}$.

\subsubsection{\bf Exercise}
Let $\id_{12}$ be the identity element in $\Mor_{\wh{\cP a \cC\caD}_\K(2)}(12, 12)$,
$$
 \id_{1,2}:=
 \Ba{c}\resizebox{5mm}{!}{
\xy
(-2,8)*{_{1}},
(-2,-8)*{_{1}},
(2,-8)*{_{2}},
(2,8)*{_{2}},
(-2,-6)*{}; (-2,6)*{} **\dir{-};
(2,-6)*{}; (2,6)*{} **\dir{-};
\endxy}\Ea
$$
Show that,  for any $n\in \N$, the family of functors,
$$
d_i: \wh{\cP a\cC\caD}_\K(n)   \lon  \wh{\cP a\cC\caD}_\K(n+1), \ \ \ \ i=0,1,\ldots, n, n+1,
$$
defined on objects and morphisms
of the category $\wh{\cP a\cC\caD}_\K(n)$ by formulae identical to the ones in
\S{\ref{5: cosimiplic struc on PaB}},
 makes the collection of categories $\{\wh{\cP a\cC\caD}_\K(n)\}_{n\geq 1}$
a (pre)cosimplicial object in the category $\mathsf{Cat(Coass_\K)}$.

\subsubsection{\bf Exercise}
Check that
$$
d_0 \bbH_{1,2}=\id_{1,2}\circ_2 \bbH_{1,2}= t_{23}
\Ba{c}\resizebox{9mm}{!}{
\xy
(-5,-8)*{_{1}},
(-5,8)*{_{1}},
(0,-8)*{_{(2}},
(0,8)*{_{(2}},
(6,-8)*{_{3)}},
(6,8)*{_{3)}},
(-5,-6)*{}; (-5,6)*{} **\dir{-};
(0,-6)*{}; (0,6)*{} **\dir{-};
(5,-6)*{}; (5,6)*{} **\dir{-};
\endxy}\Ea,\ \ \ \ \
d_1 \bbH_{1,2}= \bbH_{1,2}\circ_1\id_{1,2}= (t_{13}+t_{23})
\Ba{c}\resizebox{9mm}{!}{
\xy
(-5,-8)*{_{(1}},
(-5,8)*{_{(1}},
(0,-8)*{_{2)}},
(0,8)*{_{2)}},
(5,-8)*{_{3}},
(5,8)*{_{3}},
(-5,-6)*{}; (-5,6)*{} **\dir{-};
(0,-6)*{}; (0,6)*{} **\dir{-};
(5,-6)*{}; (5,6)*{} **\dir{-};
\endxy}\Ea,
$$
$$
d_2 \bbH_{1,2}= \bbH_{1,2}\circ_2\id_{1,2}= (t_{12}+t_{13})
\Ba{c}\resizebox{9mm}{!}{
\xy
(-5,-8)*{_{1}},
(-5,8)*{_{1}},
(0,-8)*{_{(2}},
(0,8)*{_{(2}},
(6,-8)*{_{3)}},
(6,8)*{_{3)}},
(-5,-6)*{}; (-5,6)*{} **\dir{-};
(0,-6)*{}; (0,6)*{} **\dir{-};
(5,-6)*{}; (5,6)*{} **\dir{-};
\endxy}\Ea,\ \ \ \ \
d_3 \bbH_{1,2}=\id_{1,2}\circ_1 \bbH_{1,2}= t_{12}
\Ba{c}\resizebox{9mm}{!}{
\xy
(-5,-8)*{_{(1}},
(-5,8)*{_{(1}},
(0,-8)*{_{2)}},
(0,8)*{_{2)}},
(5,-8)*{_{3}},
(5,8)*{_{3}},
(-5,-6)*{}; (-5,6)*{} **\dir{-};
(0,-6)*{}; (0,6)*{} **\dir{-};
(5,-6)*{}; (5,6)*{} **\dir{-};
\endxy}\Ea.
$$
\subsubsection{\bf Exercise}\label{6: Exercise on Phi}
 Let
$$
\Phi= \Phi(t_{12},t_{23}) \id_{1,2,3} \in  \Mor_{\wh{\cP a\cC\caD}_\K(3)}\left((12)3), 1(23)\right)
$$
 for some formal power series $\Phi(t_{12},t_{23})\in \K\langle\langle t_{12},t_{23}\rangle\rangle \subset \wh{U}(\ft_3)$. Check that
$$
d_0\Phi=\id_{1,2}\circ_2 \Phi= \Phi(t_{23},t_{34})d_0\id_{1,2,3},
\ \ \ d_1\Phi= \Phi \circ_1\id_{1,2}= \Phi(t_{12}+t_{13},t_{34})d_1\id_{1,2,3}
$$
$$
d_2\Phi= \Phi \circ_2\id_{1,2}= \Phi(t_{12}+t_{13},t_{24}+t_{34})d_2\id_{1,2,3},
$$
$$
d_3\Phi= \Phi \circ_3\id_{1,2}= \Phi(t_{12},t_{23}+t_{24})d_3\id_{1,2,3},\ \ \ \
d_4\Phi=\id_{1,2}\circ_1 \Phi =\Phi(t_{12},t_{23})d_4\id_{1,2,3}.
$$

\subsection{Grothendieck-Teichm\"uller group $GRT_1$ and its Lie algebra} The group of automorphisms of the operad $\wh{\cP a\cC\caD}_\K$ which preserve elements $\bbH_{1,2}$ and $\mathbb{X}_{1,2}$ is called the {\em graded Grothendieck-Tecihm\"uller group} and is denoted by $GRT_1$ \cite{Dr,BN,Ta2}.

\subsubsection{\bf Theorem} {\em There is a one-to-one correspondence between elements of $GRT_1$ and
grouplike elements, $\Phi(x,y)=e^{\phi(x,y)}$, of the completed universal enveloping (Hopf) algebra $\wh{U}({\ffree}_2)\simeq \K\langle\langle x,y\rangle\rangle$ of the free Lie algebra on generators $x$ and $y$ which satisfy the following equations
\Beq\label{6: TH on GRT 1st eqn - skewsym}
\Phi(x,y)=\Phi^{-1}(y,x),
\Eeq
\Beq\label{6: TH on GRT 2nd eqn - hexagon}
\Phi(t_{12},t_{23}) \Phi(t_{13},t_{23})^{-1} \Phi(t_{13},t_{12})
=1\ \ \ \mbox{in}\ \wh{U}({\ft}_3),
\Eeq
\Beq\label{6: TH on GRT 3rd eqn - pentagon}
\Phi(t_{13}+t_{23},t_{34})\Phi(t_{12},t_{23}+t_{24})=\Phi(t_{12},t_{23}) \Phi(t_{12}+t_{13},t_{24}+t_{34})\Phi(t_{23},t_{34})
\ \ \ \mbox{in}\ \wh{U}({\ft}_4).
\Eeq
}

\begin{proof} Any automorphism $F$  of the operad  $\wh{\cP a\cC\caD}_\K$ is uniquely determined
by its values on the generators $\bbH_{1,2}$, $\mathbb{X}_{1,2}$ and $\id_{1,2,3}$. These values are not arbitrary, they must respect the relations between the generators. Mac Lane coherence theorem for symmetric monoidal categories implies that all the relations between the generators follow
from the following two ones:
\Beqrn
\mbox{Pentagon}: && d_1\id_{1,2,3} \circ d_3\id_{1,2,3}=d_4\id_{1,2,3}\circ d_2\id_{1,2,3}\circ d_0\id_{1,2,3}\\
\mbox{Hexagon}: && d_1 \mathbb{X}_{1,2}= \id_{1,2,3}\  \circ\  d_0 \mathbb{X}_{1,2}\ \circ\ \left(\mbox{\scriptsize (23)}\ccdot \id_{1,2,3}^{-1}\right)\ \circ \
\mbox{\scriptsize (23)} \ccdot d_3\mathbb{X}_{1,2}\  \circ\ \mbox{\scriptsize (321)}\ccdot \id_{1,2,3}.
\Eeqrn
As an element of $GRT_1$ preserves, by definition, the generators $\bbH_{1,2}$ and $\mathbb{X}_{1,2}$, we conclude that such an element $F$ is uniquely determined by its value on $\id_{1,2,3}$,
$$
F(\id_{1,2,3})\in \Mor_{\wh{\cP a\cC\caD}_\K(3)}\left((12)3,1(23)\right)\simeq \wh{U}({\ft}_3),
$$
which must satisfy the equations,
\Beq\label{6: pentagon for F in GRT}
d_1F(\id_{1,2,3}) \circ d_3F(\id_{1,2,3})
=d_4F(\id_{1,2,3})\circ d_2F(\id_{1,2,3})\circ d_0F(\id_{1,2,3})
\Eeq
and
\Beq\label{6: hexagon for F in GRT}
d_1 X_{1,2}= F(\id_{1,2,3})\  \circ\  d_0 X_{1,2}\ \circ\ \left(\mbox{\scriptsize (23)}\ccdot F(\id_{1,2,3})^{-1}\right)\ \circ \
\mbox{\scriptsize (23)} \ccdot d_3X_{1,2}\  \circ\ \mbox{\scriptsize (321)}\ccdot F(\id_{1,2,3}).
\Eeq
By Corollary \S{\ref{6: corollary on t_3}}
$$
\wh{U}(\ft_3)=\K[[c_3]]\times \K\langle\langle t_{12},t_{23}\rangle\rangle
$$
so that
$$
F(\id_{1,2,3})=f(c_3)\Phi(t_{12},  t_{23})
$$
for some invertible  formal power series,
$$
f(c_3)=a_0+ a_1c_3 + a_2c_3^2+ \ldots ,\ \ a_0\neq 0
$$
in $c_3=t_{12}+ t_{13}+t_{23}$, and some invertible formal power series
$\Phi\in \wh{U}(\ffree_2)$. As $\id_{1,2,3}$ is grouplike in the coalgebra  $\Mor_{\wh{\cP a\cC\caD}_\K(3)}\left((12)3,1(23)\right)$, $\Phi$ is also grouplike, and hence it is of the form
$e^{\phi}$ for some $\phi\in \wh{\ffree}_2$.
\sip

The hexagon equation (\ref{6: hexagon for F in GRT}) implies,
$$
\Phi(t_{12},t_{23}) \Phi(t_{13},t_{23})^{-1} \Phi(t_{13},t_{12})=1, \ \ \ f(c_3)=1.
$$
The hexagon equation is $\bS_3$-equivariant. Applying to that equation permutation $(32)$, or, equivalently,
applying this permutation to the above equation for $\Phi$, we get,
$$
\Phi(t_{13},t_{23}) \Phi(t_{12},t_{23})^{-1} \Phi(t_{12},t_{13})=1
$$
and hence
$$
 \Phi(t_{12},t_{13}) \Phi(t_{13},t_{12})=1.
$$
This proves claims (\ref{6: TH on GRT 1st eqn - skewsym}) and (\ref{6: TH on GRT 2nd eqn - hexagon}).
The last claim follows immediately from the pentagon equation for $F(\id_{1,2,3})$
and the definition of the functors $d_i$. We leave the details to the reader.
\end{proof}

The group $GRT_1$ is prounipotent. Hence it is of the form $\exp(\fg\fr\ft_1)$
where the Lie algebra  $\fg\fr\ft_1$ of $GRT_1$ can be identified with elements $\phi(x,y)\in \wh{\ffree}_2$ which satisfy the equations
\Beq\label{grt1 pentagon}
\phi(t_{13}+t_{23},t_{34}) + \phi(t_{12},t_{23}+t_{24})=\phi(t_{12},t_{23})
+ \phi(t_{12}+t_{13},t_{24}+t_{34}) + \phi(t_{23},t_{34})
\Eeq
in $\wh{\ft}_4$, and
\Beq\label{grt1 hex1}
\phi(x,y)=-\phi(y,x)
\Eeq
\Beq\label{grt1 hex2}
\phi(x,y)+\phi(y,-x-y) + \phi(-x-y,x).
\Eeq
\subsubsection{\bf Exercise} Show that $\phi(x,y)=[x,y]$ satisfies equations
(\ref{grt1 pentagon}) and (\ref{grt1 hex1}), but not (\ref{grt1 hex2}).

\sip

H.\ Furusho \cite{Fu} has proven that $\fg\fr\ft_1$ can be identified with
elements $\phi(x,y)\in \wh{\ffree}_2$ which satisfy the pentagon equation
(\ref{grt1 pentagon}) and do not contain $[x,y]$ as a summand. Put another,
in this case hexagon equations  (\ref{grt1 hex1}) and (\ref{grt1 hex2})
follow from (\ref{grt1 pentagon}).

\subsubsection{\bf Proposition} {\em Let $\mathsf{Simp}^\bu(\ft)$ be the cosimplicial complex of the operad of infinitesimal braids (see ({\ref{6: simp t}})). Then
$$
H^1(\mathsf{Simp}^\bu(\ft))=\fg\fr\ft_1 \oplus \K
$$
where the direct summand $\K$ is generated by $[t_{12},t_{23}]\in \ft_3$.}

\sip

This is a reformulation of Proposition E.1 in \cite{Wi1}. It gives a second proof of H.\ Furusho's result discussed just above.

\subsection{Associators} We have operads $\wh{\cP a\cB}_\K$ and $\wh{\cP a \cC\caD}_\K$ in the category
$\mathsf{Cat(coAss_\K)}$. An  isomorphism of operads (if it exists)
$$
\cA: \wh{\cP a\cB}_\K \lon \wh{\cP a \cC\caD}_\K
$$
is called an {\em associator} \cite{Dr,BN}. It is not hard to show the following

\subsubsection{\bf Theorem} {\em There is a one-to-one correspondence between the set of associators and
grouplike elements, $\Phi(x,y)=e^{\phi(x,y)}$, of the completed filtered Hopf algebra $\wh{U}({\ffree}_2)$ of the free Lie algebra on generators $x$ and $y$ which satisfy the following equations
\Beq\label{6: TH on Assoc 1st eqn - skewsym}
\Phi(x,y)=\Phi^{-1}(y,x),
\Eeq
\Beq\label{6: TH on Assoc 2nd eqn - hexagon}
\Phi(t_{12},t_{23})e^{-\frac{1}{2}t_{23}} \Phi(t_{13},t_{23})^{-1}
e^{-\frac{1}{2}t_{13}}
\Phi(t_{13},t_{12})
e^{-\frac{1}{2}(t_{13}+t_{23})}
=1\ \ \ \mbox{in}\ U(\wh{\ft}_3),
\Eeq
\Beq\label{6: TH on Assoc 3rd eqn - pentagon}
\Phi(t_{13}+t_{23},t_{34})\Phi(t_{12},t_{23}+t_{24})=\Phi(t_{12},t_{23}) \Phi(t_{12}+t_{13},t_{24}+t_{34})\Phi(t_{23},t_{34})
\ \ \ \mbox{in}\ \wh{U}({\ft}_4).
\Eeq
}

It is much harder to show that for any field $\K$ of characteristic zero the set of associators,
$$
\mathsf{As}_\K:=\{\wh{\cP a\cB}_\K \rar \wh{\cP a \cC\caD}_\K\}
$$
is non-empty, see \cite{Dr} for the proof.
The equations defining associators are algebraic, but the only two explicit solutions we know involve transcendental methods in their constructions.

\subsubsection{\bf Example: Knizhnik-Zamolodchikov associator}
The first associator  was constructed by V.\ Drinfeld  with the
help of the monodromy of the flat Knizhnik-Zamolodchikov connection
and is known nowadays as  the Knizhnik-Zamolodchikov
associator $\Phi_{KZ}(x,y)$; this is a formal power series in two
non-commutative formal variable $x$ and $y$ of the form,
$$
\Phi_{KZ}(x,y) = 1+
\sum_{m,k_1,\ldots,k_{m-1}\in \N_{\geq 1}\atop
k_m\in \N_{\geq 2}}
(-1)^m\zeta (k_1, \ldots , k_m)x^{k_m-1}y \cdots x^{k_1-1}y +
\mbox{regularized terms},
$$
where $\zeta (k_1, \ldots , k_m)$ is the multiple zeta value, the real number
defined by the following converging sum,
$$
\zeta (k_1, \ldots , k_m)=\sum_{0< n_1<\ldots< n_m} \frac{1}{n_1^{k_1} \cdots n_m^{k_m} }.
$$
For $m=1$ this gives us the Riemann zeta function $\zeta(k_1)$, $k_1>1$.
There is an explicit iterative construction of the {\em regularized terms}\, in the above formula in terms of
multiple zeta values so that {\em all}\, coefficients of $\Phi_{KZ}$ are rational linear combinations of that values.

\subsubsection{\bf Example: Alekseev-Torossian associator}
The second explicit associator $\Phi_{AT}(x,y)$ was constructed by A.\ Alekseev and
C.\ Torossian in \cite{AT2} using Fulton-MacPherson's compactified configuration spaces
 of points in the complex plane and the integration theory of singular
 differential forms on semialgebraic chains. It was proved to be a Drinfeld associator in \cite{SW}. H.\ Furusho found  a method in \cite{Fu2} of computing
 coefficients of the formal power series $\Phi_{AT}(x,y)$ in terms
of rational linear combinations of iterated integrals of M.\ Kontsevich weight differential forms associated to Lie graphs (see also a  nice paper \cite{BPP} on the systematic computation of weights of the M.\ Kontsevich formality map).

\bip

{\large
\section{\bf Grothendieck-Teichm\"uller group, graph complexes and T.Willwacher theorems}
}

\mip

\subsection{Operads of graphs}
In this section we consider only graphs $\Ga$ without hairs (see  \S {\ref{4: subsec on graphs}}) which are called from  now on simply {\em graphs}. Recall that the set of vertices of $\Ga$ is denoted by $V(\Ga)$ and the set of edges by $E(\Ga)$. A graph $\Ga$ is called {\em directed}\, if each edge $e=(h,\tau(h))\in E(\Ga)$ comes equipped
with a choice of a direction, that is, with a total ordering of its set $\{h,\tau(h)\}$ of half-edges. Here are a few examples of directed graphs
$$
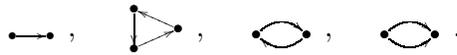
\Ba{c}\resizebox{7mm}{!}{
\xy
 (0,0)*{\bullet}="a",
(8,0)*{\bu}="b",
\ar @{->} "a";"b" <0pt>
\endxy}\Ea, \ \ \
\Ba{c}\resizebox{8mm}{!}{
\xy
 {\ar@{->}(-5,4)*{\bu};(-5,-4)*{\bu}};
   {\ar@{<-}(-5,4)*{\bu};(4,0)*{\bu}};
 {\ar@{<-}(4,0)*{\bu};(-5,-4)*{\bu}};
\endxy}\Ea, \ \ \
\Ba{c}\resizebox{8mm}{!}{\xy
   {\ar@/^0.6pc/(-5,0)*{\bu};(5,0)*{\bu}};
 {\ar@/^0.6pc/(5,0)*{\bu};(-5,0)*{\bu}};
\endxy}\Ea , \ \ \
\Ba{c}\resizebox{8mm}{!}{\xy
   {\ar@/^0.6pc/(-5,0)*{\bu};(5,0)*{\bu}};
 {\ar@/_0.6pc/(-5,0)*{\bu};(5,0)*{\bu}};
\endxy}\Ea .
$$

\sip

Let $G_{n,l}$ be the set of directed graphs $\Ga$ with $n$ vertices and $l$  edges such that
some bijections $V(\Ga)\rar [n]$ and $E(\Ga)\rar [l]$ are fixed, i.e.\ every edge and every vertex of $\Ga$ is marked. The permutation group $\bS_l$ acts on $G_{n,l}$ by relabeling the edges so that,
for each any integer $d$, it makes sense to consider a collection of $\bS_n$-modules,
$$
\caD\cG ra_{d}=\left\{\caD\cG ra_d(n):= \prod_{l\geq 0} \K \langle G_{n,l}\rangle \ot_{ \bS_l}  \sgn_l^{\ot |d-1|} [l(d-1)]   \right\}_{n\geq 1}
$$
with the group $\bS_n$ acting on $\caD\cG ra_d(n)$ by relabeling the vertices.
This is a $\Z$-graded vector space obtained by assigning to each edge of a generating graph from $G_{n,l}$ the homological degree $1-d$. Note that if $d$ is even, then each graph  $\Ga \in \caD\cG ra_d(n)$ is assumed to come equipped with a choice of ordering of edges up to an even permutation (an odd permutation acts as the multiplication by $-1$). In particular, the graph $\Ba{c}\resizebox{7mm}{!}{\xy
(-5,2)*{^1},
(5,2)*{^2},
   {\ar@/^0.6pc/(-5,0)*{\bu};(5,0)*{\bu}};
 {\ar@/^{-0.6pc}/(-5,0)*{\bu};(5,0)*{\bu}};
\endxy}
\Ea\in \caD\cG ra_{d\in 2\Z}(2)$ vanishes identically as it admits an automorphism which changed the ordering of edges by an odd permutation, i.e.\ it is equal to minus itself.


\sip

This $\bS$-module
 is an operad with respect to the following operadic composition,
$$
\Ba{rccc}
\circ_i: &  \caD\cG ra_d(n) \times \caD\cG ra_d(m) &\lon & \caD\cG ra_d(m+n-1),  \ \ \forall\ i\in [n]\\
         &       (\Ga_1, \Ga_2) &\lon &      \Ga_1\circ_i \Ga_2,
\Ea
$$
where  $\Ga_1\circ_i \Ga_2$ is defined by substituting the graph $\Ga_2$ into the $i$-labeled vertex of $\Ga_1$ and taking a sum over all possible re-attachments of dangling edges (attached before to $v_i$) to the vertices of $\Ga_2$. For example,
$$
\Ba{c}\resizebox{10mm}{!}{
\xy
(-5,2)*{^1},
(5,2)*{^2},
   {\ar@/^0.6pc/(-5,0)*{\bu};(5,0)*{\bu}};
 {\ar@/^0.6pc/(5,0)*{\bu};(-5,0)*{\bu}};
\endxy}
\Ea
 \ \ \circ_1\
 \Ba{c}\resizebox{4mm}{!}{
\xy
(-2,7)*{^1},
(-2,0)*{_2},
 {\ar@{->}(0,7)*{\bu};(0,0)*{\bu}};
\endxy}\Ea
=
 \Ba{c}\resizebox{12mm}{!}{
\xy
(-7,7)*{^1},
(-7,0)*{_2},
(5,2)*{^3},
 {\ar@{->}(-5,7)*{\bu};(-5,0)*{\bu}};
   {\ar@/^0.6pc/(-5,0)*{\bu};(5,0)*{\bu}};
 {\ar@/^0.6pc/(5,0)*{\bu};(-5,0)*{\bu}};
\endxy}\Ea
+
 \Ba{c}\resizebox{12mm}{!}{
\xy
(-7,0)*{^1},
(-7,-7)*{_2},
(5,2)*{^3},
 {\ar@{->}(-5,0)*{\bu};(-5,-7)*{\bu}};
   {\ar@/^0.6pc/(-5,0)*{\bu};(5,0)*{\bu}};
 {\ar@/^0.6pc/(5,0)*{\bu};(-5,0)*{\bu}};
\endxy}\Ea
+
 \Ba{c}\resizebox{11mm}{!}{
\xy
(-7,4)*{^1},
(-7,-4)*{_2},
(5,2)*{^3},
 {\ar@{->}(-5,4)*{\bu};(-5,-4)*{\bu}};
   {\ar@{->}(-5,4)*{\bu};(5,0)*{\bu}};
 {\ar@{->}(5,0)*{\bu};(-5,-4)*{\bu}};
\endxy}\Ea
+
 \Ba{c}\resizebox{11mm}{!}{
\xy
(-7,4)*{^1},
(-7,-4)*{_2},
(5,2)*{^3},
 {\ar@{->}(-5,4)*{\bu};(-5,-4)*{\bu}};
   {\ar@{<-}(-5,4)*{\bu};(5,0)*{\bu}};
 {\ar@{<-}(5,0)*{\bu};(-5,-4)*{\bu}};
\endxy}
\Ea
$$
Note that for $d\in 2\Z$ one has to choose of an ordering of edges in each graph summand of the composition $\Ga_1\circ_i \Ga_2$, and there is a canonical way to do by, roughly speaking, placing all the edges of $\Ga_1$ in front of the edges of $\Ga_2$.

\sip

One can define an ``undirected version" of the operad above by noticing that
the group $\bS_l \ltimes
 (\bS_2)^l$ acts on set $G_{n,l}$ of directed labelled graphs by relabeling the
edges and reversing the directions of the edges. Hence
for each fixed integer $d\in \Z$, one can consider a collection of $\bS_n$-modules,
$$
\cG ra_d(n):= \prod_{l\geq 0} \K \langle G_{n,l}\rangle \ot_{ \bS_l \ltimes
 (\bS_2)^l}  \sgn_l^{|d|}\ot \sgn_2^{\ot l|d-1|} [l(d-1)]
$$
where the group $\bS_n$ acts by relabeling the vertices. For $d$ even elements
$\cG ra_d(n)$ can be understood as {\em undirected}\, graphs, e.g.
 $$
\Ba{c}\resizebox{7mm}{!}{
\xy
(0,2)*{^1},
(8,2)*{^2},
 (0,0)*{\bullet}="a",
(8,0)*{\bu}="b",
\ar @{-} "a";"b" <0pt>
\endxy}\Ea \in \cG ra_{d\in 2\Z}(2), \ \ \
\Ba{c}\resizebox{13mm}{!}{
\begin{tikzpicture}[baseline=-.65ex]
\node[int](v) at (0,0){};
\node[int](x) at (1,0){};
\draw (v) edge [bend left] (x) node[above] {$_1$}  ;
\draw (x) edge [bend left] (v) node[above] {$_2$} ;
\end{tikzpicture}
}\Ea\equiv 0, \ \ \
\Ba{c}\resizebox{8mm}{!}{
\xy
(-7,4)*{^1},
(-7,-4)*{_2},
(4,2)*{^3},
 {\ar@{-}(-5,4)*{\bu};(-5,-4)*{\bu}};
   {\ar@{-}(-5,4)*{\bu};(4,0)*{\bu}};
 {\ar@{-}(4,0)*{\bu};(-5,-4)*{\bu}};
\endxy}\Ea\in \cG ra_{d\in 2\Z}(3)
$$
 while for $d$ odd one has identifications of the type
$$
\Ba{c}\resizebox{9mm}{!}{\xy
(-5,2)*{^1},
(5,2)*{^2},
   {\ar@/^0.6pc/(-5,0)*{\bu};(5,0)*{\bu}};
 {\ar@/^{0.6pc}/(5,0)*{\bu};(-5,0)*{\bu}};
\endxy}
\Ea
=-
\Ba{c}\resizebox{9mm}{!}{\xy
(-5,2)*{^1},
(5,2)*{^2},
   {\ar@/^0.6pc/(-5,0)*{\bu};(5,0)*{\bu}};
 {\ar@/^{-0.6pc}/(-5,0)*{\bu};(5,0)*{\bu}};
\endxy}
\Ea\ \ \text{in}\ \ \cG ra_{d\in 2\Z+1}(2).
$$
\subsubsection{\bf Exercise} Show that graphs in $\cG ra_{d\in 2\Z}$ which have multiple edges vanish identically.

\subsubsection{\bf Remark} The integer parameter $d$ has sometimes a clear geometric meaning. With operads of graphs $\caD \cG ra_d$ and $\cG ra_d$ for $d\geq 2$ one can associate de Rham field theories on various compactified configuration spaces of points in (subspaces of) $\R^d$ , the most prominent of which is the Kontsevich formality theory \cite{Ko2} for $d=2$. Other examples of $d=2$ field theories have been studied in, e.g., \cite{Me4,Sh, Wi4}; for an example of a highly non-trivial $d=3$ de Rham field theory we refer to \cite{MW3}. It seems that to get {\em non-trivial}\, (in the sense that $GRT_1$ plays a classifying role, see \S 8) de Rham field theories for $d\geq 4$ one has to work with operads/complexes of {\em multi-oriented}\, graphs which are discussed below.

\subsection{M.\ Kontsevich graph complexes and T.\ Willwacher theorems}
As discussed in \S {\ref{4: subset on Kapr-Manin}} above, for any operad $\f=\{\f(n)\}_{n\geq 1}$  in the category of graded vector spaces
the linear map
$$
\Ba{rccc}
[\ ,\ ]:&  \f^{tot} \ot \f^{tot} & \lon & \f^{tot}\\
& (a\in \cP(n), b\in \cP(m)) & \lon &
[a, b]:= \sum_{i=1}^n a\circ_i b - (-1)^{|a||b|}\sum_{i=1}^m b\circ_i a\ \in \cP(m+n-1)
\Ea
$$
makes a graded vector space
$
\f^{tot}:= \prod_{n\geq 1}\cP(n)$
into a Lie algebra \cite{KM}; moreover, these brackets induce a Lie algebra structure on the subspace
of invariants
$
\f^{tot}_\bS:=  \prod_{n\geq 1}\f(n)^{\bS_n}$. For any $d\in \Z$ one can consider a degree shifted operad $\f\{d\}$ (see \S {\ref{4: subsec on degree shoft of operads}}) and conclude
that for any operad $\f$ and any integer $d$ the associated graded vector space
$$
\f\{d\}^{tot}_\bS\simeq \prod_{n\geq 1} \f(n)\ot_{\bS_n} \sgn_n^{\ot |d|}[d(1-n)]
$$
is canonically a Lie algebra\footnote{In fact, this Lie algebra controls the deformation theory
of the trivial morphism of operads $0: \caL ie\{1-d\} \rar \f$, where $\caL ie$ is the operad of Lie algebras. This observation gives us a very useful approach to the Kontsevich graph complexes and their ribbon graphs generalizations (see \cite{Wi1,MW}) which we skip in this survey.}.
 In particular,
the graded vector spaces (the ``directed full graph complex" and, respectively, the
``full graph complex")
$$
\mathsf{dfGC}_{d}:= \prod_{n\geq 1} \caD \cG ra_{d}(n)\ot_{\bS_n} \sgn_n^{\ot |d|} [d(1-n)]
$$
and
$$
\mathsf{fGC}_{d}:= \prod_{n\geq 1} \cG ra_{d}(n)\ot_{\bS_n} \sgn_n^{\ot |d|} [d(1-n)]
$$
are Lie algebras with Lie brackets given by the substitution of a graph into a vertex of another graph as explained above.  Note that the homological degree of graph $\Ga$ from $\mathsf{dfGC}_{d}$ or $\fGC_d$ is given by
$$
|\Ga|=d(\# V(\Ga) -1) + (1-d) \# E(\Ga).
$$

 \sip

 A generator $\Ga\in \mathsf{dfGC}_{d}$ can be understood as a {\em directed}\, graph (with {\em no}\ labeling of vertices or edges) equipped with an orientation $or$ which is, by definition,
a unital vector in the following 1-dimensional Euclidian space
$$
\R_\Ga:=\left\{\Ba{ll}  \wedge^{\# E(\Ga)} \K[E(\Ga)] & \text{if}\ d\in 2\Z\\
\wedge^{\# V(\Ga)} \R[V(\Ga)] & \text{if}\ d\in 2\Z+1
\Ea
\right.
$$
Every directed graph has precisely two possible orientations, $or$ and $or^{opp}$, and one identifies $(\Ga,or)=-(\Ga,or^{opp})$. We often abbreviate $(\Ga,or)$ to $\Ga$.

\sip

 A generator $\Ga\in \mathsf{fGC}_{d}$ can be understood as an {\em undirected}\, graph (with {\em no}\ labeling of vertices or edges) equipped with an orientation $or$ which for $d$ even is the same as defined just above, while for $d$ odd it is defined as a unital vector in the 1-dimensional Euclidian vector space
$$
\R_\Ga:= \wedge^{\# V(\Ga)} \R[V(\Ga)]\ot \bigotimes_{e=(h,\tau(h))\in E(\Ga)} \wedge^2 \R[h,\tau(h)]
$$
where each edge is interpreted as the orbit consisting of two half-edges under the involution $\tau$ (see \S {\ref{4: subsec on graphs}}). Put another way, for $d$ even an orientation $or$ is a choice of ordering of edges (up to an even permutation), while for $d$ odd $or$ is a choice of ordering of vertices (up to an even permutation) and a choice of direction on each edge (up to a flip and multiplication by $-1$).
Every undirected graph has precisely two possible orientations, $or$ and $or^{opp}$, and one identifies $(\Ga,or)=-(\Ga,or^{opp})$.

\subsubsection{\bf Exercise} Show that the graph
$$
\ga_0:= 
\Ba{c}\resizebox{7mm}{!}{\xy
%
 (0,1)*{\bullet}="a",
(8,1)*{\bu}="b",
\ar @{-} "a";"b" <0pt>
\endxy}\Ea \in \fGC_d \ \ \ \text{or}\ \ \ \
\ga_0:=
\Ba{c}\resizebox{7mm}{!}{\xy
%
 (0,1)*{\bullet}="a",
(8,1)*{\bu}="b",
\ar @{->} "a";"b" <0pt>
\endxy}\Ea \in \dfGC
$$
is a Maurer-Cartan element, i.e.\ it has degree $1$ and satisfies $[\ga_0,\ga_0]=0$

\sip

Hence Lie algebras $\fGC_d$ and, respectively, $\dfGC_d$ come equipped with a compatible differential
$$
\delta \Ga:=[\Ga, \ga_0]= \sum_{v\in V(\Ga)} \left(\delta_v'\Ga -(-1)^{|\Ga|} \delta_v''\Ga\right)
$$
where $\delta_v'$ splits the vertex $v$  of $\Ga$ into two vertices
$$
\delta'_v  \Ba{c}\resizebox{10mm}{!}{\xy
%
 (0,0)*{\bullet}="a",
(-8,0)*{}="1",
(8,0)*{}="2",
(-4,6)*{}="3",
(-4,-6)*{}="4",
(4,6)*{}="5",
(4,-6)*{}="6",
\ar @{-} "a";"1" <0pt>
\ar @{-} "a";"2" <0pt>
\ar @{-} "a";"3" <0pt>
\ar @{-} "a";"4" <0pt>
\ar @{-} "a";"5" <0pt>
\ar @{-} "a";"6" <0pt>
\endxy}
\Ea=\sum_{H(v)=I'\sqcup I''\atop
\# I',\#I''\geq 0}
\Ba{c}\resizebox{20mm}{!}{\xy
%
(0,-10)*{\underbrace{\ \ \ \ \  \ \ \ \ \ }_{I'\ \text{half-edges}}},
(13,10)*{\overbrace{ \ \ \  \ \ \ \ \ }^{I''\ \text{half-edges}}},
 (0,0)*{\bullet}="a",
 (8,0)*{\bullet}="b",
(-8,0)*{}="1",
(16,0)*{}="2",
(-4,6)*{}="3",
(-4,-6)*{}="4",
(12,6)*{}="5",
(4,-6)*{}="6",
\ar @{-} "a";"b" <0pt>
\ar @{-} "a";"1" <0pt>
\ar @{-} "b";"2" <0pt>
\ar @{-} "a";"3" <0pt>
\ar @{-} "a";"4" <0pt>
\ar @{-} "b";"5" <0pt>
\ar @{-} "a";"6" <0pt>
\endxy}
\Ea
\ \ \text{and resp.}\ \
\delta_v'  \Ba{c}\resizebox{10mm}{!}{\xy
%
 (0,0)*{\bullet}="a",
(-8,0)*{}="1",
(8,0)*{}="2",
(-4,6)*{}="3",
(-4,-6)*{}="4",
(4,6)*{}="5",
(4,-6)*{}="6",
\ar @{->} "a";"1" <0pt>
\ar @{->} "a";"2" <0pt>
\ar @{<-} "a";"3" <0pt>
\ar @{->} "a";"4" <0pt>
\ar @{<-} "a";"5" <0pt>
\ar @{<-} "a";"6" <0pt>
\endxy}
\Ea:=\sum_{H(v)=I'\sqcup I''\atop
\# I',\#I''\geq 0}
\Ba{c}\resizebox{20mm}{!}{\xy
%
(0,-10)*{\underbrace{\ \ \ \ \  \ \ \ \ \ }_{I'\ \text{half-edges}}},
(13,10)*{\overbrace{ \ \ \  \ \ \ \ \ }^{I''\ \text{half-edges}}},
 (0,0)*{\bullet}="a",
 (8,0)*{\bullet}="b",
(-8,0)*{}="1",
(16,0)*{}="2",
(-4,6)*{}="3",
(-4,-6)*{}="4",
(12,6)*{}="5",
(4,-6)*{}="6",
\ar @{->} "a";"b" <0pt>
\ar @{->} "a";"1" <0pt>
\ar @{->} "b";"2" <0pt>
\ar @{<-} "a";"3" <0pt>
\ar @{->} "a";"4" <0pt>
\ar @{<-} "b";"5" <0pt>
\ar @{<-} "a";"6" <0pt>
\endxy}
\Ea
$$
 connected by a (directed) edge and redistributes the attached half-edges to $v$ (if any) along the two new vertices in all possible ways, while $\delta_v''$ attaches to $v$ a new univalent vertex
$$
\delta'_v  \Ba{c}\resizebox{10mm}{!}{\xy
%
 (0,0)*{\bullet}="a",
(-8,0)*{}="1",
(8,0)*{}="2",
(-4,6)*{}="3",
(-4,-6)*{}="4",
(4,6)*{}="5",
(4,-6)*{}="6",
\ar @{-} "a";"1" <0pt>
\ar @{-} "a";"2" <0pt>
\ar @{-} "a";"3" <0pt>
\ar @{-} "a";"4" <0pt>
\ar @{-} "a";"5" <0pt>
\ar @{-} "a";"6" <0pt>
\endxy}
\Ea=2
\Ba{c}\resizebox{10mm}{!}{\xy
%
(7,5)*{\bullet}="b",
 (0,0)*{\bullet}="a",
(-8,0)*{}="1",
(7.8,-1)*{}="2",
(-5,5.6)*{}="3",
(-4,-6)*{}="4",
(1.5,7.1)*{}="5",
(4,-6)*{}="6",
\ar @{-} "a";"1" <0pt>
\ar @{-} "a";"2" <0pt>
\ar @{-} "a";"3" <0pt>
\ar @{-} "a";"4" <0pt>
\ar @{-} "a";"5" <0pt>
\ar @{-} "a";"6" <0pt>
\ar @{-} "a";"b" <0pt>
\endxy}\Ea
\ \ \ \text{and resp.}\ \ \
\delta_v'  \Ba{c}\resizebox{10mm}{!}{\xy
%
 (0,0)*{\bullet}="a",
(-8,0)*{}="1",
(8,0)*{}="2",
(-4,6)*{}="3",
(-4,-6)*{}="4",
(4,6)*{}="5",
(4,-6)*{}="6",
\ar @{->} "a";"1" <0pt>
\ar @{->} "a";"2" <0pt>
\ar @{<-} "a";"3" <0pt>
\ar @{->} "a";"4" <0pt>
\ar @{<-} "a";"5" <0pt>
\ar @{<-} "a";"6" <0pt>
\endxy}
\Ea:=
\Ba{c}\resizebox{10mm}{!}{\xy
(7,5)*{\bullet}="b",
 (0,0)*{\bullet}="a",
(-8,0)*{}="1",
(7.8,-1)*{}="2",
(-5,5.6)*{}="3",
(-4,-6)*{}="4",
(1.5,7.1)*{}="5",
(4,-6)*{}="6",
\ar @{->} "a";"1" <0pt>
\ar @{->} "a";"2" <0pt>
\ar @{<-} "a";"3" <0pt>
\ar @{->} "a";"4" <0pt>
\ar @{<-} "a";"5" <0pt>
\ar @{->} "a";"6" <0pt>
\ar @{->} "a";"b" <0pt>
\endxy}\Ea
\ +\
\Ba{c}\resizebox{10mm}{!}{\xy
(7,5)*{\bullet}="b",
 (0,0)*{\bullet}="a",
(-8,0)*{}="1",
(7.8,-1)*{}="2",
(-5,5.6)*{}="3",
(-4,-6)*{}="4",
(1.5,7.1)*{}="5",
(4,-6)*{}="6",
\ar @{->} "a";"1" <0pt>
\ar @{->} "a";"2" <0pt>
\ar @{<-} "a";"3" <0pt>
\ar @{->} "a";"4" <0pt>
\ar @{<-} "a";"5" <0pt>
\ar @{->} "a";"6" <0pt>
\ar @{<-} "a";"b" <0pt>
\endxy}\Ea
$$
Note that if the vertex $v$ has at least one half-edge, then the two terms with univalent vertices in the sum $\delta_v'\Ga$ cancel out with the two terms
in $\delta_v''\Ga$ so that in most cases one can ignore in the above formula for $\delta$ the second term $\delta_v''$ and assume in the first term $\delta_v'$
that the summation runs over decompositions $H(v)=I'\sqcup I''$ with $\# I'\geq 1$ and $\# I''\geq 1$. For example,
the graph $\ga_0$ is a cycle, $\delta\ga_0=[\ga_0,\ga_0]=0$, but the associated cohomology class is trivial, $\ga_0=-\delta(\bu)$.

\sip

The dg Lie algebra $(\fGC_d,\delta)$ was introduced by M.Kontsevich in \cite{Ko1} in his attempt to calculate obstructions to universal quantizations of Poisson structures.

\sip

Both dg Lie algebras contain dg Lie subalgebras, $\mathsf{fcGC}_d$ and $\mathsf{dfcGC}_d$, generated by {\em connected}\, graphs.  There is an isomorphism of complexes,
$$
\fGC_d=\odot^{\bu\geq 1} \left(\mathsf{fcGC}_d [-d]\right)[d],\ \ \
\dfGC_d=\odot^{\bu\geq 1} \left(\mathsf{dfcGC}_d [-d]\right)[d]
$$
so that at the cohomology level
$$
H^\bu(\fGC_d)=\odot^{\bu\geq 1} \left(H^\bu(\mathsf{fcGC}_d) [-d]\right)[d],\ \ \
H^\bu(\dfGC_d)=\odot^{\bu\geq 1} \left(H^\bu(\mathsf{dfcGC}_d) [-d]\right)[d]
$$
Hence there is no loss of important information when working with connected graphs only.

\sip

Moreover, there is a monomorphism of graph complexes
\Beq\label{6: fcGC to dfcGC}
\mathsf{fcGC}_d \lon \mathsf{dfcGC}_d
\Eeq
which sends an undirected graph $\Ga$ into a sum of graphs
 obtained by interpreting
each edge as the sum of edges in both directions. It was proven in
\cite{Wi1}(Appendix K) that this map is a {\em quasi-isomorphism}. Hence there is no need to study the complex $\dfGC_d$ by itself\footnote{However the complex $\mathsf{dfcGC}_d$ contains two very important subcomplexes of {\em oriented}\, graphs and  of {\em sourced}\, graphs which are of great significance  in applications (see below).} so that we continue in this subsection discussing only the complex $\mathsf{fcGC}_d$. That complex splits as a direct sum
$$
\mathsf{fcGC}_d=\mathsf{fcGC}^{\leq 1}_d \oplus \mathsf{fcGC}^2_d \oplus \GC_d^\circlearrowright
$$
where $\mathsf{fcGC}^{\leq 1}$ is spanned by connected graphs having at least one vertex of valency $\leq 1$, $\mathsf{fcGC}^2_d$ is spanned by connected graphs with no vertices of valency $\leq 1$ but with at least one vertex of valency $2$, and $\GC_d^\circlearrowright$ is generated by graphs with each vertex of valency $\geq 3$. It was proven in \cite{Wi1} that the complex
$\mathsf{fcGC}^{\leq 1}_d$ is acyclic while
$$
H^\bu(\mathsf{fcGC}_{d}^2) = \bigoplus_{j\geq 1\atop j\equiv 2d+1 \mod 4} \K[d-j],
$$
where the summand $\K[d-j]$ is generated by the polytope with $j$ vertices,
that is, a connected graph with $j$ bivalent vertices, e.g.\
$\Ba{c}\resizebox{6mm}{!}{\xy
 (0,0)*{\bullet}="a",
(6,0)*{\bu}="b",
(3,5)*{\bu}="c",
\ar @{-} "a";"b" <0pt>
\ar @{-} "a";"c" <0pt>
\ar @{-} "c";"b" <0pt>
\endxy}\Ea$ for $j=3$ (in fact this particular graph vanishes identically for $d$ {\em even}\, as in this case it admits an automorphism reversing its orientation).

\sip

The complex $\GC_d^\circlearrowright$ contains a subcomplex $\GC_d$ spanned by graphs with no {\em loops}, that is, edges attached to one and the same vertex, and the inclusion is a quasi-isomorphism \cite{Wi1},
$$
H^\bu(\GC_d)=H^\bu(\GC_d^\circlearrowright).
$$
One of the major results in \cite{Wi1} is the following

\subsubsection{\bf T.Willwacher Theorem} {\em There is an isomorphism of Lie algebras,
$$
H^0(\GC_2)=\grt_1,
$$
where $\grt_1$ is the Lie algebra of the Grothendieck-Teichm\"uller group $GRT_1$ (see \S 6).
Moreover, $H^{\bu<0}(\GC_2)=0$.}

\sip

There is an explicit construction of infinitely many cohomology classes $\{[\mathfrak{w}_{2n+1}]\}_{n\geq 1}$ in  $H^0(\GC_2)$, more precisely, of their cycle representatives $\{\mathfrak{w}_{2n+1}\}_{n\geq 1}$  in $\GC_2$. The first two classes can be given explicitly as follows
$$
\mathfrak{w}_3=
 \Ba{c}\resizebox{12mm}{!}{
\xy
 (0,0)*{\bullet}="a",
(0,8)*{\bullet}="b",
(-7.5,-4.5)*{\bullet}="c",
(7.5,-4.5)*{\bullet}="d",
\ar @{-} "a";"b" <0pt>
\ar @{-} "a";"c" <0pt>
\ar @{-} "b";"c" <0pt>
\ar @{-} "d";"c" <0pt>
\ar @{-} "b";"d" <0pt>
\ar @{-} "d";"a" <0pt>
\endxy}
\Ea
$$
$$
\mathfrak{w}_5=
 \Ba{c}\resizebox{13mm}{!}{
\xy
 (0,0)*{\bullet}="0",
(0,8)*{\bullet}="1",
(-8,3)*{\bullet}="5",
(8,3)*{\bullet}="2",
(-5,-7)*{\bullet}="4",
(5,-7)*{\bullet}="3",
\ar @{-} "0";"1" <0pt>
\ar @{-} "0";"2" <0pt>
\ar @{-} "0";"3" <0pt>
\ar @{-} "0";"4" <0pt>
\ar @{-} "0";"5" <0pt>
\ar @{-} "1";"2" <0pt>
\ar @{-} "2";"3" <0pt>
\ar @{-} "3";"4" <0pt>
\ar @{-} "4";"5" <0pt>
\ar @{-} "5";"1" <0pt>
\endxy}\Ea
  \ +\ \frac{5}{2}
 \Ba{c}\resizebox{12mm}{!}{
\xy
 (1,0)*{\bullet}="0",
(0,8)*{\bullet}="1",
(-8,3)*{\bullet}="5",
(8,3)*{\bullet}="2",
(-5,-7)*{\bullet}="4",
(5,-7)*{\bullet}="3",
\ar @{-} "0";"1" <0pt>
\ar @{-} "0";"2" <0pt>
\ar @{-} "0";"4" <0pt>
\ar @{-} "1";"4" <0pt>
\ar @{-} "5";"3" <0pt>
\ar @{-} "1";"2" <0pt>
\ar @{-} "2";"3" <0pt>
\ar @{-} "3";"4" <0pt>
\ar @{-} "4";"5" <0pt>
\ar @{-} "5";"1" <0pt>
\endxy}
\Ea
$$

For higher $n$ there is an explicit transcendental formula for the cycles
$\mathfrak{w}_{2n+1}$ given in \cite{RW} which presents each cycle as a linear combination of graphs with $2n+2$ vertices and $4n+2$ edges
\Beq\label{7: w_2m+1 formulae}
\mathfrak{w}_{2n+1}=\sum_{\Ga\in G_{2n+2,4n+2}} c_\Ga \Ga
= \ \ \la_{2n+1}\
 \Ba{c}\resizebox{17mm}{!}{\xy
 (0,0)*{\bullet}="1",
(-11,11)*{\bullet}="2",
(0,15)*{\bullet}="3",
(11,11)*{\bullet}="4",
(15,0)*^+{\ldots}="5",
(11,-11)*{\bullet}="6",
(0,-15)*{\bullet}="7",
(-15,0)*{\bullet}="n",
(-11,-11)*{\bullet}="n-1",
\ar @{-} "1";"2" <0pt>
\ar @{-} "1";"3" <0pt>
\ar @{-} "1";"4" <0pt>
\ar @{-} "1";"5" <0pt>
\ar @{-} "1";"6" <0pt>
\ar @{-} "1";"7" <0pt>
\ar @{-} "1";"n" <0pt>
\ar @{-} "1";"n-1" <0pt>
\ar @{-} "2";"3" <0pt>
\ar @{-} "3";"4" <0pt>
\ar @{-} "4";"5" <0pt>
\ar @{-} "5";"6" <0pt>
\ar @{-} "6";"7" <0pt>
\ar @{-} "7";"n-1" <0pt>
\ar @{-} "n-1";"n" <0pt>
\ar @{-} "n";"2" <0pt>
\endxy}\Ea\  + \ \ldots
\Eeq
where the weights $c_\Ga$ are given by explicit converging integrals
over the configuration space of $2n$ different points in $\C\setminus \{0,1\}$ (see (23) in \cite{RW}). In particular, the coefficients $\la_{2n+1}$ of the wheel-type summands in $\mathfrak{w}_{2n+1}$ are equal to  zeta values $\frac{\zeta(2n+1)}{(2\pi i)^{(2n+1)}}$ up to a non-zero rational factor, cf.\ \cite{Go}. (In fact, from the presence in
$\mathfrak{w}_{2n+1}$ of such wheel-graphs one can conclude almost immediately that these cycles are {\em not}\, boundaries.)

\sip

\subsubsection{\bf Deligne-Drinfeld-Ihara conjecture} {\em The pro-nilpotent Lie algebra $\grt_1$ is isomorphic as a $\Z_{\geq 3}$ graded Lie algebra to the degree
completion of the free Lie algebra generated by
formal variables of degrees $2n+1$, $n\geq 1$.}

\sip
It has been proved by F.\ Brown in \cite{B} that
 $\grt_1$ does contain indeed the above mentioned free Lie algebra.

 \sip

 Thus the graph cohomology classes $[\mathfrak{w}_{2n+1}]$ might generate the whole $\grt_1$.

\sip

Denote by $\GC_d^{\geq 2}$ the dg Lie subalgebra of $\mathsf{fcGC}_d$ spanned by graphs with valency of each vertex $\geq 2$ and with no loops. Then it follows from the results discussed above that
$$
H^\bu(\mathsf{fcGC}_d)=H^\bu(\GC_d^{\geq 2})= H^\bu(\GC_d) +
\bigoplus_{j\geq 1\atop j\equiv 2d+1 \mod 4} \K[d-j].
$$
There are other differentials on the Kontsevich graph complexes $\GC_d$
which have been studied in \cite{MW0, KWZ1}.

\subsection{Oriented graph complexes} A graph $\Ga$ from the dg Lie algebra $\mathsf{dfGC}_d$
is called {\em oriented}\, if it contains no {\em wheels}, that is, directed paths of edges
forming a closed circle. The subspace
$\mathsf{fGC}_d^{or}\subset \mathsf{dfGC}_d$ spanned by oriented graphs is a dg Lie subalgebra. For example,
$$
 \Ba{c}\resizebox{11mm}{!}{
 \xy
 (0,0)*{\bullet}="a",
(12,0)*{\bu}="b",
(6,10)*{\bu}="c",
\ar @{->} "a";"b" <0pt>
\ar @{->} "a";"c" <0pt>
\ar @{<-} "c";"b" <0pt>
\endxy}\Ea
\in \fGC_d^{or} \ \ \ \ \mbox{while}\ \ \ \
 \Ba{c}\resizebox{11mm}{!}{
 \xy
 (0,0)*{\bullet}="a",
(12,0)*{\bu}="b",
(6,10)*{\bu}="c",
\ar @{->} "a";"b" <0pt>
\ar @{<-} "a";"c" <0pt>
\ar @{<-} "c";"b" <0pt>
\endxy}\Ea
\not\in \fGC_d^{or}.
$$
Let $\mathsf{dfcGC}_d^{or}$ the subcomplex of $\fGC_d^{or}$ spanned by connected graphs, and $\GCor_d\subset \mathsf{dfcGC}_d^{or}$ the subcomplex spanned by graphs
which each vertex at least bivalent and with no bivalent vertices of the form $
\Ba{c}\resizebox{10mm}{!}{\xy
 (0,1)*{}="a",
(7,1)*{\bu}="b",
 (14,1)*{}="c",
\ar @{->} "a";"b" <0pt>
\ar @{->} "b";"c" <0pt>
\endxy}\Ea$. Then \cite{Wi2}
$$
H^\bu(\mathsf{fcGC}_d^{or})=H^\bu(\GCor_d)
$$
so that there no loss of important information when working with $\GCor_d$ solely.

\subsubsection{\bf Theorem \cite{Wi2}}
$
H^\bu(\GCor_{d+1})=H^\bu(\mathsf{GC}_d^{\geq 2}).
$

\sip

Hence one has a remarkable isomorphism of Lie algebras,
$
H^0(\GCor_3)=\fg\fr\ft$. Moreover it follows from this theorem that $H^i(\GCor_3)=0$ for $i\leq -2$ and $H^{-1}(\GCor_3)$
is a 1-dimensional space generated by the graph $
\Ba{c}\resizebox{4mm}{!}{   \xy
   \ar@/^0.6pc/(0,-5)*{\bullet};(0,5)*{\bullet}
   \ar@/^{-0.6pc}/(0,-5)*{\bullet};(0,5)*{\bullet}
 \endxy}\Ea
$. This theorem tells us that in order to have some non-trivial de Rham field theory
in dimension 3 (rather than two), one has to work with {\em oriented}\, graphs rather than with simply directed (or even undirected) graphs.
We refer to \cite{MW3} for an example of such a de Rham field theory
providing us with an {\em explicit}\, universal formula for deformation quantization of Lie bialgebras.
\sip

The original argument in \cite{Wi3} does not give us an explicit relation between
the cohomology groups in the above Theorem.
Let $\widetilde{\GCor_{d+1}}$ and $\widetilde{\GC^{\geq 2}_{d}}$ be the graph complexes dual to $\GCor_{d+1}$ and $\mathsf{GC}_d^{\geq 2}$ respectively. Then one also  has $H^\bu (\widetilde{\GCor_{d+1}})=H^\bu(\widetilde{\mathsf{GC}_d^{\geq 2}})$.
An explicit construction sending a homology class in $H^\bu(\widetilde{\mathsf{GC}_d^{\geq 2}})$ to a cohomology class in
$H^\bu (\widetilde{\GCor_{d+1}})$ has been found in \cite{Z1}. It can be used to find an explicit $3$-dimensional oriented incarnation $\mathfrak{w}_3^{or}\in H^0(\GCor_3)$ of, for example, the tetrahedron cohomology class $\mathfrak{w}_3\in H^0(\GC_2)$ discussed above; $\mathfrak{w}_3^{or}$ is a linear combination of oriented graphs with $7$ vertices and $9$ edges.

\sip

In fact, the complex $\GCor_3$ controls
the deformation theory of the properad of Lie bialgebras so that $GRT_1$ is essentially a symmetry group of that properad \cite{MW2}. In string topology and the theory of moduli spaces of algebraic curves it is important to work with {\em involutive}\, Lie bialgebras, and their deformation theory is controlled by the following deformation of the complex $\GCor_3$.

\subsubsection{\bf A deformation of $\GCor_3$}
Consider a Lie algebra $(\GCor_3[[\hbar]], [\ ,\ ])$, where
$\GCor_3[[\hbar]]$ is the (completed) topological vector space spanned by formal power series in a formal parameter $\hbar$ of homological degree $2$, and $[\ ,\ ]$ are the Lie brackets obtained from the standard ones in
$\GCor_d$ by the  $\K[[\hbar]]$-linearity. It was shown in \cite{CMW} that the formal power series
$$
\ga_\hbar:= \sum_{k=1}^\infty \hbar^{k-1} \underbrace{
\Ba{c}\resizebox{6mm}{!}  {\xy
(0,0)*{...},
   \ar@/^1pc/(0,-5)*{\bullet};(0,5)*{\bullet}
   \ar@/^{-1pc}/(0,-5)*{\bullet};(0,5)*{\bullet}
   \ar@/^0.6pc/(0,-5)*{\bullet};(0,5)*{\bullet}
   \ar@/^{-0.6pc}/(0,-5)*{\bullet};(0,5)*{\bullet}
 \endxy}
 \Ea}_{k\ \mathrm{edges}}
$$
is a Maurer-Cartan element in the Lie algebra $(\mathsf{f}\sG\sC_3^{or}[[\hbar]], [\ ,\ ])$ and hence makes the latter into a {\em differential}\, Lie algebra with the differential
$
\delta_\hbar \Ga:=[\Ga, \ga_\hbar].
$
It was proven in \cite{CMW} that
 $H^0( \GC_3^{or}[[\hbar]], \delta_\hbar)\simeq
H^0(\GC_3^{or},\delta)\simeq \grt_1$ as Lie algebras. Moreover, $H^i( \GC_3^{or}[[\hbar]], \delta_\hbar)=0$
for all $i\leq -2$ and $H^{-1}( \GC_3^{or}[[\hbar]], \delta_\hbar)$ is a 1-dimensional vector space  class generated by the formal power series
$
\sum_{k=2}^\infty (k-1)\hbar^{k-2}\underbrace{
\Ba{c}\resizebox{6mm}{!}  
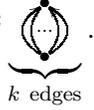
{\xy
(0,0)*{...},
   \ar@/^1pc/(0,-5)*{\bullet};(0,5)*{\bullet}
   \ar@/^{-1pc}/(0,-5)*{\bullet};(0,5)*{\bullet}
   \ar@/^0.6pc/(0,-5)*{\bullet};(0,5)*{\bullet}
   \ar@/^{-0.6pc}/(0,-5)*{\bullet};(0,5)*{\bullet}
 \endxy}
 \Ea}_{k\ \mathrm{edges}}.
$

\subsection{Sourced graphs}  A graphs $\Ga$ from the dg Lie algebra of connected oriented graphs $\mathsf{dfcGC}_d$
is called {\em sourced}\, if it contains at least one vertex $v$ of valency $\geq 2$ with no incoming edges, i.e.\ at least one vertex of the form
$$
v=\Ba{c}\resizebox{12mm}{!}
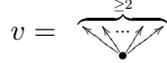
{\xy
(0,9)*{\overbrace{\ \ \ \ \ \ \ \ \ \ \ \ \ \ }^{\geq 2} },
(0,5)*{...},
 (0,0)*{\bullet}="a",
(-8,6)*{}="1",
(-4,6)*{}="2",
(4,6)*{}="3",
(8,6)*{}="4",
\ar @{->} "a";"1" <0pt>
\ar @{->} "a";"2" <0pt>
\ar @{->} "a";"3" <0pt>
\ar @{->} "a";"4" <0pt>
\endxy}
\Ea
$$
Such vertex is called a {\em source}.
The subspace $\mathsf{sfGC}_d\subset \mathsf{dfcGC}_d$ spanned by sourced graphs is a dg Lie subalgebra.
Consider a smaller subspace  $\mathsf{sGC}_d\subset \mathsf{sfGC}_d$ spanned by sourced graphs with all vertices at least bivalent
and with no bivalent vertices of the form $
\Ba{c}\resizebox{10mm}{!}{\xy
 (0,1)*{}="a",
(7,1)*{\bu}="b",
 (14,1)*{}="c",
\ar @{->} "a";"b" <0pt>
\ar @{->} "b";"c" <0pt>
\endxy}\Ea$. It was proven by in \cite{Z1} that
$$
H^\bu(\mathsf{sfGC}_{d+1})=H^\bu(\mathsf{sGC}_{d+1})=H^\bu(\GC_d^{\geq 2})
$$
so that one gets one more incarnation of the Grothendieck-Teichm\"uller Lie algebra in dimension $3$
$$
\grt_1=H^0(\mathsf{sGC}_{3})
$$
in terms of sourced graphs. Similarly one can work with directed graphs with at least one {\em target}\, vertex
of the form
$\Ba{c}\resizebox{10mm}{!}{\xy
(0,-9)*{\underbrace{\ \ \ \ \ \ \ \ \ \ \ \ \ \ }_{\geq 2} },
(0,-5)*{...},
 (0,0)*{\bullet}="a",
(-8,-6)*{}="1",
(-4,-6)*{}="2",
(4,-6)*{}="3",
(8,-6)*{}="4",
\ar @{<-} "a";"1" <0pt>
\ar @{<-} "a";"2" <0pt>
\ar @{<-} "a";"3" <0pt>
\ar @{<-} "a";"4" <0pt>
\endxy}
\Ea$.

\subsection{Multi-oriented graph complexes}\label{2: subsec on DRGC}
The Grothendieck-Teichm\"uller Lie algebra appears naturally as a single entity in dimension $d=2$ as the cohomology group $H^0(\GC_2)$ of the ``ordinary" graph complex.
In dimension 3 it reappears as  the cohomology group $H^0(\GCor_3)=H^0(\mathsf{sGC}_{3})$ of the {\em oriented}\, or {\em
sourced}\, graph complex, i.e.\ we have a sequence of isomorphism of Lie algebras,
$$
\grt_1=H^0(\GC_2)=H^0(\GCor_3)=H^0(\mathsf{sGC}_{3}).
$$
What kind of graph complex $\GC_d^?$ gives us an incarnation of $\grt_1$ in dimension $d=4, 5,\ldots$,
$$
\grt_1=H^0(\GC_2)=H^0(\GCor_3)=H^0(\mathsf{sGC}_{3})=H^0(\GC_4^?)=H^0(\GC_5^?)=\ldots?
$$
The answer to that question was found by Marko \v Zivkovi\' c in \cite{Z1,Z2} who introduced and studied {\em multi-oriented sourced}\, graph complexes.
 Their algebro-geometric interpretation was found in \cite{Me5} --- they control quantizations and deformation theory of
 strong homotopy (even or odd) Lie bialgebra structures, in particular of Poisson structures, on  infinite-dimensional spaces spaces {\em with branes}.

 \sip

 Let us go back for a moment to the Lie algebra of directed connected graphs $\mathsf{dfcGC}_d$ and define, for any integer $k\in \N_{\geq 1}$, its extension $\sfd^k\fcGC_d$ as a $\Z$-graded vector space spanned by directed graphs $\Ga$ with each edge decorated with
 $k-1$ new directions labelled  by $2,3,\ldots, k$, the value $1$ being reserved for the original direction\footnote{It is more suitable to understand new directions as {\em colored}\, ones, say the original direction  $1$ has black colour, the new direction $2$ is red so on, because we do {\em not}\, assume that the set $[k]$ of all directions on each edge is {\em ordered}.  Hence we call often the direction on an edge labeled by integer $l\in [k]$ the $l$-{\em coloured}\, direction. This new decoration of edges has no impact
 on the degree of graphs which is given by the standard formula
 $|\Ga|=d(\# V(\Ga) -1) + (1-d) \# E(\Ga)$, $\forall \Ga\in \sfd^k\fcGC_d$.},
  and each new direction can take only two ``values" --- it can agree or disagree with the original direction,

  $$
\begin{tikzpicture}[baseline=-1ex]
\node[] (a) at (0,-0.05) {};
\node[] (b) at (1.5,-0.05) {};
\draw (a) edge[->] (b);
\end{tikzpicture}
\rightsquigarrow
\Ba{c}
\begin{tikzpicture}[baseline=-2ex]
\draw (0,0) edge[rightblack] node[above] {\ \ \ $\scriptstyle 1$}(0.9,0);
\draw (0.9,0) edge[leftblack] node[above] {\ \ $\scriptstyle 2$ }(1.5,0);
\draw (1.5,0) edge[] node[above] {\ $...$}(1.9,0);
\draw (1.9,0) edge[leftblack] node[above] {\ \ \ $\scriptstyle k$ }(2.3,0);
\draw (2.2,0) edge[](2.7,0);
\end{tikzpicture}
\Ea
$$
Note that the Lie bracket in $\mathsf{dfcGC}_d$ does not change the number and the ``internal structure" of edges in graphs so that exactly the same formula as in $\mathsf{d}\fcGC_d$ gives us a Lie bracket in  $\mathsf{d}^k\fcGC_d$ for any $k\geq 1$. It is not hard to see that the element
$$
\ga_k:=
\Ba{c}
\begin{tikzpicture}[baseline=-2ex]
\draw (0,0) edge[rightblack] node[above] {\ \ \ $\scriptstyle 1$}(0.9,0);
\draw (0.9,0) edge[leftblack] node[above] {\ \ $\scriptstyle 2$ }(1.5,0);
\draw (1.5,0) edge[] node[above] {\ $...$}(1.9,0);
\draw (1.9,0) edge[leftblack] node[above] {\ \ \ $\scriptstyle k$ }(2.3,0);
\draw (2.2,0) edge[](2.7,0);
\end{tikzpicture}\Ea
+
\Ba{c}
\begin{tikzpicture}[baseline=-2ex]
\draw (0,0) edge[rightblack] node[above] {\ \ \ $\scriptstyle 1$}(0.9,0);
\draw (0.9,0) edge[rightblack] node[above] {\ \ $\scriptstyle 2$ }(1.5,0);
\draw (1.5,0) edge[] node[above] {\ $...$}(1.9,0);
\draw (1.9,0) edge[leftblack] node[above] {\ \ \ $\scriptstyle k$ }(2.3,0);
\draw (2.2,0) edge[](2.7,0);
\end{tikzpicture}\Ea
+
\Ba{c}
\begin{tikzpicture}[baseline=-2ex]
\draw (0,0) edge[rightblack] node[above] {\ \ \ $\scriptstyle 1$}(0.9,0);
\draw (0.9,0) edge[leftblack] node[above] {\ \ $\scriptstyle 2$ }(1.5,0);
\draw (1.5,0) edge[] node[above] {\ $...$}(1.9,0);
\draw (1.9,0) edge[rightblack] node[above] {\ \ \ $\scriptstyle k$ }(2.3,0);
\draw (2.2,0) edge[](2.7,0);
\end{tikzpicture}
\Ea
+
\ldots \ \ (2^{k-1}\ \text{terms})
$$
where the summation runs over all possible way to decorate the directed edge with new directions $2,3,\ldots k$,
is a Maurer-Cartan element and hence makes $\mathsf{d}^k\fcGC_d$ into a {\em dg}\, Lie algebra
with the differential $\delta\Ga:=[\Ga,\ga_k]$.
Of course, the case $k=1$ gives nothing new, $\mathsf{d}^1\fcGC_d\equiv \mathsf{dfcGC}_d$. In fact, the general case $k\in \Z_{\geq 2}$ also does not give us immediately anything really new as well: the natural monomorphisms which add to each edge of a graph one extra direction in all possible ways (cf.\ (\ref{6: fcGC to dfcGC})),
$$
\mathsf{fcGC}_d \lon \mathsf{d}^1\fcGC_d \lon \mathsf{d}^2\fcGC_d \lon \mathsf{d}^3\fcGC_d \lon \ldots
$$
are quasi-isomorphisms. However we can consider now lots of interesting subcomplexes in the multidirected graph complex  $(\mathsf{d}^k\fcGC_d, \delta)$. Given any two non-negative integers $p$ and $q$ with $p+q\leq k$,
and any two disjoint subsets $I_p$ and $I_q$ of the set $[k]$ of cardinalities $p$ and $q$ respectively (the particular choice of such subsets plays no role), one can define a subcomplex
$$
\mathsf{s}^p\sfo^q\sfd^k\fcGC_d \subset \mathsf{d}^k\fcGC_d
$$
of {\em $k$-directed $p$-sourced and $q$-oriented graphs}\, as the span
of those  multioriented graphs $\Ga$ which satisfy the conditions
\Bi
\item[(i)] $\Ga$ is {\em sourced}\, with respect to every direction in the subset $I_p\subseteq [k]$, i.e.\ for each direction $c\in I_p$ there is at least one vertex $v\in V(\Ga)$ which has valency $\geq 2$ and no incoming edges with respect to  $c$;
\item[(ii)] $\Ga$ is {\em oriented}\, with respect to every direction in the subset $I_q\subset [k]$, i.e.\ for each direction $c\in I_p$ the graph $\Ga$ has no {\em closed}\, directed paths of edges ({\em wheels})  with respect to $c$.

\Ei
By analogy to the previous examples, there is no loss of generality when working with a subcomplex
$$
\mathsf{s}^p\sfo^q\sfd^k\GC_d\subset \mathsf{s}^p\sfo^q\sfd^k\fcGC_d
$$
generated by graphs whose vertices are at least bivalent as this inclusion is a quasi-isomorphism. The major result of \cite{Z2} is the following

\subsubsection{\bf M.\ \v Zivkovi\'c Theorem} {\em For any integer $k\geq 1$ and any non-negative integers $p$ and $q$ with $p+q\leq k$ there is an isomorphism of cohomology groups}
$$
H^\bu(\GC_d^{\geq 2})=H^\bu(\mathsf{s}^p\sfo^q\sfd^k\GC_{d+p+q}).
$$

\sip

Therefore extra directions which are neither oriented nor sourced can be omitted
as irrelevant, i.e.\  one can set $k=p+q$ without loss of important information  and work only with dg Lie algebras
$$
\mathsf{s}^p\sfo^q\GC_{d+p+q}:=\mathsf{s}^p\sfo^q\sfd^{p+q}\GC_{d+p+q}.
$$
 M.\ \v Zivkovi\'c Theorem gives us infinitely many graph incarnations of the Grothendieck-Teichm\"uller Lie algebra,
$$
\grt_1=H^0(\GC_2)=H^0(\mathsf{s}^p\sfo^q\GC_{2+p+q})\ \ \ \forall\ p,q\in \Z_{\geq 0},
$$
and hence a clue (cf.\ \cite{Me5}) to what could be a non-trivial deformation quantization theory in dimension $d\geq 4$ in which the Grothendieck-Teichm\"uller group plays a classifying role (the cases $d=2$ and $d=3$ are well understood by now, see below).

\subsubsection{\bf Remark on hairy graphs} There is an important version of the graph complexes $\GC_d$ which is based on graphs with hairs \cite{KWZ2,Wi5, Wi6}. These graph complexes control the rational homotopy groups of
of long embeddings (modulo
immersion) of $\R^m$ in $\R^n$. Very recently, their marked version has been used in the study of cohomology groups of moduli spaces $\cM_{g,n}$ of genus $g$ algebraic curves with $n$ punctures \cite{CGP2}.

\bip

{\large
\section{\bf Some applications of the theory of Drinfeld's associators, $GRT_1$\\ and graph complexes}
}

\sip
\subsection{Universal quantizations of Lie bialgebras}
 A  {\em Lie bialgebra}\, is a  graded vector space $V$
equipped with two linear maps,
$$
[\ , \ ]: \wedge^2 V
\rightarrow V \ \ \ , \ \ \vartriangle: V\rightarrow V\wedge V ,
$$
such that the first map
$[\ ,\ ]$ makes  $V$ into a Lie algebra,
the second map
$\vartriangle$ makes $V$ into a Lie coalgebra,  and  the compatibility condition
$$
\vartriangle [a, b] = \sum a_1\otimes [a_2, b] +  [a,
b_1]\otimes b_2 - (-1)^{|a||b|}( [b, a_1]\otimes a_2
+ b_1\otimes [b_2, a]),
$$
holds for any $a,b\in V$ with $\vartriangle a=:\sum a_1\otimes a_2$, $\vartriangle b=:\sum
b_1\otimes b_2$.
This algebraic structure was introduced by V.\ Drinfeld \cite{D1}
in the context of studying universal deformations of the standard
Hopf algebra structure on universal enveloping algebras. More precisely, consider the symmetric tensor algebra $\odot^\bu V$ equipped with the standard graded commutative and graded cocommutative bialgebra structure $(\cdot, \Delta_0)$, where $\cdot$ is the canonical multiplication on $\odot^\bu V$ and $\Delta_0$ is the coproduct on $\odot^\bu V$ uniquely determined by the condition that the elements of $V\subset \odot^\bu V$ are primitive. Assume the topological vector space $\odot^\bu(V)[[\hbar]]$, $\hbar$ being a formal parameter of degree zero,
has a continuous bialgebra structure $(\star_\hbar, \Delta_\hbar)$ of the form,
$$
\Ba{rccc}
\star_\hbar: & \odot^\bu V \bigotimes \odot^\bu V & \lon & \odot^\bu V[[\hbar]]
\\
   &   (f(x),g(x)) & \lon & f *_\hbar g = f\cdot g + \sum_{k\geq 1}^\infty \hbar^k P_k(f,g)
\Ea
$$
$$
\Ba{rccc}
\Delta_\hbar: &  \odot^\bu V& \lon & \odot^\bu V \bigotimes \odot^\bu V [[\hbar]]
\\
   &   f(x) & \lon & \Delta_\hbar f = \Delta_0 f + \sum_{k\geq 1}^\infty \hbar^k Q_k(f)
\Ea
$$
where all operators $P_k$ are bidifferential, and $Q_k$ are co-bidifferential (that is, dual to bi-differential operators on polynomial functions, see, e.g., \S 4.2 in \cite{Me1} for their explicit description). It is not hard to see that the bialgebra conditions on these operators modulo terms $O(\hbar^2)$ imply that the first order
deformations $P_1$ and $Q_1$ are uniquely determined by some Lie bialgebra structure $(\vartriangle, [\ ,\ ])$ on $V$. In this case the data $(\star_\hbar, \Delta_\hbar)$ on $\odot^\bu(V)[[\hbar]]$ is called a {\em deformation quantization}\, of that Lie bialgebra structure $(\vartriangle, [\ ,\ ])$ on $V$. Thus Lie bialgebra structures control infinitesimal deformations of the standard Hopf algebra structure on $\odot^\bu V$.

\sip

{\em The Drinfeld quantization problem}: given any Lie bialgebra structure
on a vector space $V$, does there always exist its deformation quantization
$(\star_\hbar, \Delta_\hbar)$? Put another way, given infinitesimal deformations $P_1$ and $Q_1$ of the standard bialgebra on $\odot^\bu V$, can we always find suitable operators $P_2,P_3,\ldots$, $Q_2,Q_3, \ldots$ which make $\odot^\bu(V)[[\hbar]]$ into a (not necessarily commutative cocommutative) bialgebra?

\sip

Surprisingly enough, the problem can not be solve by  a trivial inductive procedure\footnote{The obstructions to such a universal iterative procedure belong to $H^1(\GC_2)$. A ``folklore" conjecture says that $H^1(\GC_2)=0$, but one hast to choose an associator to make every degree 1 cycle in $\GC_2$ into a coboundary; thus an iterative procedure can exist but it can not be trivial.}.
It was proven by P.\ Etingof and D.\ Kazhdan in \cite{EK} that, {\em given a choice of a Drinfeld associator, there does exist a corresponding {\bf\it universal}} (in the sense, for {\em any}\, vector space $V$ and {\em any}\, Lie bialgebra structure on $V$) {\em solution to the Drinfeld quantization problem.} It was proven in \cite{Me1} that any solution of the Drinfeld quantization problem extends to a quasi-isomorphism --- a {\em formality map} --- of completed $\caL ie_\infty$ algebras,
\Beq\label{8: formality lieb F_V}
F_V: \wh{\odot^{\bu\geq 1}} ( V[-1] \oplus  V^*[-1]) \lon
\prod_{m,n\geq 1} \Hom\left(\bigotimes^m (\odot^\bu V), \bigotimes^n (\odot^\bu V)\right)[2-m-n]
\Eeq
where the l.h.s.\ is equipped with the Poisson type Lie bracket
induced by the standard paring $V\ot V^*\rar \K$, and the r.h.s.\ is
the Gerstenhaber-Schack complex of the graded commutative cocommutative Hopf algebra $\odot^\bu V$ equipped with the $\caL ie_\infty$ structure uniquely determined by a choice of a minimal resolution $\cA ss\cB_\infty$ of the properad $\cA ss\cB$ controlling (associative) bialgebras. In the other direction: any quasi-isomorphism  as above gives a solution of the
the Drinfeld quantization problem.

\sip

Since we are interested in {\em universal}\, solutions, it makes sense to search for a reformulation of this problem which does not refer to a particular vector space at all. The theory of operads (see \S {\ref{4: subsec on def of operad}}) is very effective in the study of the homotopy theory of algebraic operations on vector spaces with many inputs but only {\em one} output. It can not be applied to such algebraic structures as Lie bialgebras, but there is a very nice extension of that theory called the theory {\em of props,  properads,
and  their wheeled versions}. We refer to \cite{V} for an excellent introduction into the theory of props and properads.
Roughly speaking, the theory prop(erad)s is based on (connected) decorated oriented graphs with legs.
The deformation theory of {\em morphisms}\, of properads and props was developed in \cite{MV}. There is a properad $\LB$ (see, e.g., \cite{MW-f,MW3,V} and references cited therein) whose representations in a vector space $V$ are precisely Lie bialgebra structures on $V$ and the homotopy theory of such structures is controlled by the minimal resolution $\HoLB$ of $\LB$. There is also a properad $\cA ss\cB$  whose representations in a vector space $W$ are precisely bialgebra structures on $W$ and the homotopy theory of such structures is controlled by a minimal resolution $\cA ss\cB_\infty$ of $\cA ss\cB$
(which exists but is not unique \cite{Ma}).
V.\ Drinfeld deformation quantization problem (in the extended {\em formality}\, version as explained above) can be reformulated as existence of a morphism of dg props \cite{MW-f}
\Beq\label{8: formality map lien props}
F: \cA ss\cB_\infty \lon \caD\wh{\HoLB}
\Eeq
satisfying certain non-triviality conditions. Here $\wh{\HoLB}$ is the genus completion of $\HoLB$ and
$$
\caD: \text{category of dg props} \ \lon \ \text{category of dg props}
$$
a polydifferential (exact) endofunctor which creates out of any prop $\cP$ another prop $\caD \cP$ with the property that there is a one-to-one correspondence
between representations of $\cP$ in a vector space $V$ and representations of $\caD \cP$ in
$\odot^\bu V$ given in terms of polydifferential operators. The Etingof-Kazhdan theorem can be used to show that for any Drinfeld associator
there exists a morphism $F$ of dg props as above. To classify all solutions to the Drinfeld problem, i.e.\ to classify  all possible homotopy non-trivial formality maps as in (\ref{8: formality map lien props}), one has to compute
the cohomology of the deformation complex of any map $F$ as above,
$$
\Def\left(\cA ss\cB_\infty \stackrel{F}{\rar} \caD\wh{\HoLB}\right).
$$
This was done in \cite{MW-f} where it was proven that there is an isomorphism of cohomology groups
$$
H^{\bu+1}\left( \Def\left(\cA ss\cB_\infty \stackrel{F}{\rar} \caD\wh{\HoLB}\right)\right)=H^\bu(\fGC_2')\oplus \K
$$
where $\fGC_2'$ is a version of the full graph complex $\fGC_2$ to which we added by hand a new element $\emptyset$ concentrated in degree zero, ``a graph with no vertices and edges". More precisely,
$$
 \fGC'_{2}:=\widehat{\odot^\bu}\left((\GC_{2}^{\geq 2}\oplus \K )[-2]\right)[2],
$$
the summands $\K$ being generated by $\emptyset$. The formal class $\emptyset$ takes care
for (homotopy non-trivial) rescaling operation of the prop under considerations. The main point of this extension is that the Lie bracket of $\emptyset$  with elements $\Ga$  of $\GC_{2}^{\geq 2}$ is defined as the multiplication of $\Ga$  by twice the number of its loops.

This result implies in particular that
$$
H^{1}\left( \Def\left(\cA ss\cB_\infty \stackrel{F}{\rar} \caD\wh{\HoLB}\right)\right)=\grt_1 \oplus \K,
$$
i.e.\ the Grothendieck-Teichm\"uller group $GRT=GRT_1\rtimes \K^*$ acts faithfully and transitively on the set of solutions of the Drinfeld quantization problem which in turn implies that this set can be identified with the set of Drinfeld associators.

\sip

An explicit transcendental formula for a universal quantization of Lie bialgebras has been
constructed in \cite{MW3}.

\subsection{Universal quantizations of Poisson structures}
Let  $C^\infty(\R^n)$ be the commutative algebra of smooth functions in $\R^n$.
A {\em star product} on $C^\infty(\R^n)$ is a continuous  associative product
on $C^\infty(\R^n)[[\hbar]]$, $\hbar$ being a formal parameter,
$$
\Ba{rccc}
*_\hbar: & C^\infty(\R^n) \times C^\infty(\R^n) & \lon & C^\infty(\R^n)[[\hbar]]\\
   &   (f(x),g(x)) & \lon & f *_\hbar g = fg + \sum_{k\geq 1}^\infty \hbar^k P_k(f,g)
\Ea
$$
where all operators $P_k$ are bidifferential. It is not hard to check that the associativity condition
 on $*_\hbar$ implies that $\pi(f,g):= B_1(f,g) - B_1(g,f)$ is a Poisson structure in $\R^n$; then
$*_\hbar$ is called a {\em deformation quantization}\, of $\pi\in \cT_{poly}(\R^n)$.

\sip

The {\em deformation quantization problem}\, addresses the question: given a Poisson structure $\pi$ on $\R^n$, does there exist a
star product $*_\hbar$ on $C^\infty(\R^n)$ which is a deformation quantization of $\pi$? A spectacular solution of this problem was given by Maxim Kontsevich \cite{Ko2} in the form of a explicit  map between the two sets

$$
\xymatrix{\left\{
\Ba{c}{\mbox{$\Ba{c} \ \mathrm{Poisson} \ \\
 \  \mathrm{structures\ in}\ \R^n\   \Ea$}}\Ea\right\} \ \ \ \ \ \
 \ar[r]^{\mathit{depends\ on}}_{\mathit{associators}} &
\ \ \ \ \ \
\
\left\{\Ba{c}{\mbox{$\Ba{c} \ \mathrm{Star\ products} \ \\
\ *_\hbar\ \mathrm{in}\ C^\infty(\R^n)[[\hbar]]  \Ea$}}
\Ea\right\}
 }
$$
given by transcendental formulae.
In fact a stronger statement was proven --- the {\em formality theorem}\,
which says that for any $n$ there is an explicit $\caL ie_\infty$ quasi-isomorphism of dg Lie algebras,
$$
F_K:\cT_{poly}(\R^n) \lon C^\bu(C^\infty(\R^n), C^\infty(\R^n))
$$
where $\cT_{poly}(\R^n)$ is the Lie algebra of polyvector fields on $\R^n$
equipped with the Schouten-Nijenhuis bracket, and $C^\bu(C^\infty(\R^n), C^\infty(\R^n))$ is the Hochschild complex of the algebra  $C^\infty(\R^n)$, that is, the dg Lie algebra controlling deformations of $C^\infty(\R^n)$ as an associative (but not necessarily commutative) $\R$-algebra. Moreover, the formality theorem holds true for any manifold $M$, not necessarily for $\R^n$ \cite{Ko2}.
Later D.\ Tamarkin has proven  the existence theorem for deformation quantizations which exhibited a key role of Drinfeld's associators \cite{Ta2}, and V.\ Dolgushev \cite{Do} has proven that the set of homotopy classes of universal formality maps can be identified with the set of Drinfeld's associators so that $GRT$ acts faithfully and transitively of homotopy classes of formality maps.

\sip

This classification result can be restated in a very short and compact form without any reference to any particular manifold as follows \cite{AsMe}. Let $c\cA ss_\infty$ be the dg free operad of strongly homotopy curved associative algebras,
$\HoLB_{1,0}^\circlearrowright$ the wheeled closure \cite{MMS} of the  prop $\HoLB_{1,0}$ of $(1,0)$ Lie bialgebras and
$$
\f: \text{category of dg (wheeled) props} \ \lon \ \text{category of dg operads}
$$
the polydifferential functor introduced in \cite{MW2} (which is an operadic part of the functor $\caD$ mentioned in the previous subsection). Then the Kontsevich formality map $F_K$ can be understood as a morphism of dg operads
$$
F_K: c\cA ss_\infty \lon \f(\wh{\HoLB}_{1,0}^\circlearrowright)
$$
so that it makes sense to study a deformation complex of the Kontsevich map
$$
\Def\left(c\cA ss_\infty \stackrel{F_K}{\lon} \f(\wh{\HoLB}_{1,0}^\circlearrowright)\right)
$$
It was proven in \cite{AsMe} that there exists isomorphism of cohomology groups
$$
H^{\bu}(\fGC_2')= H^{\bu+1}\left(\Def\left(c\cA ss_\infty \stackrel{F_K}{\lon} \f(\HoLB_{1,0}^\circlearrowright)\right)\right)
$$
implying
$$
\grt=H^1(\fGC_2')=H^1\left(\Def\left(c\cA ss_\infty \stackrel{F_K}{\lon} \f(\HoLB_{1,0}^\circlearrowright)\right)\right)
$$
Hence the Grothendieck-Teichm\"uller Lie algebra $\grt_1$ controls {\em all}\,
homotopy non-trivial infinitesimal deformations of the Kontsevich map $F_K$. Each such infinitesimal deformation can be exponentiated to a genuine deformation of $F_K$ implying the fundamental classifying role of $GRT_1$
in the theory of universal deformation quantizations of Poisson structures. A concrete algorithm for computing weights in the M.\ Kontsevich formula for $F_K$  was developed in \cite{BPP} where it was proven that all such weights are rational linear combinations of multiple zeta values.

\subsubsection{\bf Action of $GRT$ on Lie bialgebras}
One can consider a degree shifted version of the notion of Lie bialgebra. For any integers $c,d\in \Z$ one defines \cite{MW2} a {\em Lie $(c,d)$-bialgebra}\, as a graded vector space equipped with linear maps
$$
[\ , \ ]: \odot^2 (V[d])
\rightarrow V[d+1] \ \ \ , \ \ \vartriangle: V[c]\rightarrow \odot^2(V[c])[1-2c] ,
$$
 making $V$ into a degree shifted
Lie algebra and degree shifted Lie coalgebra satisfying an analogue of the Drinfeld compatibility condition. There exist a prop(erad) $\LB_{c,d}$ which governs these structures and its minimal resolution $\HoLB_{c,d}$
which governs their homotopy theory \cite{MW2}. The case $c=1$, $d=1$ gives
us ordinary Lie bialgebras discussed above, while the case $c=1$, $d=0$
is important in the theory of Poisson structures on manifolds as there is a
one-to-one correspondence \cite{Me-w} between representation of the pro(erad) $\HoLB_{1,0}$ on a graded vector space $V$ and formal graded Poisson structures on $V$ viewed as a formal manifold which vanish at the distinguished point $0\in V$.

 The deformation theory \cite{MV} of the identity map $\Id: {\wHoLBcd}\rar \wHoLBcd$ gives us a clue to the symmetry group of the genus completion ${\wHoLBcd}$ of the dg prop $\HoLB_{c,d}$ (and hence of $\wh{\LB_{c,d}}$). This problem was settled in \cite{MW3} where it was proven that there is a quasi-isomorphism (up to one rescaling class) of complexes
$$
\GCor_{c+d+1} \lon
\Def\left(\wHoLBcd\stackrel{\Id}{\rar} \wHoLBcd\right)[1]
$$
implying
$$
H^\bu(\GC_{c+d}^{\geq 2})\oplus \K =  H^\bu(\GCor_{c+d+1})\oplus \K= H^{\bu+1}\left(\wHoLBcd\stackrel{\Id}{\rar} \wHoLBcd\right)
$$
so that
$$
\grt_1\oplus \K= H^{1}(\wh{\HoLB}\stackrel{\Id}{\rar} \wh{\HoLB})
$$
Hence
the Grothendieck-Teichm\"uller group $GRT=GRT_1\rtimes \K^*$ is essentially the automorphism group of the  completed prop(erad) $\wh{\LB}$.

%

\sip

There is  a further generalization of the notion of Lie $(c,d)$-bialgebra to the {\em multioriented}\, case, and $GRT$ acts on such structures for any $c,d\in \N$ with  $c+d\geq 3$ \cite{A, Me5}.

\subsection{Solutions of the Kashiwara-Vergne conjecture} Let $\wh{\fass}_2$
and $\wh{\mathsf{lie}}_2$ be the completed free associative and, respectively, Lie algebra generated by  formal variables $x$ and $y$ (see \S {\ref{2: ass_n and Lie_n}}). Consider the quotient space
$$
\wh{Cyc}_2:=\wh{\fass}_2/\langle AB - BA \mid \forall\, A,B\in \wh{\fass}_2\rangle\equiv\wh{\fass}_2/[\wh{\fass}_2,\wh{\fass}_2]
$$
the (completed) vector space spanned by cyclic words in two letters. There is a canonical projection $tr: \wh{\fass}_2\rar \wh{Cyc}_2$. Note that every element $A\in \wh{\fass}_2$ has a unique decomposition,
$$
A=A_0 + \p_x(A)x + \p_y(A) y
$$
for some $A_0\in \K$, $\p_x(A), \p_y(A)\in \wh{\fass}_2$.

\sip

Recall the Bernoulli power series series
$$
\frac{x}{1 + e^x}
= 1+ \frac{x}{2}
+ \sum_{n\geq 1}\frac{b_{2n}}{2n!}x^{2n} =: 1+
\frac{x}{2}
+\mathfrak{b}(x)
$$
and the Baker-Cambell-Hausdorf power series $\mathfrak{bch}(x,y)$ defined in (\ref{2: ch(c,y) series}). A triple $(A,B,g)$ consisting of two Lie series  $A,B\in \wh{\mathsf{lie}}_2$
and a formal power series $g(x)=\sum_{n\geq 2} g_nx^n$
is called a {\em solution of the (generalized) Kashiwara-Vergne problem}\, if they satisfy the equations
$$
x+y - \mathfrak{bch}(x,y)=(1-e^{-\text{ad}_x})A + (e^{\text{ad}_y} -1)B\ \ \ \ \text{in}\ \wh{\mathsf{lie}}_2
$$
and
$$
tr\left(\p_x(A)x + \p_y(B)y\right)=\frac{1}{2} tr\left( g(x) - g(\mathfrak{bch}(x,y))+ g(y)\right).
$$
Given a Lie group $G$ with the Lie algebra $\fg$, if a solution of the $KV$ problem exists, then the $G$-invariant
harmonic analysis on $\fg$ is related with $G$-invariant harmonic analysis on the group $G$ itself by the standard exponential
map. Moreover, if a solution exists, then $g_{even}:=\sum_{n\geq 1} g_{2n}x^{2n}$ must be equal to the Bernoulli series $\mathfrak{b}(x)$.

\sip

Existence of solution of the $KV$ problem was established in \cite{AM}. An alternative solution of the $KV$ problem which classified all such solutions in terms of so called {\em $KV$ associators}\, was given in \cite{AT}. Any Drinfeld associator is a $KV$ associator, and there is a conjecture which says that this association is one-to-one.

\subsection{Formality theorem in the Goldman-Turaev theory}
Let  $\Sigma$ be a Riemann surface of genus $g$ with $N+1$ boundary components, $N\in \N$. The group algebra ${\K}\langle\pi_1(\Sigma)\rangle$
of the   fundamental group $\pi_1(\Sigma)$  is canonically filtered by powers
of the augmentation ideal (see \S {\ref{2: Quillen's construction}}) and hence
admits a canonical completion $\wh{\pi_1(\Sigma)}$. The associated
(completed) vector space spanned by conjugacy classes in $\pi_1(\Sigma)$,
$$
\widehat{\fg[\Sigma]}:= \frac{\widehat{\K\langle\pi_1(\Sigma)\rangle}}{[ \widehat{\K\langle\pi_1(\Sigma)}, \widehat{\K\langle\pi_1(\Sigma)\rangle}]}
$$
 that is,  by free homotopy classes of loops in $\Sigma$, has a canonical {\em Goldman-Turaev Lie bialgebra structure}. It is filtered, and the associated
 graded vector space
 $$
 \mathrm{gr}\widehat{\fg[\Sigma]} := \frac{\widehat{\ot^\bu} H_1(\Sigma)}{[ \widehat{\ot^\bu} H_1(\Sigma), \widehat{\ot^\bu} H_1(\Sigma)]},
$$
where $H_1(\Sigma)$ is the first homology group of $\Sigma$ over $\K$, has an induced Lie bialgebra structure which admits a rather simple combinatorial description.
The {\em formality theorem}  \cite{AKKN1,AKKN2,Mas} (see also references cited therein) says that for any solution of the $KV$ problem, in particular, for any Drinfeld associator
there is an associated isomorphism
$$
\Theta:\  \widehat{\fg[\Sigma]}  \lon \ \
\mathrm{gr}\widehat{\fg[\Sigma]}
$$
of Lie bialgebras. In the case $g=0$ and $\K=\C$ a nice explicit formula
for such an isomorphism was constructed in \cite{AN} with the help of the
Knizhnik-Zamolodchikov connection. In particular, the group $GRT_1$
acts as automorphisms on $\mathrm{gr}\widehat{\fg[\Sigma]}$ for any Riemann surface $\Sigma$.

\subsection{Cohomology of moduli spaces of algebraic curves and $\grt_1$}
Let $\cM_g$ be the moduli space of Riemann surfaces of genus $g$. The authors of \cite{CGP} constructed a remarkable monomorphism of cohomology groups
$$
H^\bu(\GC_2) \lon \prod_{g}H^{\bu+2g}_c(\cM_g;\K)= \prod_g (H^{4g-6-\bu}(\cM_g;\K))^*
$$
implying an injection of the Grothendieck-Teichm\" uller Lie algebra
$$
\grt_1 \lon  \prod_{g}H^{2g}_c(\cM_g;\K)
$$
and hence the following estimation \cite{CGP} on the dimension of cohomology groups:
$$
\dim H^{4g-6}(\cM_g,\Q)> \beta^g + \text{constant},
$$
for any $\beta<\beta_0$, where $\beta_0\approx 1.3247\ldots$ is the real root
of $t^3 - t -1=0$.

\bip

\bip

\def\cprime{$'$}


\begin{thebibliography}{10}

\bibitem[AKKN1]{AKKN1} A.\ Alekseev, N.\ Kawazumi, Y.\ Kuno and F.\ Naef,
{\em The Goldman-Turaev Lie bialgebra in
genus zero and the Kashiwara-Vergne problem}, Adv.\ Math. {\bf 326} (2018) 1-53

\bibitem[AKKN2]{AKKN2} A.\ Alekseev, N.\ Kawazumi, Y.\ Kuno and F.\ Naef,
{\em The Goldman-Turaev Lie bialgebra and the Kashiwara-Vergne problem in higher genera}, arxiv arXiv:1703.0581.

\bibitem[AN]{AN} A.\ Alekseev and F.\ Naef,
{\em Goldman-Turaev formality from the Knizhnik-Zamolodchikov connection},
C.\ R.\ Math. Acad. Sci. Paris {\bf 355} (2017), no. 11, 1138-1147.







\bibitem[AlMe]{AM} A.\ Alekseev and E.\ Meinrenken. {\em On the Kashiwara-Vergne conjecture},
Invent. Math. {\bf 164} (2006), {615--634}.


\bibitem[AT1]{AT2} A.\ Alekseev and C.\ Torossian,
{\em
Kontsevich deformation quantization and flat connections},
Comm.\ Math.\ Phys.\  {\bf 300} (2010) 47-64.


\bibitem[AT2]{AT}
A. Alekseev, C. Torossian,
{\em The Kashiwara-Vergne conjecture and Drinfeld's associators},
Ann.\ of Math.\ (2) {\bf 175} (2012), no.\ 2, 415-463.


\bibitem[An]{A} A.\ Andersen, {\em Haired and multioriented graph complexes
with applications to algebra and geometry},  PhD thesis, University of Luxembourg, 2020.

 \bibitem[AnMe]{AsMe}
A.\ Andersson and S. A.\ Merkulov, {\em From deformation theory of wheeled props to classification of Kontsevich formality maps},  preprint arXiv:1911.09089 (2019)

\bibitem[BPP]{BPP} P.\ Banks, E.\ Panzer and B.\ Pym, {\em Multiple zeta values in deformation quantization}. Invent.\ Math. (2020). https://doi.org/10.1007/s00222-020-00970-x



\bibitem[B-N]{BN}  D.\ Bar-Natan, {\em On associators and the Grothendieck-Teichmuller group. I}. Selecta Math. (N.S.)
{\bf 4} (1998), no. 2, 183-212.

\bibitem[B]{B} F.\ Brown, {\em Mixed Tate Motives over $Spec(\Z)$}, Annals of Math.\ {\bf 175} (2012), no. 2, 949-976.

\bibitem[Br]{Br} A.\ Brochier, {\em Introduction to Drinfeld's associator theory}, preprint.

\bibitem[CMW]{CMW} R.\ Campos, S.\ Merkulov and T.\ Willwacher {\em The Frobenius properad is Koszul},  Duke Math.\ J. {\bf 165}, No.1 (2016), 2921-2989.

\bibitem[CGP1]{CGP} M.\ Chan, S.\ Galatius and S.\ Payne,
{\em Tropical curves, graph complexes, and top weight
cohomology of $\cM_g$}, preprint arXiv:1805.10186 (2018)

\bibitem[CGP2]{CGP2} M.\ Chan, S.\ Galatius and S.\ Payne,
{\em Topology of moduli spaces of tropical curves with marked points},
 arXiv:1903.07187 (2019)

\bibitem[CW]{CW} F.R.\ Cohen and J.\ Wu, {\em On braid groups and homotopy groups}, Geometry \& Topology Monographs {\bf 13}
 (2008) 169-193.

\bibitem[Do]{Do} V.\ Dolgushev, {\em Stable Formality Quasi-isomorphisms for Hochschild Cochains I},  (2011), available at arXiv:1109.6031.

\bibitem[Dr1]{D1}
V.\ Drinfeld,
{\em Hamiltonian structures on Lie groups, Lie bialgebras and the geometric
meaning of the classical Yang-Baxter equations}, Soviet Math. Dokl. {\bf 27} (1983) 68--71.

\bibitem[Dr2]{Dr}
V. Drinfeld, {\em On quasitriangular quasi-Hopf algebras and a group closely connected
with $Gal(\bar{Q}/Q)$}, Leningrad Math. J. {\bf 2}, No.\ 4 (1991),  829-860.

\bibitem[Ei]{Ei} D.\ Eisenbud,
Commutative algebra, with a view towards algebraic geometry. Springer, 2004

\bibitem[EE]{EE}
B.\ Enriquez and P.\ Etingof,
\newblock {\em On the invertability of quantization functors},
Advances in Math.\ {\bf 197} (2005), 430-479


\bibitem[EK]{EK} P.\ Etingof and D.\ Kazhdan.
\newblock {\em Quantization of Lie bialgebras, I}.
\newblock{ Selecta Math. (N.S.)} {\bf 2} (1996), 1-41.



\bibitem[ES]{ES} P.\ Etingof and O.\ Schiffmann,
Lectures on Quantum Groups, International Press, 2002.

\bibitem[FV]{FV} E.\ Fadell and J.\ Van Buskirk, {\em  The braid groups of $E^2$ and $S^2$}. Duke
Math.\ J., {\bf 29} (1962) 243-257.


\bibitem[FHT]{FHT} Y.\ Felix, S.\ Halperin and J.-C.\ Thomas, {\em Rational homotopy theory}, Springer, 2001.


\bibitem[Fr]{Fr}  B.\ Fresse, Homotopy of Operads and Grothendieck-Teichm\"uller Groups,
 (volume 217 of the series Mathematical Surveys and Monographs) AMS, 2017

\bibitem[Fu1]{Fu} H.\ Furusho, {\em Pentagon and hexagon equations}, Annals of Mathematics,
{\bf 171} (2010), No. 1, 545-556.

\bibitem[Fu2]{Fu2} H.\ Furusho, {\em On the coefficients of
the Alekseev-Torossian associator}, Journal of Algebra {\bf 506} (2018) 364-378.

\bibitem[G]{Go}  A.\ B.\ Goncharov, {\em Hodge correlators}, arXiv:0803.0297 (2008)

\bibitem[Ha]{Ha} R.\ Hain, {\em Completions of mapping class groups and the cycle $C-C^{-}$}. In :
Mapping class groups and moduli spaces of Riemann surfaces (G\"ottingen,
1991/Seattle, WA, 1991), Contemp. Math. {\bf 150} (1993) 75-105.


\bibitem[HLTV]{HLTV} R.\ M. Hardt, P.\ Lambrechts, V. Turchin, and I.\ Voli\'{c},
{\em Real homotopy theory of semi-algebraic sets}, preprint arXiv:0806.0476 (2008).


\bibitem[KM]{KM} M.\ Kapranov and Yu.I.\ Manin, {\em Modules and Morita theorem for operads}. Amer. J. Math. {\bf 123} (2001), no. 5, 811-838.

\bibitem[Ko]{Koh} M.\ Kohno, {\em S\'erie de Poincar\'e-Koszul associee aux groupes de tresses pures Deformation quantization}, Invent. Math. {\bf 82} (1985), no. 1, 57-75.



 \bibitem[K1]{Ko1} M. Kontsevich, {\em Formality Conjecture}, D. Sternheimer et al. (eds.),
Deformation Theory and Symplectic
Geometry, Kluwer 1997, 139-156.

 \bibitem[K2]{Ko2} M.\ Kontsevich, {\em Deformation quantization
 of Poisson manifolds}, Lett.\ Math.\ Phys. {\bf 66} (2003), 157-216.

\bibitem[K3]{Ko3} M.\ Kontsevich, {\em
Operads and motives in deformation quantization}, Lett.\ Math.\ Phys.
{\bf 48}(1) (1999), 35-72.

\bibitem[KWZ1]{KWZ1}  A.\ Khoroshkin, T.\ Willwacher and
M.\ \v Zivkovi\'c, {\em Differentials on graph complexes}, preprint arXiv:1512.04710 (2015)


\bibitem[KWZ2]{KWZ2}  A.\ Khoroshkin, T.\ Willwacher and
M.\ \v Zivkovi\'c, {\em Differentials on graph complexes II - hairy
graphs}, preprint arXiv:1508.01281 (2016).


\bibitem[K]{Kn} K.\ P.\ Knudsen, {\em Relative completions and the cohomology
of linear groups over local rings}, J.\ London Math. Soc. {\bf 65} (2002) 183-203.

\bibitem[L]{L} P.\ Lochak, {\em Fragments of nonlinear Grothendieck-Teichm\"uller theory}. Woods Hole Mathematics,
225-262, Series on Knots and Everything 34, World Scientfic, 2004.

\bibitem[LS]{LS} P.\ Lochak and L.\ Schneps, {\em  Open problems in Grothendieck-Teichm\"uller theory}, preprint 2011.

\bibitem[LV]{LV} J.-L.\ Loday  and B.\ Vallette, Algebraic operads. Springer 2012.





\bibitem[Ma]{Ma} M.\ Markl,
\newblock {\em A resolution (minimal model) of the prop for bialgebras},
J.\ Pure Appl.\ Algebra {\bf 205} (2006), no. 2, 341-374.






\bibitem[MMS]{MMS}  M.\ Markl, S.\ Merkulov and S.\ Shadrin, {\em Wheeled props and the master
equation}, preprint math.AG/0610683, J.\ Pure and Appl.\ Algebra {\bf 213} (2009), 496-535.

\bibitem[MSS]{MSS}
M.~Markl, S.~Shnider, and J.~D. Stasheff.
\newblock {Operads in Algebra, Topology and Physics}, volume~96 of {
  Mathematical Surveys and Monographs}.
\newblock American Mathematical Society, Providence, Rhode Island, 2002.

\bibitem[Mas]{Mas}
G.\  Massuyeau, {\it Formal descriptions of Turaev's loop operations}, Quantum Topol. {\bf 9} (2018) 39-117.

\bibitem[May]{May}  J.P. May. The Geometry of Iterated Loop Spaces, volume 271 of Lecture
Notes in Mathematics. Springer-Verlag, New York, 1972.

\bibitem[Me1]{Me-w} S.A.\ Merkulov,
 {\em Graph complexes with
  loops and wheels}. In:
  ``Algebra, Arithmetic
and Geometry - Manin Festschrift" (eds. Yu.\ Tschinkel and Yu.\ Zarhin),
Progress in Mathematics, Birkha\"user (2010) 311-354.



 \bibitem[Me2]{Me4} S.A.\ Merkulov, {\em Operads, configuration spaces and quantization}.
 In: ``Proceedings of Poisson 2010, Rio de Janeiro", Bull.\ Braz.\ Math.\ Soc., New Series {\bf 42}(4) (2011), 1-99.

\bibitem[Me3]{Me1} S.A.\ Merkulov, {\em Formality theorem for
 quantizations of Lie bialgebras},  Lett.\ Math.\ Phys. {\bf 106} (2016) 169-195

 \bibitem[Me4]{Me5} S.A.\ Merkulov, {\em Multi-oriented props and homotopy algebras with branes},   Lett.\ Math.\ Phys.  {\bf 110} (2020)  1425-1475.






\bibitem[MW1]{MW0} S.\ Merkulov and T.\ Willwacher, {\em Grothendieck-Teichm\"uller and Batalin-Vilkovisky}, Lett.\  Math.\ Phys.  {\bf 104} (2014) No.\ 5,  625-634.

 \bibitem[MW2]{MW} S.\ Merkulov and T.\ Willwacher, {\em Props of ribbon graphs, involutive Lie bialgebras and moduli spaces of curves}, preprint  arXiv:1511.07808 (2015) 51pp.

\bibitem[MW3]{MW-f} S.A. Merkulov and T.\ Willwacher, {\em Classification of universal formality maps
for quantizations of Lie bialgebras}, Compositio Mathematica, {\bf 156} (2020) 2111-2148.

\bibitem[MW4]{MW2} S.\ Merkulov and T.\ Willwacher, {\em Deformation theory of  Lie bialgebra properads},   In: Geometry and Physics: A Festschrift in honour of Nigel Hitchin, Oxford University Press 2018, pp. 219-248

\bibitem[MW5]{MW3} S.\ Merkulov and T.\ Willwacher, {\em An explicit two step quantization of Poisson structures and Lie bialgebras},  Comm.\ Math.\ Phys. {\bf 364} (2018),  505-578




\bibitem[MV]{MV}  S.A.\ Merkulov and  B.\ Vallette,
{\em Deformation theory of representations of prop(erad)s I \& II},
{ Journal f\"ur die reine und angewandte Mathematik (Qrelle)  634}, 51-106,
 \& {\bf 636}, 123-174 (2009)


\bibitem[RW]{RW}
C.\ Rossi and T.\ Willwacher, {\em P.\ Etingof's conjecture about Drinfeld associators}, preprint
arXiv:1404.2047 (2014)



\bibitem[Se]{Se} P.\ Severa, {\em
Formality of the chain operad of framed little disks},  Lett.\ Math.\ Physics {\bf 93} (2010)  29-35.




  \bibitem[SW]{SW} P.\ Severa and T.\ Willwacher, {\em
  Equivalence of formalities of the little discs operad}, Duke
Math.\ J.\ {\bf 160} (2011), no. 1, 175-206.



\bibitem[Sh]{Sh} B.\ Shoikhet, {\em  An $L_\infty$ algebra structure on polyvector fields
}, preprint arXiv:0805.3363, (2008).




\bibitem[T1]{Ta1} D.E.\ Tamarkin, {\em Formality of Chain Operad of Small Squares
},  Lett.\ Math.\ Phys. {\bf 66} (2003), no. 1-2, 65-72.


\bibitem[T2]{Ta2} D.E.\ Tamarkin, {\em Another proof of M. Kontsevich formality
theorem},  math.QA/9803025, Lett.\ Math.\ Phys. {\bf 66}
(2003) 65-72.



\bibitem[TW]{TW}  V.\ Turchin and  T.\ Willwacher, {\em
Commutative hairy graphs and representations of $Out(F_r)$},
J. Toplogy {\bf 10}, Issue 2 (2017) 386-411

\bibitem[V]{V}
B.\ Vallette,  {\em A Koszul duality for
props}, Trans.\ Amer. Math. Soc., {\bf 359} (2007), 4865--4943.


 \bibitem[W1]{Wi1} T.\ Willwacher, {\em M.\ Kontsevich's graph complex and the Grothendieck-Teichmueller Lie algebra},
Invent. Math. {\bf 200} (2015) 671-760.

\bibitem[W2]{Wi3} T.\ Willwacher, {\em
Stable cohomology of polyvector fields}, Math. Res. Letters. {\bf 21} (2014) 1501-1530

\bibitem[W3]{Wi2} T.\ Willwacher, {\em Oriented graph complexes},   Comm. Math. Phys. {\bf 334} (2015), no. 3, 1649-1666.

\bibitem[W4]{Wi6} T.\ Willwacher,  {\em Deformation quantization and the Gerstenhaber structure on the homology of knot spaces},
    preprint arXiv:1506.07078 (2015)

\bibitem[W5]{Wi4} T.\ Willwacher, {\em The Homotopy Braces Formality Morphism}, Duke Math.\ J.\ {\bf 165} (2016), no.
10, 1815-1964.

\bibitem[W6]{Wi5} T.\ Willwacher, {\em
Pre-Lie pairs and triviality of the Lie bracket on the twisted hairy graph complexes}, preprint arXiv:1702.04504 (2017)

\bibitem[Z1]{Z1} M.\ \v Zivkovi\' c, {\em Multi-oriented graph complexes and quasi-isomorphisms between them I: oriented graphs}, High.\ Struct.\ 4 (2020), no. 1, 266-283.

\bibitem[Z2]{Z2} M.\ \v Zivkovi\' c, {\em Multi-oriented graph complexes and quasi-isomorphisms between them II: sourced graphs}, preprint
 arXiv:1712.01203 (2017)

 \end{thebibliography}
\end{document}